\documentclass[prodmode,acmtoms]{acmsmall} 

\usepackage[ruled]{algorithm2e}
\usepackage[squaren]{SIunits}
\usepackage{amsmath}
\usepackage{amssymb}
\usepackage{multirow}
\usepackage{url}
\usepackage{epsfig}
\usepackage{color}

\newcommand{\R}{\mathbb{R}}
\newcommand{\C}{\mathbb{C}}

\newcommand{\Z}{\mathbb{Z}}
\newcommand{\N}{\mathbb{N}}

\newcommand{\ip}[2]{\langle #1,#2 \rangle}
\newcommand{\abs}[1]{\left| #1 \right|}

\newcommand{\cH}{\mathcal{H}}
\newcommand{\bitem}{\begin{itemize}}
\newcommand{\eitem}{\end{itemize}}
\newcommand{\benum}{\begin{enumerate}}
\newcommand{\eenum}{\end{enumerate}}
\newcommand{\beq}{\begin{equation}}
\newcommand{\eeq}{\end{equation}}
\newcommand{\absip}[2]{| \langle#1,#2\rangle |}
\def\cSH{{\mathcal{S} \mathcal{H}}}
\DeclareMathOperator*{\supp}{supp}

\usepackage{listings}
\usepackage{microtype}
\definecolor{bgcolor}{rgb}{0,0,0}
\definecolor{string}{rgb}{1,0,0}
\definecolor{darkblue}{rgb}{0,0,.75}
\definecolor{darkred}{rgb}{.75,0,0}
\definecolor{darkgreen}{rgb}{0,.75,0}
\definecolor{grey}{rgb}{0.9,0.9,0.9}

\lstnewenvironment{MatlabCode}[1][]
{\microtypesetup{activate=false}      
\lstset{
   language= Matlab,
   basicstyle=\ttfamily,                     
   basicstyle=\scriptsize,
   keywordstyle=\color{darkblue},
   commentstyle=\color{darkgreen},
   stringstyle=\color{string},
   backgroundcolor=\color{grey},      
   showstringspaces=false,                  
   captionpos=b,                                 
   numbers   =   none,                              
   xleftmargin=.04\textwidth,
   #1}
}
{}

\renewcommand{\algorithmcfname}{ALGORITHM}
\SetAlFnt{\small}
\SetAlCapFnt{\small}
\SetAlCapNameFnt{\small}
\SetAlCapHSkip{0pt}
\IncMargin{-\parindent}

\begin{document}

\markboth{G. Kutyniok et al.}{ShearLab 3D: Faithful Digital Shearlet Transforms based on Compactly Supported Shearlets}

\title{ShearLab 3D: Faithful Digital Shearlet Transforms based on Compactly Supported Shearlets}
\author{
Gitta Kutyniok
\affil{Technische Universit\"at Berlin}
Wang-Q Lim
\affil{Technische Universit\"at Berlin}
Rafael Reisenhofer
\affil{Technische Universit\"at Berlin}}

\begin{abstract}
Wavelets and their associated transforms are highly efficient when approximating and analyzing
one-dimensional signals. However, multivariate signals such as images or videos typically exhibit
curvilinear singularities, which wavelets are provably deficient of sparsely approximating
and also of analyzing in the sense of, for instance, detecting their direction. Shearlets
are a directional representation system extending the wavelet framework, which overcomes those
deficiencies. Similar to wavelets, shearlets
allow a faithful implementation and fast associated transforms. In this paper, we will
introduce a comprehensive carefully documented software package coined ShearLab 3D (\url{www.ShearLab.org})
and discuss its algorithmic details. This package provides MATLAB code for a novel faithful
algorithmic realization of the 2D and 3D shearlet
transform (and their inverses) associated with compactly supported universal shearlet systems
incorporating the option of using CUDA. We will present extensive numerical experiments in 2D and
3D concerning denoising, inpainting, and feature extraction, comparing the performance of ShearLab
3D with similar transform-based algorithms such as curvelets, contourlets, or surfacelets. In
the spirit of reproducible reseaerch, all scripts are accessible on \url{www.ShearLab.org}.
\end{abstract}





\acmformat{Gitta Kutyniok, Wang-Q Lim, and Rafael Reisenhofer, 2014.
ShearLab 3D: Faithful Digital Shearlet Transform based on Compactly Supported Shearlets.}

\begin{bottomstuff}
This work is supported in part by the Einstein Foundation Berlin, by Deutsche Forschungsgemeinschaft (DFG) Grant SPP-1324
KU 1446/13 and DFG Grant KU 1446/14, by the DFG Collaborative Research Center TRR 109 ``Discretization in Geometry and
Dynamics'', and by the DFG Research Center {\sc Matheon} ``Mathematics for key technologies'' in Berlin.

Author's addresses: Gitta Kutyniok, Wang-Q Lim, {and} Rafael Reisenhofer, Department of Mathematics, Technische Universit\"at Berlin, 
10623 Berlin, Germany.
\end{bottomstuff}

\maketitle

\section{Introduction}

Wavelets have had a tremendous success in both theoretical and practical applications such as, for instance,
in optimal schemes for solving elliptic PDEs or in the compression standard JPEG2000. A wavelet system
is based on one or a few generating functions to which isotropic scaling operators and translation
operators are applied to. One main advantage of wavelets is their ability to deliver highly sparse
approximations of 1D signals exhibiting singularities, which makes them a powerful maximally flexible
tool in being applicable to a variety of problems such as denoising or detection of singularities.
But similarly important for applications is the fact that wavelets admit a faithful digitalization of the
continuum domain systems with efficient algorithms for the associated transform computing the respective
wavelet coefficients (cf. \cite{Dau92,Mal08}).

However, each multivariate situation starting with the 2D situation differs significantly from the
1D situation, since now not only (0-dimensional) point singularities, but in addition typically also
(1-dimensional) curvilinear singularities appear; one can think of edges in images or shock fronts
in transport dominated equations. Unfortunately, wavelets are deficient to adequately handle such
data, since they are themselves isotropic -- in the sense of not directional based -- due to their
isotropic scaling matrix. Thus, it was proven in \cite{CD2004} that wavelets do not provide optimally
sparse approximations of 2D functions governed by curvilinear singularities in the sense of the decay rate of the
$L^2$-error of best $N$-term approximation. This causes problems for any application requiring
sparse expansions such as any imaging methodology based on compressed sensing (cf. \cite{DDEK12}).
Moreover, being associated with just a scaling and a translation parameter, wavelets are, for instance,
also not capable of detecting the orientation of edge-like structures.

\subsection{Geometric Multiscale Analysis}

These problems were the reason that within applied harmonic analysis the research area of geometric multiscale
analysis arose, whose main goal consists in developing representation systems which efficiently capture
and sparsely approximate the geometry of objects such as curvilinear singularities of 2D functions. One
approach pursued was to introduce a 2D model situation coined {\em cartoon-like functions} consisting of
functions compactly supported on the unit square while being $C^2$ except for a closed $C^2$ discontinuity curve.
A representation system was referred to as `optimally sparsely approximating cartoon-like functions' provided
that it provided the optimally achievable decay rate of the $L^2$-error of best $N$-term approximation.

A first breakthrough could be reported in 2004 when Cand\`{e}s and Donoho introduced the system of
{\em curvelets}, which could be proven to optimally sparsely approximating cartoon-like functions while
forming Parseval frames \cite{CD2004}.
Moreover, this system showed a performance superior to wavelets in a variety of applications (see, for instance,
\cite{HWHM08,SMF10}). However, curvelets suffer from the fact that in addition to a parabolic scaling operator
and translation operator, the rotation operator utilized as a means to change the orientation can not be
faithfully digitalized; and the implementation could therefore not be made consistent with the continuum domain theory
\cite{CDDY2006}. This led to the introduction of {\em contourlets} by Do and Vetterli \cite{DV2005}, which can
be seen as a filterbank approach to the curvelet transform. However, both of these well-known approaches do
not exhibit the same advantages wavelets have, namely a unified treatment of the continuum and digital situation,
a simple structure such as certain operators applied to very few generating functions, and a theory for compactly
supported systems to guarantee high spatial localization.

\subsection{Shearlets and Beyond}

Shearlets were introduced in 2006 to provide a framework which achieves these goals. These systems are
indeed generated by very few functions to which parabolic scaling and translation operators as well
as shearing operators to change the orientation are applied \cite{KL2012}. The utilization of
shearing operators ensured --~due to consistency with the digital lattice~-- that the continuum and
digital realm was treated uniformly in the sense of the continuum theory allowing a faithful
implementation (cf. for band-limited generators \cite{KSZ12}). A theory for compactly supported shearlet frames
is available \cite{KKL2012}, showing that although presumably Parseval frames can not be derived,
still the frame bounds are within a numerically stable range. Moreover, compactly supported shearlets
can be shown to optimally sparsely approximating cartoon-like functions \cite{KL11}. Also, to date a
3D theory is available \cite{KLL2012}.

Very recently, two extensions of shearlet theory were explored. One first extension is the theory of
{\em $\alpha$-molecules} \cite{GKKS2013}. This approach extends the theory of parabolic molecules
introduced in \cite{GK14}, which provides a framework for systems based on parabolic scaling such
as curvelets and shearlets to analyze their sparse approximation properties. $\alpha$-molecules
are a parameter-based framework including systems based on different types of scaling such as,
in particular, wavelets and shearlets, with the parameter $\alpha$ measuring the degree of anisotropy.
As a subfamily of this general framework so-called {\em $\alpha$-shearlets} (also sometimes called
{\em hybrid shearlets}) were studied in \cite{KLL2012} (cf. also \cite{Kei12}), which can be
regarded as a parametrized family ranging from wavelets ($\alpha = 2$) to shearlets ($\alpha = 1$).
Again, the frame bounds can be controlled and optimal sparse approximation properties -- now for a
parametrized model situation -- are proven, also for compactly supported systems.

A further extension are {\em universal shearlets}, which allow a different type of scaling (parabolic,
etc.) at each scaling level of $\alpha$-shearlets by setting $\alpha = (\alpha_j)_j$ with $j$ being
the scale, thereby achieving maximal flexibility \cite{GK15}. This approach has so far been only
analyzed for band-limited generators, deriving properties on the sequence $(\alpha_j)_j$ such that
the resulting system forms a Parseval frame.

\subsection{Contributions}

Several implementations of shearlet transforms are available to date, and we refer to Subsections \ref{subsec:previousimplements}
and \ref{subsec:comparison} for more details. Most of those focus on the (2D) band-limited case. In this situation, a Parseval frame
can be achieved providing immediate numerical stability and a straightforward inverse transform. However, from an application
point of view, those approaches typically suffer from high complexity, various artifacts, and insufficient spatial localization.
The faithful algorithmic realization suggested in \cite{Lim2010} by one of the authors was the first to focus on compactly supported
shearlets, achieving also low complexity by utilizing separable shearlet generators. Interestingly, this approach could then be
improved in \cite{Lim2013} by utilizing non-separable compactly supported shearlet generators. The key idea behind this seemingly
unreasonable approach is that the classical band-limited generators -- whose Fourier transforms have wedge-like support -- leading
to Parseval frames can be much better approximated by non-separable compactly supported functions than by separable ones.

ShearLab 3D builds on the approach from \cite{Lim2013}, and extends it in two ways, namely to universal shearlets as well as to
the 3D situation. In addition, elaborate numerical experiments comparing ShearLab 3D with the current state-of-the-art algorithms
in geometric multiscale analysis are provided, more precisely with the Nonsubsampled Shearlet Transform in 2D \& 3D \cite{ELL08},
the Nonsubsampled Contourlet Transform in 2D \cite{dCZD2006}, the Fast Discrete Curvelet Transform in 2D \cite{CD1999}, the
Surfacelet Transform in 3D \cite{DL2007}, and finally also with the classical Stationary Wavelet Transform in 2D.
The algorithms were tested with respect to denoising and inpainting
in 2D and 3D as well as separation in 2D (separation in 3D requires additional systems which renders a comparison unfair). Each
time carefully chosen performance measures were specified and served as an objective basis for the comparisons. It could be
shown that indeed ShearLab 3D outperforms the other algorithms in most tasks, also concerning speed. It is interesting to note
that with respect to video denoising ShearLab 3D as a universally applicable tool is only marginally beaten by the specifically
for this task designed BM3D algorithm \cite{MBFE2012}.

In the spirit of reproducible research \cite{DMSSU09}, ShearLab 3D as well as the codes for all comparisons are freely accessible on
\url{www.ShearLab.org}.

\subsection{Outline}

This paper is organized as follows. We first review the main definitions and results from frame and shearlet
theory in Section \ref{sec:shearlets}. Section \ref{sec:digitransform} is then devoted to a detailed discussion
of the faithful algorithmic realization of the shearlet transform (and its inverse), which is implemented in
ShearLab 3D. After a short review of the wavelet transform in Subsection \ref{subsec:wavelets}, the 2D
(Subsection \ref{subsec:2Dtrafo}) and 3D (Subsection \ref{subsec:3Dtrafo}) algorithms are discussed,
each time first in the parabolic and then in the general situation, and finally in Subsection \ref{subsec:inverse}
the inverse transform is presented. The actual MATLAB implementation in ShearLab 3D is described in
Section \ref{sec:implementation}. Finally, elaborate numerical experiments are presented (Section \ref{sec:experiments})
ranging from denoising to video inpainting in Subsections \ref{subsec:denoising} to \ref{subsec:v_inpainting}
each time carefully comparing (Subsection \ref{subsec:discussion}) ShearLab 3D with the state-of-the-art algorithms described
in Subsection \ref{subsec:comparison}.

\section{Discrete Shearlet Systems} \label{sec:shearlets}

In this section, we will state the main definitions and necessary results about shearlet systems for $L^2(\R^2)$ and $L^2(\R^3)$.
We start though with a brief review of frame theory which -- with the notion of frames -- provides a functional analytic concept
to generalize the setting of orthonormal bases. In fact, shearlet systems do not form bases, but require this extended concept.

\subsection{Review: Frame Theory}

Frame theory is nowadays used when redundant, yet stable expansions are required. A sequence $(\varphi_i)_{i \in I}$
in a Hilbert space $\cH$ is called {\em frame for $\cH$}, if there exist constants $0 < A \le B < \infty$ such that
\[
A \|x\|^2 \le \sum_{i \in I} \absip{x}{\varphi_i}^2 \le B \|x\|^2 \quad \mbox{for all } x \in \cH.
\]
If $A=B$ is possible, the frame is referred to as {\em tight}, in case $A=B=1$ as a {\em Parseval frame}. One application
of frames is the analysis of elements in a Hilbert space, which is achieved by the {\em analysis operator} given by
\[
T : \cH \to \ell_2(I), \quad T(x) = (\ip{x}{\varphi_i})_{i \in I}.
\]
Reconstruction of each element $x \in \cH$ from $Tx$ is possible by the {\em frame reconstruction formula} given by
\beq \label{eq:framerec}
x = \sum_{i \in I} \ip{x}{\varphi_i} S^{-1}\varphi_i,
\eeq
where $Sx=\sum_{i \in I} \ip{x}{\varphi_i} \varphi_i$ is the {\em frame operator} associated with the frame
$(\varphi_i)_{i \in I}$. We remark that for tight frames, the frame operator is just a multiple of the identify,
hence easily invertible. For general frames, it might be difficult to use \eqref{eq:framerec} in practise,
since inverting $S$ might be numerically unfeasible, in particular if $B/A$ is large. In those cases, the
so-called {\em frame algorithm}, for instance, described in \cite{Chr03}, can be applied. Moreover, in \cite{Gr93},
the {\em Chebyshev method} and the {\em conjugate gradient methods} were introduced, which are significantly
better adapted to frame theory leading to faster convergence than the frame algorithm.

\subsection{Universal Shearlet Systems} \label{subsec:discreteshearlets}

We next introduce 2D and 3D universal shearlet systems, each time first discussing the parabolic case and then
extending it to the general case. For more information, we refer to \cite{KL2012} and \cite{GK15}.

\subsubsection{2D Situation}

\begin{figure}[ht]
\hspace*{1cm}
\begin{center}
\includegraphics[width=4cm]{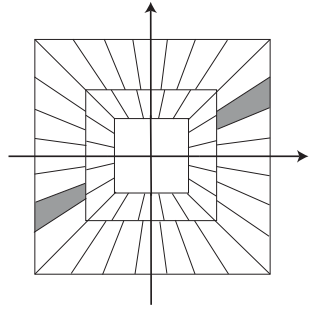}
\end{center}
\caption{Tiling of the frequency plane induced by the Parseval frame of classical shearlets. }
\label{fig:tilingfrequency}
\end{figure}

Shearlet systems can be regarded as consisting of certain generating functions whose resolution is changed by
a parabolic scaling matrix $A_{2^j}$ or $\tilde{A}_{2^j}$ defined by
\[
A_{2^j} = \begin{pmatrix} 2^j & 0 \\ 0 & 2^{j/2} \end{pmatrix} \quad \mbox{and} \quad
\tilde{A}_{2^j} = \begin{pmatrix} 2^{j/2} & 0 \\ 0 & 2^{j} \end{pmatrix},
\]
whose orientation is changed by a shearing matrix $S_k$ or $S_k^T$ defined by
\[
S_k = \begin{pmatrix} 1 & k \\ 0 & 1\end{pmatrix},
\]
and whose position is changed by translation. More precisely. a shearlet system -- sometimes in this form also
referred to as {\em cone-adapted} due to the fact that it is adapted to a cone-like partition in frequency
domain (cf. Figure \ref{fig:tilingfrequency}) -- is defined as follows.

\begin{definition} \label{defi:2Dshearletsystem}
For $\phi, \psi, \tilde{\psi} \in L^2(\R^2)$ and $c=(c_1,c_2) \in (\R_+)^2$, the \emph{shearlet system}
$SH(\phi,\psi,\tilde{\psi}; c)$ is defined by
\[
SH(\phi,\psi,\tilde{\psi}; c) = \Phi(\phi; c_1) \cup \Psi(\psi; c) \cup \tilde{\Psi}(\tilde{\psi}; c),
\]
where
\begin{align*}
  \Phi(\phi; c_1) &= \{\phi_m = \phi(\cdot- c_1 m) : m \in \Z^2\},\\
  \Psi(\psi; c) &= \{\psi_{j,k,m} = 2^{\frac34 j} \psi(S_k A_{2^j}\cdot  - M_c m) : j \ge 0,  |k| \leq \lceil2^{j/2}\rceil, m \in \Z^2\},\\
  \tilde{\Psi}(\tilde{\psi}; c) &= \{\tilde{\psi}_{j,k,m} =  2^{\frac34 j} \tilde{\psi}(S_k^T\tilde{A}_{2^j}\cdot  -\tilde{M}_c m): j \ge 0,
  |k| \leq \lceil 2^{j/2} \rceil, m \in \Z^2 \},
\end{align*}
with
\[
M_c =
\begin{pmatrix}
c_1 & 0 \\ 0 & c_2
\end{pmatrix}
\quad \text{and} \quad \tilde{M}_c =
\begin{pmatrix}
c_2 & 0 \\ 0 & c_1
\end{pmatrix}.
\]
\end{definition}

Its associated transform maps functions to the sequence of shearlet coefficients, hence is merely the associated
analysis operator.

\begin{definition}
Set $\Lambda = \N_0 \times \{-\lceil2^{j/2}\rceil, \ldots, \lceil2^{j/2}\rceil\} \times \Z^2.$
Further, let $SH(\phi,\psi,\tilde{\psi}; c)$ be a shearlet system and retain the notions from
Definition \ref{defi:2Dshearletsystem}. Then the associated {\em shearlet transform} of
$f \in L^2(\R^2)$ is the mapping defined by
\[
f \to \cSH_{\phi,\psi,\tilde{\psi}} f(m',(j,k,m),(\tilde{j},\tilde{k},\tilde{m}))
= (\langle f,\phi_{m'}\rangle,\langle f,\psi_{j,k,m}\rangle,\langle f,\tilde{\psi}_{\tilde{j},\tilde{k},\tilde{m}}\rangle),
\]
where
\[
(m',(j,k,m),(\tilde{j},\tilde{k},\tilde{m})) \in \Z^2 \times \Lambda \times \Lambda.
\]
\end{definition}

A (historically) first class of generators are so-called classical shearlets, which are defined
as follows. Let $\psi \in L^2(\R^2)$ be defined by
\begin{equation}\label{eq:bandlimit}
\hat{\psi}(\xi) = \hat{\psi}(\xi_1,\xi_2) = \hat{\psi}_1(\xi_1) \, \hat{\psi}_2(\tfrac{\xi_2}{\xi_1}),
\end{equation}
where $\psi_1 \in L^2(\R)$ is a discrete wavelet in the sense that it satisfies the discrete Calder\'{o}n condition, given by
\[
\sum_{j \in \Z}|\hat\psi_1(2^{-j}\xi)|^2 = 1 \quad \mbox{for a.e. } \xi \in \R,
\]
with $\hat{\psi}_1 \in C^\infty(\R)$ and  $\supp \hat{\psi}_1 \subseteq [-\frac1{2},-\frac1{16}] \cup [\frac1{16},\frac12]$,
and $\psi_2 \in L^2(\R)$ is a bump function in the sense that
\[
\sum_{k = -1}^{1}|\hat\psi_2(\xi+k)|^2 = 1 \quad \mbox{for a.e. } \xi \in [-1,1],
\]
satisfying $\hat{\psi}_2 \in C^\infty(\R)$ and $\supp \hat{\psi}_2 \subseteq [-1,1]$. Then $\psi$ is called
a {\em classical shearlet}. With small modifications of the boundary elements, classical shearlets lead to
a Parseval frame for $L^2(\R^2)$ \cite{GKL05}. The induced tiling of the frequency plane is illustrated in
Figure~\ref{fig:tilingfrequency}.

Some years later, compactly supported shearlets have been studied. It was shown that a large class of
compactly supported generators yield shearlet frames with controllable frame bounds \cite{KKL2012}.
One for numerical algorithms particularly interesting special case are separable generators
given by $\psi = \psi_1 \otimes \phi_1$, which generate shearlet frames provided that the 1D
wavelet function $\psi_1$ and the 1D scaling function $\phi_1$ are sufficiently smooth and $\psi_1$ has
sufficient vanishing moments.

Let us now turn to the more flexible universal shearlets, which were introduced in \cite{GK15}. Their
definition requires an extension of the scaling matrix to insert a parameter $\alpha \in (0,2)$
measuring the degree of anisotropy. For this, let $A_{\alpha, 2^j}$ and $\tilde{A}_{\alpha, 2^j}$
be defined by
\[
A_{\alpha, 2^j} = \begin{pmatrix} 2^j & 0 \\ 0 & 2^{\alpha j/2} \end{pmatrix} \quad \mbox{and} \quad
\tilde{A}_{2^j} = \begin{pmatrix} 2^{\alpha j/2} & 0 \\ 0 & 2^{j} \end{pmatrix},
\]
A universal shearlet system can then be defined as follows.

\begin{definition} \label{defi:2Dshearletsystem_alpha}
For $\phi, \psi, \tilde{\psi} \in L^2(\R^2)$,  $\alpha = (\alpha_j)_j$, $\alpha_j \in (0,2)$ for each scale $j$, and $c=(c_1,c_2) \in (\R_+)^2$,
the \emph{universal shearlet system} $SH(\phi,\psi,\tilde{\psi}; \alpha, c)$ is defined by
\[
SH(\phi,\psi,\tilde{\psi}; \alpha, c) = \Phi(\phi; c_1) \cup \Psi(\psi; \alpha, c) \cup \tilde{\Psi}(\tilde{\psi}; \alpha, c),
\]
where
\begin{align*}
  \Phi(\phi; c_1) &= \{\phi_m = \phi(\cdot- c_1 m) : m \in \Z^2\},\\
  \Psi(\psi; \alpha, c) &= \{\psi_{j,k,m} = 2^{\frac{\alpha_j+1}{4} j} \psi(S_k A_{\alpha_j, 2^j}\cdot  - M_c m) : j \ge 0,  |k| \leq \lceil2^{j(\alpha_j-1)/2}\rceil, m \in \Z^2\},\\
  \tilde{\Psi}(\tilde{\psi}; \alpha, c) &= \{\tilde{\psi}_{j,k,m} =  2^{\frac{\alpha_j+1}{4} j} \tilde{\psi}(S_k^T\tilde{A}_{\alpha_j, 2^j}\cdot  -\tilde{M}_c m): j \ge 0,
  |k| \leq \lceil 2^{j(\alpha_j-1)/2} \rceil, m \in \Z^2 \}.
\end{align*}
\end{definition}

Let us now briefly discuss the situation when all $\alpha_j$ coincide, i.e., $\alpha_0 := \alpha_j$ for all scale $j$.
In this case, $SH(\phi,\psi,\tilde{\psi}; \alpha, c) = SH(\phi,\psi,\tilde{\psi}; c)$ if $\alpha_0 = 1$.
Moreover, if $\alpha_0 = 2$ for any scale $j$, then $SH(\phi,\psi,\tilde{\psi}; \alpha, c)$ becomes an isotropic
wavelet system. It should be also mentioned that for $\alpha_0 \to 0$, the associated universal shearlet
system approaches the system of ridgelets \cite{CD99}.

The associated transform is then defined similarly as in the parabolic situation.

\begin{definition}
Set $\Lambda = \N_0 \times \{-\lceil 2^{j(\alpha_j-1)/2} \rceil, \ldots, \lceil 2^{j(\alpha_j-1)/2} \rceil\} \times \Z^2.$
Further, let $SH(\phi,\psi,\tilde{\psi}; \alpha, c)$ be a shearlet system and retain the notions from
Definition \ref{defi:2Dshearletsystem_alpha}. Then the associated {\em universal shearlet transform}
of $f \in L^2(\R^2)$ is the mapping defined by
\[
f \to \cSH_{\phi,\psi,\tilde{\psi}} f(m',(j,k,m),(\tilde{j},\tilde{k},\tilde{m}))
= (\langle f,\phi_{m'}\rangle,\langle f,\psi_{j,k,m}\rangle,\langle f,\tilde{\psi}_{\tilde{j},\tilde{k},\tilde{m}}\rangle),
\]
where
\[
(m',(j,k,m),(\tilde{j},\tilde{k},\tilde{m})) \in \Z^2 \times \Lambda \times \Lambda.
\]
\end{definition}

In \cite{GK15} it was shown that there exists an abundance of scaling sequences $\alpha = (\alpha_j)_j$ such
that with a small modification of classical shearlets the system $SH(\phi,\psi,\tilde{\psi}; \alpha, c)$
yields a Parseval frame for $L^2(\R^2)$.

\subsubsection{3D Situation}

Turning now to the 3D situation and starting again with the parabolic case, it is apparent that the 2D parabolic
scaling matrix $A_{2^j}$ can be extended either by $\text{diag}(2^j, 2^{j/2}, 2^{j})$ or $\text{diag}(2^j, 2^{j/2}, 2^{j/2})$.
The first case however generates `needle-like shearlets' for which a frame property seems highly unlikely. Hence
the second case is typically considered.


For the sake of brevity, we next immediately present the general case of 3D universal shearlets. Then, for
$\alpha \in (0,2)$, we set
\[
A_{\alpha, 2^j} = \begin{pmatrix} 2^{ j} & 0 & 0 \\ 0 & 2^{\alpha j/2} & 0 \\ 0 & 0 & 2^{\alpha j/2} \end{pmatrix} \hspace*{-0.1cm}, \;
\tilde{A}_{\alpha, 2^j} = \begin{pmatrix} 2^{\alpha j/2} & 0 & 0\\ 0 & 2^{j} & 0 \\ 0 & 0 & 2^{\alpha j/2}\end{pmatrix} \hspace*{-0.1cm}, \; \mbox{and} \;
\breve{A}_{\alpha, 2^j} = \begin{pmatrix} 2^{\alpha j/2} & 0 & 0\\ 0 & 2^{\alpha j/2} & 0 \\ 0 & 0 & 2^{j}\end{pmatrix}.
\]
The shearing matrices are now associated with a parameter $k = (k_1,k_2) \in \Z^2$ and defined by
\[
S_{k} = \begin{pmatrix} 1 & k_1 & k_2 \\ 0 & 1 & 0 \\ 0 & 0 & 1 \end{pmatrix}, \quad
\tilde{S}_{k} = \begin{pmatrix} 1 & 0 & 0\\ k_1 & 1 & k_2 \\ 0 & 0 & 1\end{pmatrix}, \quad \mbox{and} \quad
\breve{S}_{k} = \begin{pmatrix} 1 & 0 & 0\\ 0 & 1 & 0 \\ k_1 & k_2 & 1\end{pmatrix}.
\]
Finally, for $c_1,c_2 \in \R_+$, translation will be given by $M_c = \text{diag}(c_1,c_2,c_2)$,
$\tilde{M}_c = \text{diag}(c_2,c_1,c_2)$, and $\breve{M}_c = \text{diag}(c_2,c_2,c_1)$. With these
notions, a 3D universal shearlet system can then be defined as follows.

\begin{definition} \label{defi:3Dshearletsystem_alpha}
For $\phi, \psi, \tilde{\psi}, \breve{\psi} \in L^2(\R^2)$,  $\alpha = (\alpha_j)_j$, $\alpha_j \in (0,2)$ for each scale $j$, and
$c=(c_1,c_2) \in (\R_+)^2$, the \emph{universal shearlet system} $SH(\phi,\psi,\tilde{\psi}, \breve{\psi}; \alpha, c)$ is defined by
\[
SH(\phi,\psi,\tilde{\psi}, \breve{\psi}; \alpha, c) = \Phi(\phi; c_1) \cup \Psi(\psi; \alpha, c) \cup \tilde{\Psi}(\tilde{\psi}; \alpha, c)
\cup \breve{\Psi}(\breve{\psi}; \alpha, c),
\]
where
\begin{align*}
  \Phi(\phi; c_1) &= \{\phi_m = \phi(\cdot- c_1 m) : m \in \Z^2\},\\
  \Psi(\psi; \alpha, c) &= \{\psi_{j,k,m} = 2^{\frac{\alpha_j+1}{4} j} \psi(S_k A_{\alpha_j, 2^j}\cdot  - M_c m) : j \ge 0,  |k| \leq \lceil2^{j(\alpha_j-1)/2}\rceil, m \in \Z^2\},\\
  \tilde{\Psi}(\tilde{\psi}; \alpha, c) &= \{\tilde{\psi}_{j,k,m} =  2^{\frac{\alpha_j+1}{4} j} \tilde{\psi}(S_k^T\tilde{A}_{\alpha_j, 2^j}\cdot  -\tilde{M}_c m): j \ge 0,
  |k| \leq \lceil 2^{j(\alpha_j-1)/2} \rceil, m \in \Z^2 \},\\
  \breve{\Psi}(\breve{\psi}; \alpha, c) &= \{\breve{\psi}_{j,k,m} =  2^{\frac{\alpha_j+1}{4} j} \breve{\psi}(S_k^T\breve{A}_{\alpha_j, 2^j}\cdot  -\breve{M}_c m): j \ge 0,
  |k| \leq \lceil 2^{j(\alpha_j-1)/2} \rceil, m \in \Z^2 \}.
\end{align*}
\end{definition}

The associated transform is then defined as follows.

\begin{definition}
Set $\Lambda = \N_0 \times \{-\lceil 2^{j(\alpha_j-1)/2}\rceil, \ldots, \lceil 2^{j(\alpha_j-1)/2} \rceil\}^2 \times \Z^3.$
Further, let $SH(\phi,\psi,\tilde{\psi},\breve{\psi}; \alpha, c)$ be a shearlet system and retain the notions from
Definition \ref{defi:3Dshearletsystem_alpha}. Then the associated {\em universal shearlet transform}
of $f \in L^2(\R^2)$ is the mapping, which maps $f$ to
\[
\cSH_{\phi,\psi,\tilde{\psi},\breve{\psi}} f(m',(j,k,m),(\tilde{j},\tilde{k},\tilde{m}),(\breve{j},\breve{k},\breve{m}))
= (\langle f,\phi_{m'}\rangle,\langle f,\psi_{j,k,m}\rangle,\langle f,\tilde{\psi}_{\tilde{j},\tilde{k},\tilde{m}}\rangle,\langle f,\breve{\psi}_{\breve{j},\breve{k},\breve{m}}\rangle),
\]
where
\[
(m',(j,k,m),(\tilde{j},\tilde{k},\tilde{m}),(\breve{j},\breve{k},\breve{m})) \in \Z^3 \times \Lambda \times \Lambda \times \Lambda.
\]
\end{definition}

\begin{figure}[ht]
\hspace*{1cm}
\begin{center}
\includegraphics[width=4cm]{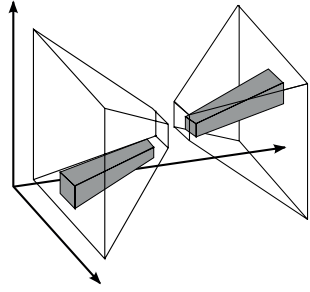}
\end{center}
\caption{Tiling of the frequency domain induced by 3D shearlets. }
\label{fig:3d_shearlets}
\end{figure}

In the situation of parabolic scaling it was shown in \cite{KLL2012} that a large class of
compactly supported generators yields shearlet frames with controllable frame bounds. This is
in particular the case for separable generators $\psi = \psi_1  \otimes \phi_1 \otimes \tilde{\phi}_1$
for which $\psi_1$ is a 1D wavelets with sufficiently many vanishing moments and
$\phi_1, \tilde{\phi}_1$ are sufficiently smooth 1D scaling functions. Still in the parabolic case,
the situation of band-limited (classical) 3D shearlets was studied in \cite{Guo2010}, which -- similar
to the 2D situation -- form a Parseval frame for $L^2(\R^3)$ with a small modification of the elements 
near the seam lines. An illustration of how 3D shearlets tile the frequency domain is provided in
Figure \ref{fig:3d_shearlets}.

\subsection{Previous Implementations}\label{subsec:previousimplements}

We now briefly review previous implementations of the shearlet transform. It should be emphasized though
that all implementations so far only focussed on the parabolic case, i.e., $\alpha_j = 1$ for all
scales $j$.

\subsubsection{Fourier Domain Approaches}\label{subsubsec:bandlimited}

Band-limited shearlet systems provide a precise partition of the frequency plane due to the fact that
the Fourier transforms of all elements are compactly supported and that they form a tight frame. Hence
it seems most appropriate to implement the associated transforms via a {\em Fourier domain approach},
which aims to directly produce the same frequency tiling.

A first numerical implementation using this approach was discussed in \cite{ELL08} as a cascade of a subband
decomposition based on the Laplacian Pyramid filter, which was then followed by a step containing directional
filtering using the Pseudo-Polar Discrete Fourier Transform.
Another approach was suggested in \cite{KSZ12}. The main idea here is to employ a carefully weighted Pseudo-Polar
transform with weights ensuring (almost) isometry. This step is then followed by appropriate windowing and
the inverse FFT applied to each windowed part.

\subsubsection{Spatial Domain Approaches}\label{subsubsec:compactsupport}

We refer to a numerical realization of a shearlet transform as a {\em spatial domain approach} if the filters
associated with the transform are implemented by a convolution in the spatial domain. Whereas Fourier-based
approaches were only utilized for band-limited shearlet transforms, the range of spatial domain-based approaches
is much broader and basically justifiable for a transform based on any shearlet system. We now present the main
contributions.

In the paper \cite{ELL08} already referred to in Subsection \ref{subsubsec:bandlimited}, also a spatial domain
approach is discussed. This implementation utilizes directional filters, which are obtained as approximations of
the inverse Fourier transforms of digitized band-limited window functions. A numerical realization specifically
focussed on separable shearlet generators $\psi^{\text{sep}}$ given by $\psi^{\text{sep}} = \psi_1 \otimes \phi_1$
-- which includes certain compactly supported shearlet frames (cf. Subsection \ref{subsec:discreteshearlets}) --
was derived in \cite{Lim2010}.  This algorithm  enables the application of fast transforms separably along both axes,
even if the corresponding shearlet transform is not associated with a tight frame. The most faithful, efficient,
and numerically stable (in the sense of closeness to tightness) digitalization of the shearlet transform was derived
in \cite{Lim2013} by utilizing non-separable compactly supported shearlet generators, which best approximate the
classical band-limited generators. In fact, the implementation of the universal shearlet transform we will discuss
in this paper will be based on this work.

There exist two other approaches, which though have not been numerically tested yet. The one introduced in \cite{KS09}
explores the theory of subdivision schemes inserting directionality, and leads to a version of the shearlet transform
which admits an associated multiresolution analysis structure. In close relation to this algorithmic realization,
in \cite{HKS11} a general unitary extension principle is proven, which -- for the shearlet setting -- provides
equivalent conditions for the filters to lead to a shearlet frame.

\section{Digital Shearlet Transform}\label{sec:digitransform}

In this section, we will introduce and discuss the algorithms which are implemented in ShearLab 3D. We will start with a brief review
of the digital wavelet transform in Subsection \ref{subsec:wavelets}, which parts of the digital shearlet transform will
be based upon. In Subsection \ref{subsec:2Dtrafo} the 2D forward shearlet transform will be discussed, first the
parabolic version, followed by a general version with freely chosen $\alpha$. The 3D version is then presented in
Subsection \ref{subsec:3Dtrafo}. Finally, the inverse shearlet transform both for 2D and 3D is detailed in Subsection
\ref{subsec:inverse}.

\subsection{Digital Wavelet Transform}\label{subsec:wavelets}

We start by recalling the notion of the discrete Fourier transform of a sequence $\{a(n)\}_{n \in \mathbb Z^d} \in \ell^2(\mathbb Z^d)$,
which is defined by
\[
\hat{a}(\xi) = \sum_{n \in \mathbb Z^d}a(n) e^{-2\pi i n \cdot \xi}.
\]
In the sequel, we will also use the notion $\overline{a}(n) = a(-n)$ for $n \in \mathbb Z^d$. We further wish to refer the reader
not familiar with wavelets to the books \cite{Dau92,Mal08}.

Let now $\psi^1$ and $\phi^1$ $\in L^2(\mathbb R)$ be a wavelet and an associated scaling function, respectively, satisfying the two scale
relations
\begin{equation}\label{eq:refine01}
\phi^1(x_1) = \sum_{n_1 \in \mathbb Z} h(n_1) \sqrt{2}\phi^1(2 x_1 - n_1)
\end{equation}
and
\begin{equation}\label{eq:refine02}
\psi^1(x_1) = \sum_{n_1 \in \mathbb Z} g(n_1) \sqrt{2}\phi^1(2 x_1 - n_1).
\end{equation}
Then, for each $j > 0$ and $n_1 \in \Z$, the associated wavelet function $\psi^1_{j,n_1} \in L^2(\R)$ is defined by
\[
\psi^1_{j,n_1}(x_1) = 2^{j/2}\psi^1(2^j x_1 - n_1).
\]

Now let $f$ be a function on $\R$, for which we assume an expansion of the type
\[
f^{\text{1D}}(x_1) = \sum_{n_1 \in \mathbb Z} f^{\text{1D}}_{J}(n_1) 2^{J/2}\phi_1(2^J x_1 - n_1)
\]
for a fixed, sufficiently large $J > 0$. To derive a formula for the associated wavelet coefficients, for each $j>0$, let
$\{h_j(n_1)\}_{n_1 \in \mathbb Z}$ and $\{g_j(n_1)\}_{n_1 \in \mathbb Z}$ denote the Fourier coefficients of the trigonometric polynomials
\begin{equation}\label{eq:lowhigh}
\hat{h}_{j}(\xi_1) = \prod_{k=0}^{j-1}\hat{h}(2^k \xi_1)
\quad \text{and} \quad
\hat{g}_{j}(\xi_1) = \hat{g}\bigl(\frac{2^{j}\xi_1}{2}\bigr)\hat{h}_{j-1}(\xi_1)
\end{equation}
with $\hat{h}_0 \equiv 1$. Then using \eqref{eq:refine01} and \eqref{eq:refine02} and setting $w_j \equiv g_{J-j}$ as well as
$\Phi^1(n_1) = \langle \phi^1(\cdot),\phi^1(\cdot-n_1)\rangle$, the wavelet coefficients $\langle f^{\text{1D}},\psi^1_{j,n_1}\rangle$
can be computed by the discrete formula given by
\[
\langle f^{\text{1D}},\psi^1_{j,m_1}\rangle = \overline{w}_{j}*(f^{\text{1D}}_J*\phi^1)(2^{J-j}\cdot m_1).
\]
In particular, when $\phi_1$ is an orthonormal scaling function, this expression reduces to
\begin{equation}\label{eq:owt}
\langle f^{\text{1D}},\psi^1_{j,m_1}\rangle = (\overline{w}_{j}*f^{\text{1D}}_J)(2^{J-j}\cdot m_1).
\end{equation}
This last formula is the common {\em 1D digital wavelet transform}.

This transform can now be easily extended to the multivariate case by using tensor products. Similar to the 1D case, we 
consider a 2D function $f$ given by
\begin{equation}\label{eq:code2}
f(x) = \sum_{n \in \mathbb Z^2} f_J(n) 2^J \phi(2^Jx_1-n_1,2^Jx_2-n_2),
\end{equation}
where $\phi(x) = \phi_1\otimes\phi_1(x)$. Assume that $\phi_1$ is an orthonormal scaling function. Then, using  \eqref{eq:owt},
it is a straightforward calculation to show that the {\em 2D digital wavelet transform}, which computes wavelet coefficients for $f$,
is of the form
\begin{eqnarray*}
\langle f,\psi^1_{j,m_1} \otimes \phi^1_{j,m_1} \rangle &=& (\overline{g_{J-j}\otimes h_{J-j}} * f )(2^{J-j}\cdot m),
\\
\langle f,\phi^1_{j,m_1} \otimes \psi^1_{j,m_2} \rangle &=& (\overline{h_{J-j}\otimes g_{J-j}} * f )(2^{J-j}\cdot m),
\end{eqnarray*}
and
\[
\langle f,\psi^1_{j,m_1} \otimes \psi^1_{j,m_1} \rangle = (\overline{g_{J-j}\otimes g_{J-j}} * f )(2^{J-j}\cdot m).
\]

\subsection{2D Shearlet Transform} \label{subsec:2Dtrafo}

We will now describe the algorithmic realization of the transform associated with a universal compactly supported
shearlet system $SH(\phi,\psi,\tilde{\psi}; \alpha, c)$ (as introduced in Subsection \ref{subsec:discreteshearlets}) used in ShearLab 3D. We remark that the subset
$\Phi(\phi; c_1)$ is merely the scaling part coinciding with the wavelet scaling part. Moreover, it will be sufficient to consider
shearlets from $\Psi(\psi; \alpha, c)$ as the same arguments apply to $\tilde{\Psi}(\tilde{\psi}; \alpha, c)$ except for switching the order of variables.

We start with the digitalization of the universal shearlet transform for the parabolic case, i.e., $\alpha_j = 1$ for all scales $j$, which was
initially introduced in \cite{Lim2013}. We then extend this algorithmic approach to the situation of a general parameter $\alpha = (\alpha_j)_j$,
allowing $\alpha_j$ to differ for each scale $j$.

\subsubsection{The Parabolic Case} \label{subsubsec:parabolic}

\begin{figure}[h]
\begin{center}
\includegraphics[width=1.4in,height=1.2in]{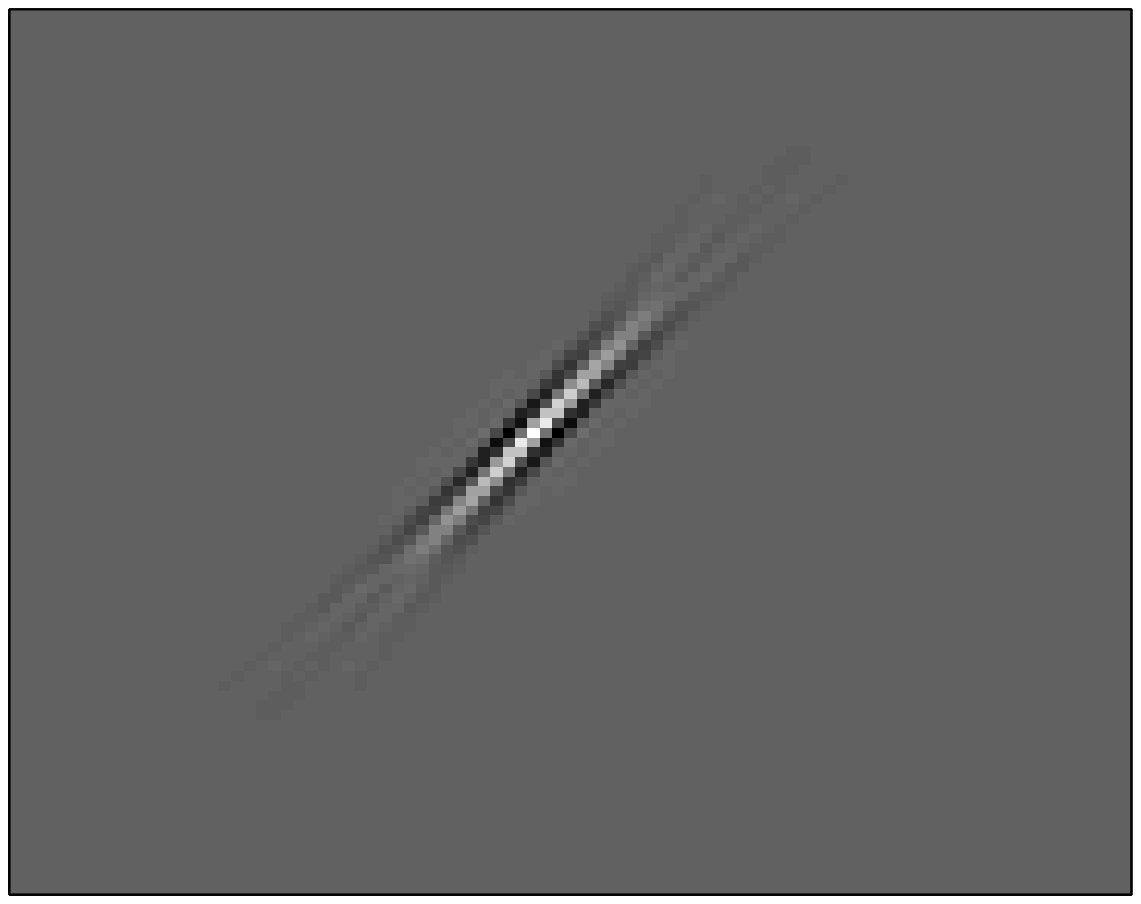}
\hspace{-20pt}
\includegraphics[width=1.4in,height=1.2in]{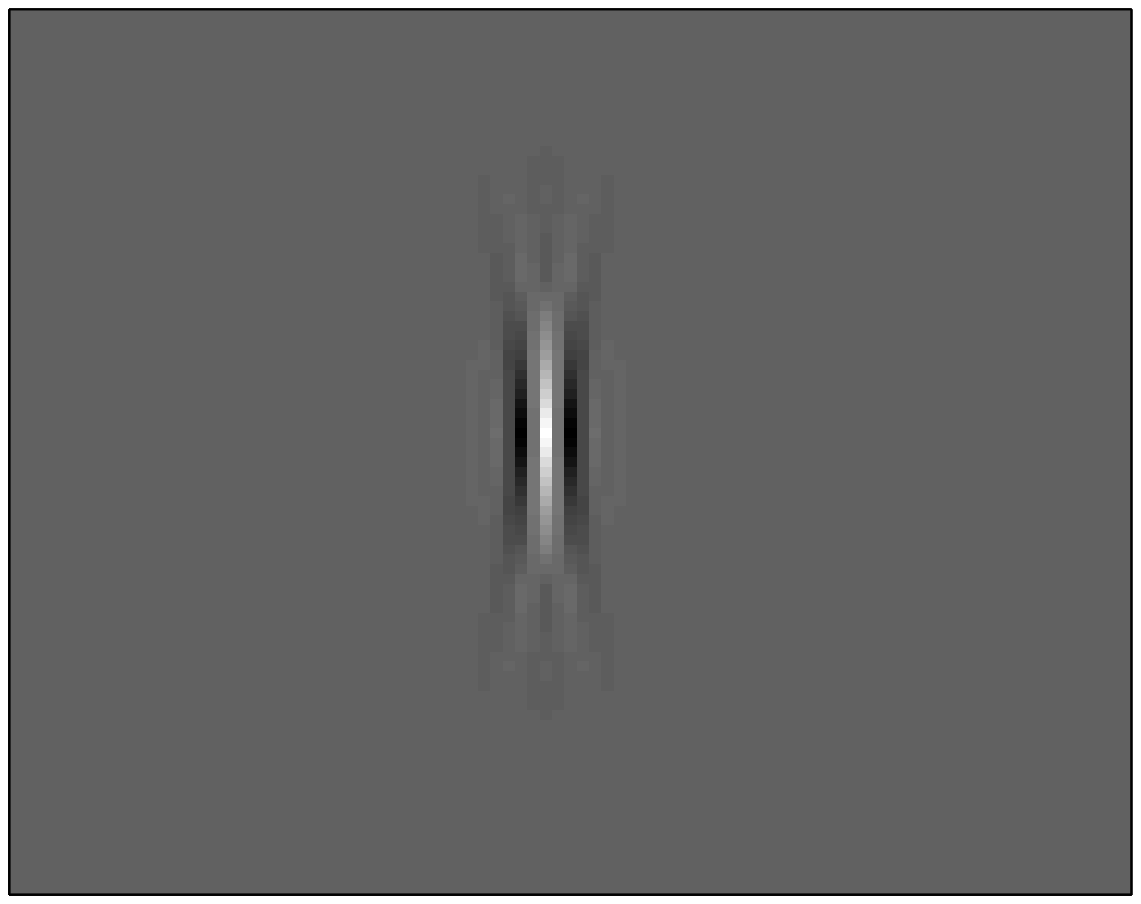}
\hspace{-20pt}
\includegraphics[width=1.4in,height=1.2in]{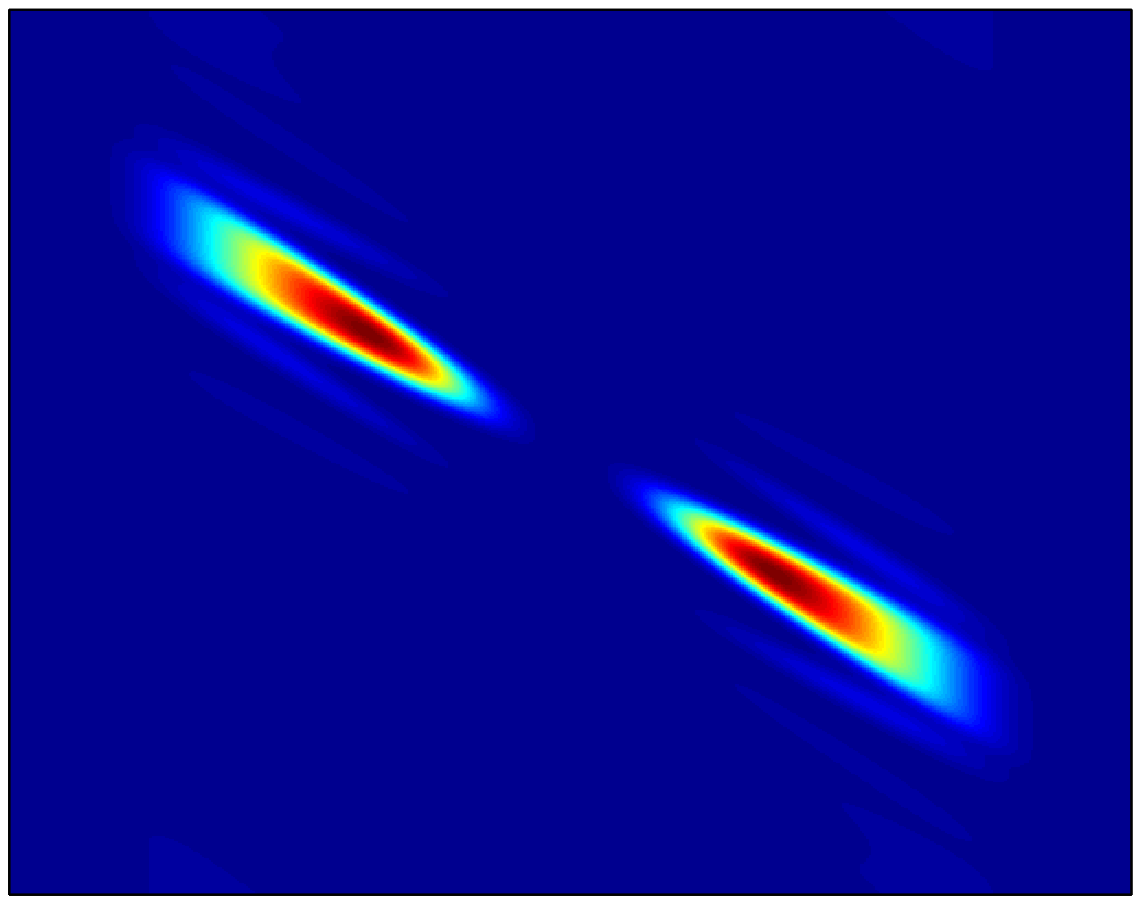}
\hspace{-20pt}
\includegraphics[width=1.4in,height=1.2in]{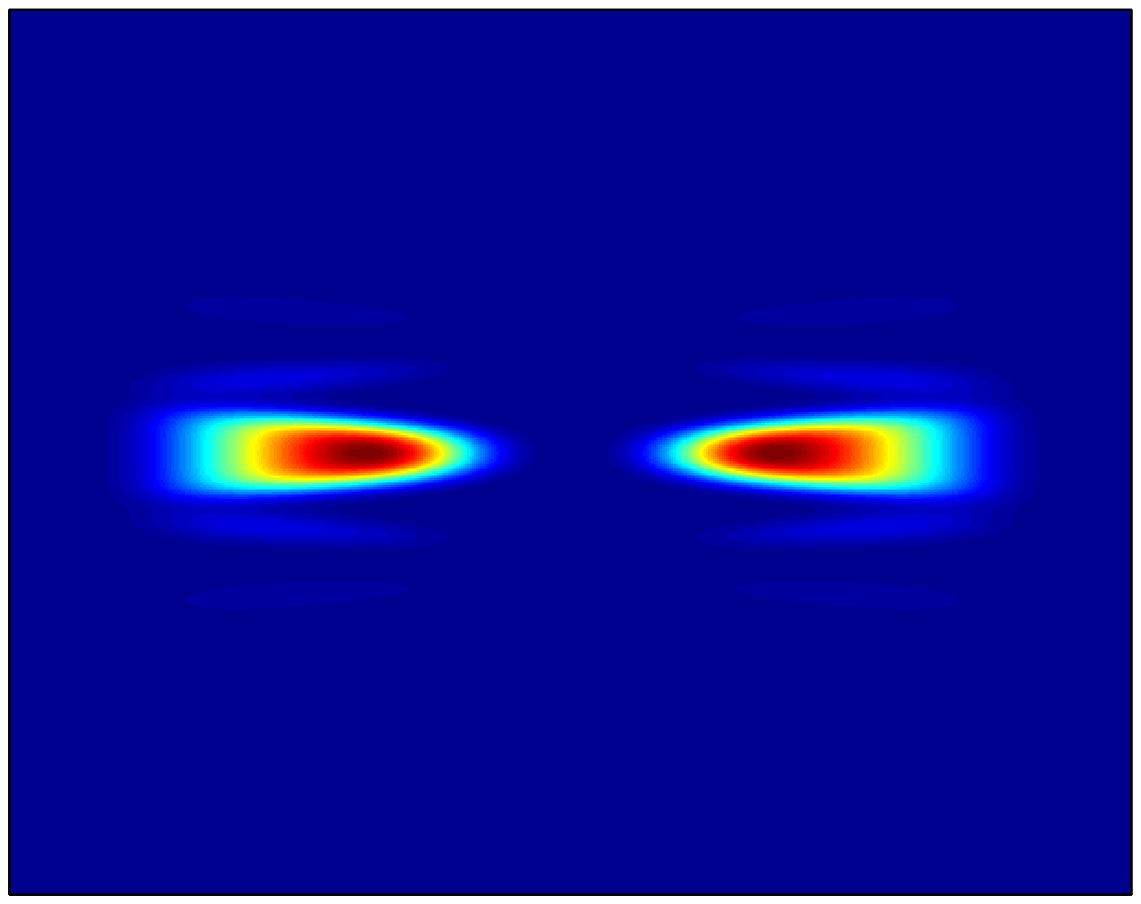}
\put(-325,-10){(a)}
\put(-235,-10){(b)}
\put(-145,-10){(c)}
\put(-55,-10){(d)}
\caption{ Shearlets $\psi_{j,k,m}$ in the time/frequency domain. (a)--(b) : Shearlets in the spatial domain. (c)--(d) : Shearlets in the frequency domain.}
\label{fig:shearlets}
\end{center}
\end{figure}

First, the shearlet generator $\psi$ needs to be chosen. In Subsection \ref{subsubsec:compactsupport}, the choice of a separable shearlet
generator $\psi = \psi_1 \otimes \phi_1$ was discussed, which generated a shearlet frame provided that the 1D wavelet function $\psi_1$ and
the 1D scaling function $\phi_1$ are sufficiently smooth and $\psi_1$ has sufficient vanishing moments. However, significantly improved
numerical results can be achieved by choosing a {\em non-separable} generator $\psi$ such as
\begin{equation}\label{eq:nonsep}
{\hat{\psi} (\xi) = P(\xi_1/2,\xi_2)\hat \psi^{\text{sep}}(\xi)},
\end{equation}
where the trigonometric polynomial $P$ is a 2D fan filter (cf. \cite{DV2005}, \cite{dCZD2006}). The reason for this is the fact, that
non-separability allows the Fourier transform of $\psi$ to have a wedge shaped essential support, thereby well approximating the
Fourier transform of a classical shearlet. In particular, with a suitable choice for a 2D fan filter $P$, we have
\[
P(\xi_1/2,\xi_2)\hat\phi_1(\xi_2) \approx \hat \psi_2(\frac{\xi_2}{\xi_1})
\]
and
\[
P(\xi_1/2,\xi_2)\hat\psi_1(\xi_1)\hat\phi_1(\xi_2) \approx \hat \psi_1(\xi_1)\hat \psi_2(\frac{\xi_2}{\xi_1})
\]
where $\hat \psi_1(\xi_1)\hat \psi_2(\frac{\xi_2}{\xi_1})$ is the classical shearlet generator defined in \eqref{eq:bandlimit}. 
Compared with separable compactly supported shearlet generators, this property does indeed
not only improve the frame bounds of the associated system, but also  improves the directional selectivity significantly. Figure
\ref{fig:shearlets} shows some non-separable compactly supported shearlet generators both in time and frequency domain.
One can in fact construct compactly supported functions $\psi_1$ and $\phi_1$, and a finite 2D fan filter $P$ such that
\begin{equation}\label{eq:cond1}
\text{inf}\{|\hat \psi_1(\xi_1)|^2 : 1/2 < |\xi_1| < 1\} > \delta_1, \quad
\text{inf}\{|\hat \phi_1(\xi_1)|^2 : -1/2 < \xi_1 < 1/2 \} > \delta_2,
\end{equation}
and
\begin{equation}\label{eq:cond2}
\text{inf}\{|P(\xi)|^2 : 1/4 < |\xi_1| < 1/2, |\xi_2/\xi_1| < 1\} > \delta_3
\end{equation}
with some $\delta_1,\delta_2$ and $\delta_3 > 0$. This implies
\[
|\hat \phi(\xi)|^2 + \sum_{j \ge 0} \sum_{|k| \leq \lceil 2^{j/2} \rceil} \Bigl( |\hat \psi_{j,k,0}(\xi)|^2+|\hat{\tilde{\psi}}_{j,k,0}(\xi)|^2\Bigr) > \min(\delta_2^2,\delta_1\delta_2\delta_3)
\]
with $\phi = \phi_1 \otimes \phi_1$ and $\tilde{\psi}(x_1,x_2) = \psi(x_2,x_1)$.
This inequality provides a lower frame bound provided that $\psi_1$ and $\phi_1$ decay sufficiently fast in frequency and $\psi_1$ 
has sufficient vanishing moments. One can also obtain an upper frame bound with $\psi$ generated by those 1D functions. We 
refer to \cite{KKL2012} for more details.

The task is now to derive a digital formulation for the computation of the associated shearlet coefficients
$\langle f,\psi_{j,k,m}\rangle$ for $j = 0,\dots,J-1$ of a function $f$ given as in \eqref{eq:code2}, where
\[
 {\psi_{j,k,m}(x) = 2^{\frac{3}{4}j}\psi(S_kA_{2^j}x-M_{c_j}m)}
\]
with the sampling matrix given by $M_{c_j} = \text{diag}(c_1^j,c_2^j)$. Without loss of generality, we will from
now on assume that $j/2$ is integer; otherwise $\lfloor j/2 \rfloor$ would need to be taken.

We first observe that
\[
S_kA_{2^j} = A_{2^{j}}S_{k/2^{j/2}},
\]
which implies that
\begin{equation} \label{eq:shear}
\psi_{j,k,m}(\cdot) = \psi_{j,0,m}(S_{k/2^{j/2}}\cdot)
\end{equation}
Thus, our strategy for discretizing $\psi_{j,k,m}$ consists of the following two parts:
\begin{enumerate}
\item[(I)] Faithful discretization of $\psi_{j,0,m}$ using the structure of the multiresolution analysis associated with \eqref{eq:nonsep}.
\item[(II)] Faithful discretization of the shear operator $S_{k/2^{j/2}}$.
\end{enumerate}

We start with part (I), which will require the digital wavelet transform introduced in Subsection \ref{subsec:wavelets}.
Without loss of generality, we assume $M_{c_j} = Id$. The general case can be treated similarly. First, using
\eqref{eq:refine01}, \eqref{eq:refine02} and \eqref{eq:nonsep}, we obtain
\begin{eqnarray} \nonumber
\hat \psi_{j,0,m}(\xi) &=& 2^{-\frac{3}{4}j}e^{-2 \pi i m \cdot A^{-1}_{2^j}\xi}\hat \psi(A_{2^j}^{-1}\xi) \\
&=& 2^{-\frac{3}{4}j-1}e^{-2 \pi i n \cdot A^{-1}_{2^j}\xi}P(A^{-1}_{2^j}Q^{-1}\xi)
\hat g(2^{-j-1}\xi_1)\hat h(2^{-j/2-1}\xi_2)\hat \phi(A^{-1}_{2^{j}}2^{-1}\xi), \label{eq:discrete1}
\end{eqnarray}
where $Q = \text{diag}(2,1)$. From now on, we assume that $\phi$ is an orthonormal scaling function, i.e.,
\begin{equation}\label{eq:ortho}
\sum_{n \in {\mathbb Z}^2} |\hat \phi (\xi + n)|^2 = 1.
\end{equation}
Applying  \eqref{eq:refine01} and \eqref{eq:refine02} iteratively, we obtain
\[
\hat \phi(A^{-1}_{2^{j}}2^{-1}\xi) = 2^{-J+\frac{3}{4}j+1}\prod_{\ell = 1}^{J-j-1}\hat h (2^{-j-1-\ell}\xi_1)\prod_{\ell = 1}^{J-j/2-1}\hat h (2^{-j/2-1-\ell}\xi_2)\hat \phi(2^{-J}\xi).
\]
Inserted in \eqref{eq:discrete1}, it follows that
\[
\hat \psi_{j,0,m}(\xi) = 2^{-J}e^{-2 \pi i m \cdot A^{-1}_{2^j}\xi}{P(A_{2^j}^{-1}Q^{-1}\xi)}\hat g_{J-j}\otimes \hat h_{J-j/2}(2^{-J}\xi)\hat \phi(2^{-J}\xi).
\]
Assuming the function $f$ to be of the form \eqref{eq:code2}, hence its Fourier transform is
\[
\hat f (\xi) = 2^{-J} \hat f_J(2^{-J}\xi)\hat \phi (2^{-J}\xi),
\]
we conclude that
\begin{eqnarray*}
\lefteqn{\langle f,\psi_{j,0,m}\rangle}\\ &=& 2^{-2J}\int_{\mathbb R^2}\hat f_J (2^{-J}\xi)|\hat \phi(2^{-J}\xi)|^2
e^{2 \pi i m \cdot A^{-1}_{2^j}\xi}{P^*(A_{2^j}^{-1}Q^{-1}\xi)}\hat{W}^*_j(2^{-J}\xi)\hat \phi^*(2^{-J}\xi)d\xi
\end{eqnarray*}
where $W_j = g_{J-j}\otimes h_{J-j/2}$. Letting $\eta = 2^{-J}\xi$ and using \eqref{eq:ortho},
\begin{eqnarray*}
\langle f,\psi_{j,0,m}\rangle &=& \int_{\mathbb R^2} \hat f_J(\eta)|\hat \phi(\eta)|^2
 e^{2 \pi i m \cdot A^{-1}_{2^j}2^J\eta}P^*(2^JA_{2^j}^{-1}Q^{-1}\eta)\hat{W}^*_j(\eta)d\eta
 \\
 &=& \int_{[0,1]^2} \hat f_J(\eta) \sum_{n \in \mathbb Z^2}|\hat\phi(\eta+n)|^2
 e^{2 \pi i m \cdot A^{-1}_{2^j}2^J\eta}P^*(2^JA_{2^j}^{-1}Q^{-1}\eta)\hat{W}^*_j(\eta)d\eta \\
 &=& \int_{[0,1]^2} \hat f_J(\eta)
 e^{2 \pi i m \cdot A^{-1}_{2^j}2^J\eta}P^*(2^JA_{2^j}^{-1}Q^{-1}\eta)\hat{W}^*_j(\eta)d\eta.
\end{eqnarray*}
Thus, letting $p_{j}(n)$ be the Fourier coefficients of $P(2^{J-j-1}\xi_1,2^{J-j/2}\xi_2)$ with a 2D fan filter,
\begin{equation}\label{eq:noshear}
\langle f,\psi_{j,0,m}\rangle = (f_J*(\overline{p_j*W_j}))(A_{2^j}^{-1}2^Jm).
\end{equation}
We remark that in case $p_j \equiv 1$ this coincides with the 2D wavelet transform associated with
the anisotropic scaling matrix $A_{2^{j}}$ and a separable wavelet generator, while we obtain the
anisotropic  wavelet transform with a nonseparable wavelet generator in case $p_j$ is a nonseparable filter.
Note that in case of $M_{c_j} \neq Id$, \eqref{eq:noshear} can be easily extended to obtain
\begin{equation}\label{eq:ext_noshear}
\langle f,\psi_{j,0,m}\rangle = (f_J*(\overline{p_j*W_j}))(A_{2^j}^{-1}2^JM_{c_j}m),
\end{equation}
for which the sampling matrix $M_{c_j}$ should be chosen so that $A_{2^j}^{-1}2^JM_{c_j}m \in \mathbb Z^2$.

\begin{figure}[h]
\begin{center}
\includegraphics[width=1.4in,height=1.2in]{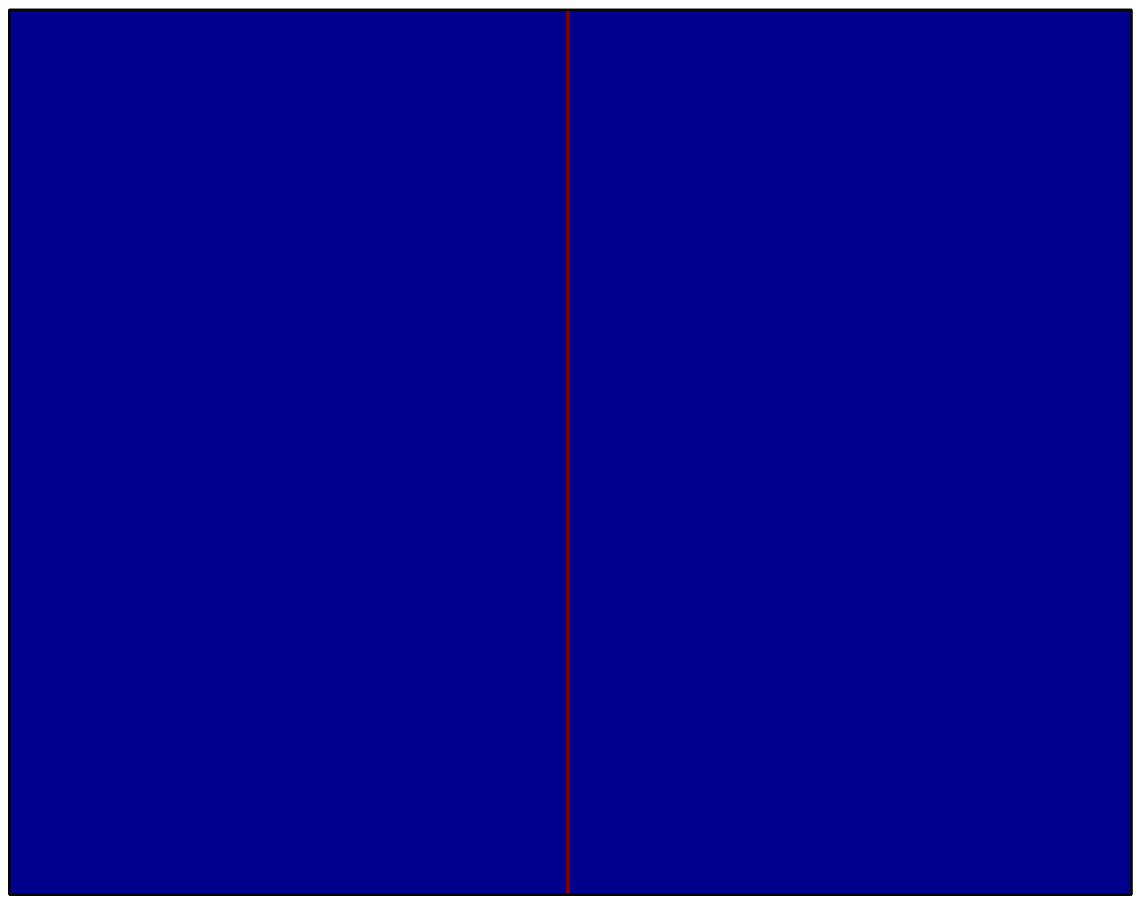}
\includegraphics[width=1.4in,height=1.2in]{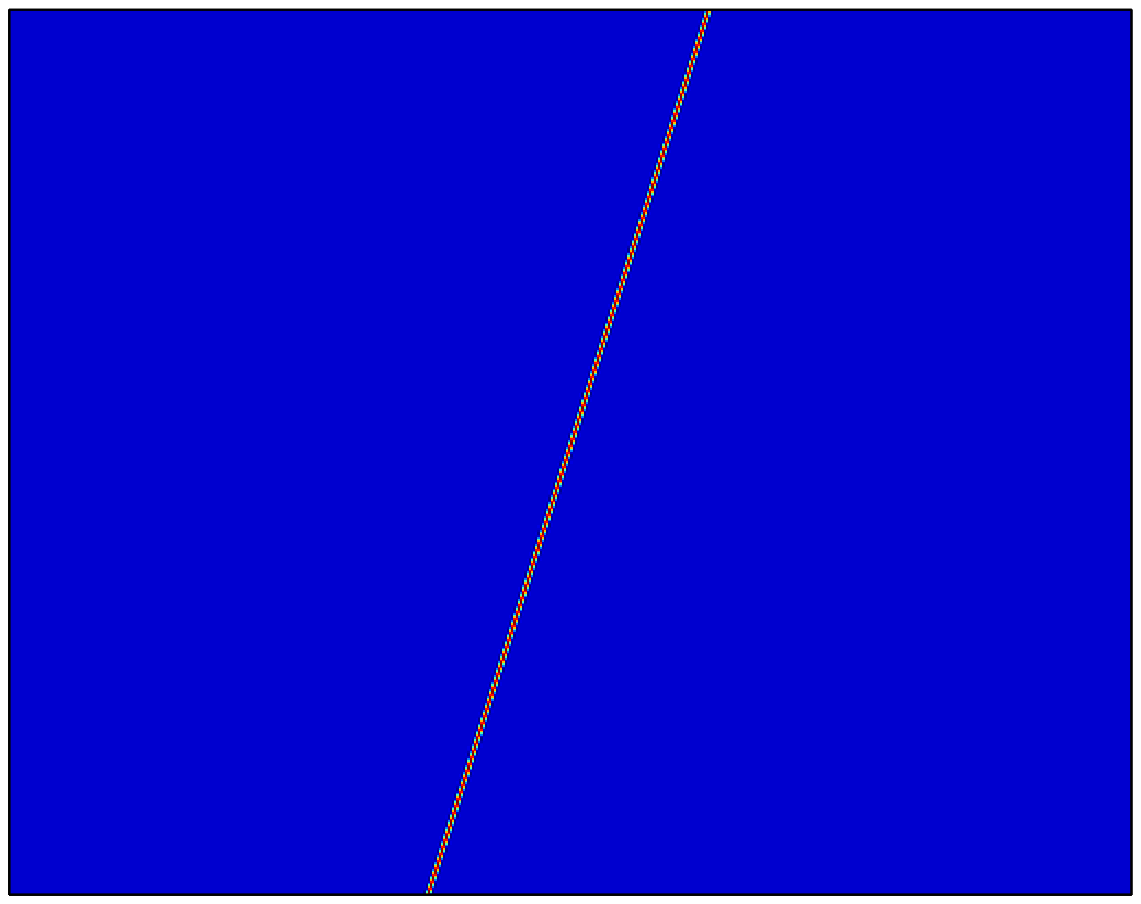}
\includegraphics[width=1.4in,height=1.2in]{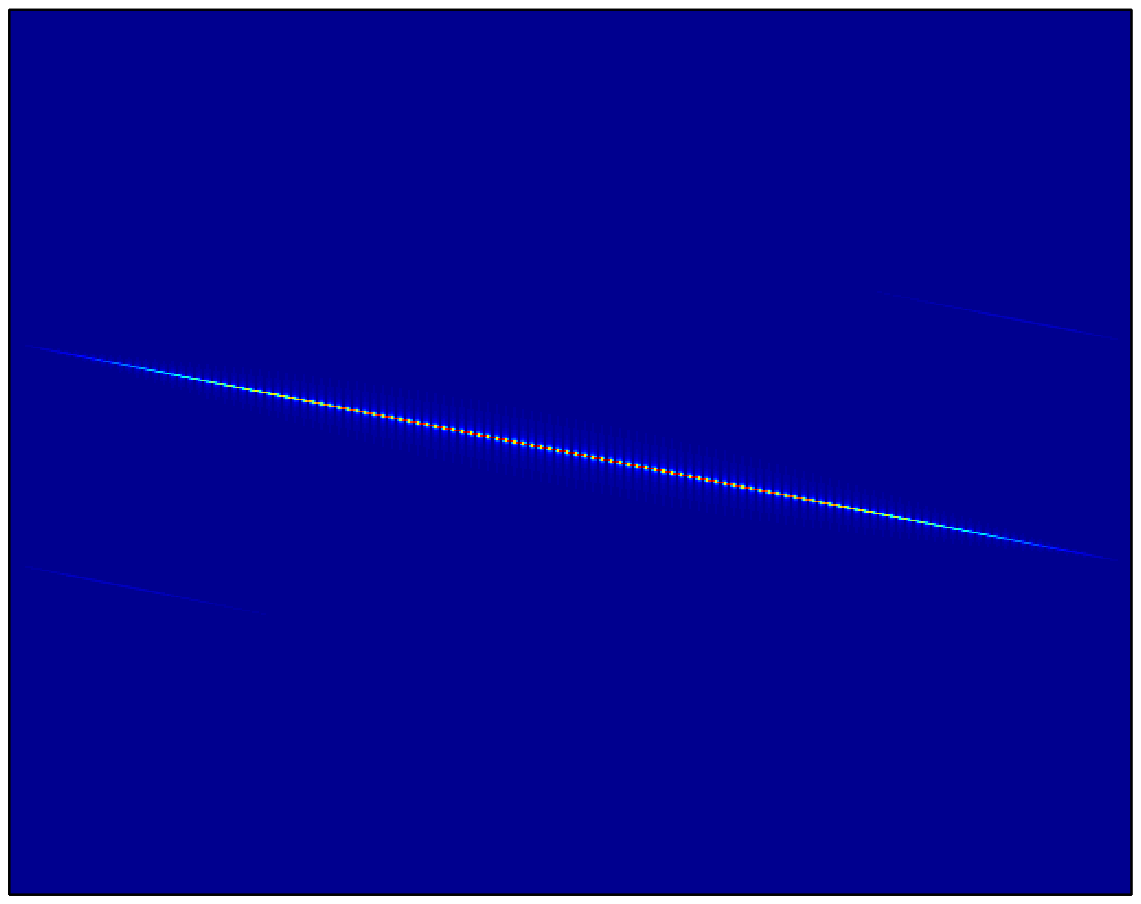}
\put(-255,-10){(a)}
\put(-155,-10){(b)}
\put(-55,-10){(c)}
\caption{(a) $f_d$ : Digital image of $1_{\{x : x_1 = 0\}}$. (b) $S^d_{-1/4}(f_d)$ : Sheared image. (b) Magnitude of the DFT of $S^d_{-1/4}(f_d)$. }\label{fig:shearop}
\end{center}
\end{figure}

We next turn to part (II), i.e., to faithfully digitize the shear operator $S_{2^{-j/2}k}$, which will then provide
an algorithm for computing $\langle f,\psi_{j,k,m} \rangle$ by using \eqref{eq:shear} combined with \eqref{eq:noshear}.
We however face the problem that in general, the shear matrix $S_{2^{-j/2}k}$ does not preserve the regular grid $\mathbb Z^2$, i.e.,
\[
S_{2^{-j/2}k}(\mathbb Z^2) \neq \mathbb Z^2.
\]
One approach to resolve this problem is to refine the regular grid $\mathbb Z^2$ along the horizontal axis $x_1$ by a factor
of $2^{j/2}$. With this modification, the new grid $2^{-j/2}\mathbb Z \times \mathbb Z$ is now invariant under the shear operator
$S_{2^{-j/2}k}$, since
\[
2^{-j/2}\mathbb Z\times \mathbb Z = Q^{-j/2}(\mathbb Z^2) = Q^{-j/2}(S_k(\mathbb Z^2)) =  S_{2^{-j/2}k}(2^{-j/2}\mathbb Z\times \mathbb Z).
\]
Thus the operator $S_{2^{-j/2}k}$ is indeed well defined on the refined grid $2^{-j/2} \mathbb Z \times \mathbb Z$, which
provides a natural discretization of $S_{2^{-j/2}k}$. This observation gives rise to the following strategy for computing
sampling values of $f(S_{k2^{-j/2}}\cdot)$ for given samples $f_J \in \ell^2(\mathbb Z^2 )$ from the function
$f \in L^2(\mathbb R^2)$. For this, let $\uparrow 2^{j/2}$, $\downarrow 2^{j/2}$, and $*_1$ be the 1D upsampling, downsampling,
and convolution operator along the horizontal axis $x_1$, respectively. First, compute interpolated sample values $\tilde{f}_J$
from $f_J(n)$ on the refined  grid $2^{-j/2}\mathbb Z \times \mathbb Z$, which is invariant under $S_{k2^{-j/2}}$, by
\begin{equation}\label{eq:dshear1}
\tilde{f}_J := ((f_J)_{\uparrow 2^{j/2}} *_{1} h_{j/2}).
\end{equation}
Recall that on this new grid $2^{-j/2}\mathbb Z \times \mathbb Z$, the shear operator $S_{k2^{-j/2}}$ becomes $S_k$ with integer
entries. This allows $\tilde{f}_J$ now to be resampled by $S_k$, followed by reversing the previous convolution and upsampling,
i.e.,
\begin{equation}\label{eq:dshear2}
S^d_{2^{-j/2}k}(f_J) := \Bigl(((\tilde{f}_J)(S_k\cdot)) *_1 \overline{h}_{j/2} \Bigr)_{\downarrow 2^{j/2}}.
\end{equation}
which combined with \eqref{eq:dshear1} performs the application of $S_{k2^{-j/2}}$ to discrete data $f_J$. Figure~\ref{fig:shearop}
illustrates how this approach effectively removes otherwise appearing aliasing effect.

Finally, combining \eqref{eq:ext_noshear} with the digital shear operator $S^d_{2^{-j/2}k}$ just defined by \eqref{eq:dshear1}
and \eqref{eq:dshear2}, yields a faithful digital shearlet transform as follows.

\begin{definition}\label{defi:digitSH}
Let $f_J \in \ell^2(\mathbb Z^2)$ be the scaling coefficients given in \eqref{eq:code2}. Then the {\em digital shearlet transform}
associated with $\Psi(\psi; c)$ is defined by
\[
DST^{\text{2D}}_{j,k,m}(f_J) = (\overline{\psi^d_{j,k}}*f_J)(2^JA^{-1}_{2^j}M_{c_j}m) \quad \mbox{for } j = 0,\dots,J-1,
\]
where
\[
\psi^d_{j,k} = S^d_{k/2^{j/2}}(p_j * W_j),
\]
with the shearing operator defined by \eqref{eq:dshear1} and \eqref{eq:dshear2},
and the sampling matrix $M_{c_j}$ chosen so that $2^JA^{-1}_{2^j}M_{c_j}m \in \mathbb Z^2$.
\end{definition}

\subsubsection{The General Case} \label{subsubsec:nonparabolic}

The digitalization of the shearlet transform associated with the special case of a classical cone-adapted discrete shearlet system
as defined in Definition \ref{defi:digitSH} shall now be extended to  universal shearlet systems, where
the parabolic scaling matrices are generalized to $A_{\alpha_j, 2^j} = \text{diag}(2^j, 2^{\alpha_j j/2})$ and
$\tilde{A}_{\alpha_j, 2^j} = \text{diag}(2^{\alpha_j j/2},2^j)$ with $\alpha_j \in(0,2)$ (cf. Subsection \ref{subsec:discreteshearlets}).
In this general setting, the range of the shearing parameter $k\in\Z$ is given as $\abs{k}\leq \lceil2^{(2-\alpha_j)j/2}\rceil$
in each cone for each scale $j\in\N_0$. The associated shearlet coefficients are then given by
\begin{equation}\label{eq:alpha_shear}
\langle f,\psi_{j,k,m} \rangle = \langle f(S_{k/2^{(2-\alpha_j)j/2}}\cdot),\psi_{j,0,m}(\cdot)\rangle,
\end{equation}
where $\psi_{j,0,m}(\cdot) = 2^{\frac{(2+\alpha_j)j}{4}}\psi(A_{\alpha_j,2^j}\cdot-m)$.

From now on, we retain notations from the previous section, otherwise we specify them. To digitalize \eqref{eq:alpha_shear},
we slightly change the range of shearings $\abs{k}\leq \lceil2^{(2-\alpha_j)j/2}\rceil$ to  $|k| \leq 2^{\lceil{(2-\alpha_j)j/2}\rceil}$,
the reason being that then the number of shearings is determined by dyadic scales. Thus, also to have an integer scaling matrix,
we now consider the slightly modified version of \eqref{eq:alpha_shear} given by
\begin{equation}\label{eq:adshear}
\langle f(S_{k/2^{d_j}}\cdot),\psi_{j,0,m}(\cdot)\rangle \quad \text{with } d_j := \lceil (2-\alpha_j)j/2 \rceil,
\end{equation}
where
\[
\psi_{j,0,m}(\cdot) = 2^{\frac{2j-d_j}{2}}\psi(A_{d_j}\cdot - m) \quad \mbox{and} \quad A_{d_j} := \text{diag}(2^j,2^{j-d_j}).
\]

Again we have to digitalize parts (I) and (II). Starting with part (I), the coefficients $\langle f,\psi_{j,0,m}\rangle$
can be digitalized similar to \eqref{eq:ext_noshear}, except for changing $A_{2^j}$ to $A_{d_j}$, i.e., changing the scaling
parameter $j/2$ to $j-d_j$. Then the resulting discretization for $\psi_{j,0,m}$ essentially follows the $\alpha_j$ scaling
operator $A_{\alpha_j, 2^j}$ in the sense that $j-d_j \approx \alpha_j (j/2)$ for sufficiently large $j$. In short, we obtain
\begin{equation}\label{eq:ext_noshear_alpha}
\langle f,\psi_{j,0,m}\rangle = (f_J*(\overline{p_j*W_j}))(A_{d_j}^{-1}2^JM_{c_j}m),
\end{equation}
for which we modify $p_j$ and $W_j$ so that $p_{j}(n)$ are now the Fourier coefficients of $P(2^{J-j-1}\xi_1,2^{J-(j-d_j)}\xi_2)$
with a 2D fan filter $P$ and $W_j = g_{J-j}\otimes h_{J-(j-d_j)}$ and the sampling matrix $M_{c_j}$ chosen so that
$A_{d_j}^{-1}2^JM_{c_j}m \in \mathbb Z^2$.

Concerning part (II), we just need to slightly modify \eqref{eq:dshear1} and \eqref{eq:dshear2} to become
\begin{equation}\label{eq:dshear2_alpha}
S^d_{k/2^{d_j}}(f_J) := \Bigl(((\tilde{f}_J)(S_k\cdot)) *_1 \overline{h}_{d_j} \Bigr)_{\downarrow 2^{d_j}}
\end{equation}
where
\begin{equation}\label{eq:dshear1_alpha}
\tilde{f}_J = ((f_J)_{\uparrow 2^{d_j}} *_{1} h_{d_j})
\end{equation}
and filter coefficients $h_{d_j}$ are defined as in \eqref{eq:lowhigh}.

Thus, concluding, we derive the following digitalization of the universal shearlet transform associated with
the shearlets from $\Psi(\psi; \alpha, c)$.

\begin{definition}\label{defi:digitSH_alpha}
Let $f_J \in \ell^2(\mathbb Z^2)$ be the scaling coefficients given in \eqref{eq:code2}. Then the {\em digital shearlet transform}
associated with $\Psi(\psi; \alpha, c)$ is defined by
\[
DST^{\text{2D}}_{j,k,m}(f_J) = (\overline{\psi^d_{j,k}}*f_J)(2^JA^{-1}_{d_j}M_{c_j}m) \quad \mbox{for } j = 0,\dots,J-1,
\]
where
\begin{equation}\label{eq:2dshear_filters}
\psi^d_{j,k} = S^d_{k/2^{j/2}}(p_j * W_j)
\end{equation}
with the shearing operator defined by \eqref{eq:dshear1_alpha} and \eqref{eq:dshear2_alpha}, $p_j$ and $W_j$ modified
as explained above, and the sampling matrix $M_{c_j}$ chosen so that $2^JA^{-1}_{d_j}M_{c_j}m \in \mathbb Z^2$.
\end{definition}

\subsection{3D Shearlet Transform} \label{subsec:3Dtrafo}

Following the approach of the 2D case, we choose the 3D shearlet generator $\psi$ by
\begin{equation} \label{eq:generator3D}
\hat \psi(\xi) = \Bigl(P\bigl({\xi_1}/{2},\xi_2\bigr)\hat\psi_1(\xi_1)\hat \phi_1(\xi_2)\Bigr)\cdot\Bigl( P\bigl({\xi_1}/{2},\xi_3\bigr)\hat \phi_1(\xi_3) \Bigr)
\end{equation}
so that -- as it was also the idea in the 2D case -- $\hat \psi$ well approximates the 3D band-limited shearlet generator
in the sense that
\[
P\bigl({\xi_1}/{2},\xi_2\bigr)\hat\psi_1(\xi_1)\hat \phi_1(\xi_2) \approx
\hat \psi_1(\xi_1) \hat \phi_1(\xi_2/\xi_1) \quad \text{and} \quad
P\bigl({\xi_1}/{2},\xi_3\bigr)\hat \phi_1(\xi_3) \approx
\hat \phi_1(\xi_3/\xi_1).
\]
Similar to the 2D case, we may choose compactly supported functions $\psi_1$ and $\phi_1$, and a finite 2D fan filter satisfying \eqref{eq:cond1} and \eqref{eq:cond2}, respectively.
Then we have
\[
|\hat \phi(\xi)|^2 + \sum_{j \ge 0}\sum_{|k| \le \lceil 2^{ j/2 } \rceil}
\Bigl(|\hat \psi_{j,k,0}(\xi)|^2 + |\hat{\tilde{\psi}}_{j,k,0}(\xi)|^2 + |\hat{\breve{\psi}}_{j,k,0}(\xi)|^2\Bigr)  > \min(\delta_2^3,\delta_1\delta^2_2\delta^2_3),
\]
where $|k| = \max (|k_1|,|k_2|)$, $\phi = \phi_1 \otimes \phi_1 \otimes \phi_1$, $\tilde{\psi}_{j,k,m}(x_1,x_2,x_3) = \psi_{j,k,m}(x_2,x_1,x_3)$ and
$\breve{\psi}_{j,k,m}(x_1,x_2,x_3) = \psi_{j,k,m}(x_3,x_1,x_2)$.
This inequality provides a lower frame bound provided that $\psi_1$ and $\phi_1$ decay sufficiently fast in frequency and $\psi_1$ has sufficient vanishing moments. Also, an upper frame bound can be obtained from the fast decay rate of $\hat \psi$  generated by those 1D functions $\psi_1$ and $\phi_1$. We refer to \cite{KLL2012} for more details.

We next discuss a digitalization of the associated shearlet coefficients $\langle f,\psi_{j,k,m} \rangle$ again only from $\Psi(\psi; \alpha, c)$,
where $f$ is given by
\begin{equation} \label{eq:code3D}
f(x) = \sum_{n \in \mathbb Z^3} f_J(n)2^{J\cdot 3/2}\phi_1\otimes\phi_1\otimes \phi_1(2^Jx-n).
\end{equation}

\subsubsection{The Parabolic Case}

We start with the parabolic case, i.e., $\alpha_j = 1$ for all scales $j$. Recalling the 3D parabolic scaling matrix
$A_{2^j} = \text{diag}(2^j,2^{j/2},2^{j/2})$ from Subsection~\ref{subsec:discreteshearlets}, we now consider
\[
\hat{\Phi}_{j,k_1}(\xi_1,\xi_2) = (\hat{\phi}_1\cdot P)\Bigl(Q^{-1}\Bigl(S^{T}_{-k_1}A^{-1}_{2^{j}}(\xi_1,\xi_2)^{T}\Bigr)\Bigr),
\]
\[
\hat{\Phi}_{j,k_2}(\xi_1,\xi_3) = (\hat{\phi}_1\cdot P)\Bigl(Q^{-1}\Bigl(S^{T}_{-k_2}A^{-1}_{2^{j}}(\xi_1,\xi_3)^{T}\Bigr)\Bigr),
\]
and
\[
\gamma_{j,k,m}(\xi) = \bigl(\frac{m_1}{2^{j}}-k_1\frac{m_2}{2^{j}}-k_2\frac{m_3}{2^{j}}\bigr)\xi_1+m_2 \frac{\xi_2}{2^{j/2}} + m_3  \frac{\xi_3}{2^{j/2}}.
\]
Using the generator \eqref{eq:generator3D}, our 3D shearlets $\psi_{j,k,m}$ are then defined by
\begin{equation} \label{eq:3Dshearlet}
\hat{\psi}_{j,k,m}(\xi) = \frac{1}{2^{j}}\hat{\psi}_1\bigl({\xi_1}/{2^{j}}\bigr)\hat{\Phi}_{j,k_1}(\xi_1,\xi_2)\hat{\Phi}_{j,k_2}
(\xi_1,\xi_3)e^{-2\pi i \gamma_{j,k,m}(\xi)}.
\end{equation}
Since $\Phi_{j,k_1}$ and $\Phi_{j,k_2}$ are functions of the form of 2D shearlets, they might be discretized similar
as in Definition \ref{defi:digitSH} -- though omitting the convolution with the high-pass filter $g_{J-j}$ --, which gives
($*_{x_i}$ denoting 1D convolution along the $x_i$ axis)
\[
\Phi^{\text{d}}_{j,k_1}(n_1,n_2) = \Bigl(S^{{d}}_{k_12^{-j/2}}(h_{J-j/2}*_{x_2}p_j)\Bigr)(n_1,n_2)
\]
and
\[
\Phi^{\text{d}}_{j,k_2}(n_1,n_3) = \Bigl(S^{{d}}_{k_22^{-j/2}}(h_{J-j/2}*_{x_3}p_j)\Bigr)(n_1,n_3).
\]
Finally, by \eqref{eq:refine02}, a 1D wavelet $2^{j/2}\psi_1(2^j\cdot)$ can be digitalized by the 1D filter $g_{J-j}$.
This gives rise to 3D digital shearlet filters $\psi^{\text{d}}_{j,k}$ specified in the definition below, which
discretize $\psi_{j,k,m}$ from \eqref{eq:3Dshearlet}. Summarizing, our digitalization of the shearlet transform associated
with the shearlets from $\Psi(\psi; c)$ (i.e., the parabolic case) is defined as follows.

\begin{definition}\label{defi:3DdigitSH}
Let $f_J \in \ell^2(\mathbb Z^3)$ be the scaling coefficients given in \eqref{eq:code3D}, and retain the definitions and
notions of this subsection. Then the {\em digital shearlet transform} associated with $\Psi(\psi; c)$ is defined by
\[
DST^{\text{3D}}_{j,k}(f_J)(m) = (f_J*\overline{\psi}^{\text{d}}_{j,k})(\tilde{m}) \quad \mbox{for } j = 0,\dots,J-1,
\]
where the 3D digital shearlet filters $\psi^{\text{d}}_{j,k}$ are defined by
\[
\hat{\psi}^{\text{d}}_{j,k}(\xi) =  \hat{g}_{J-j}(\xi_1)\hat{\Phi}^{\text{d}}_{j,k_1}(\xi_1,\xi_2)\hat{\Phi}^{\text{d}}_{j,k_2}(\xi_1,\xi_3)
\]
and
\[
\tilde{m} = (2^{J-j}c^j_1m_1,2^{J-j/2}c^j_2m_2,2^{J-j/2}c^j_3m_3)
\]
with the sampling constants $c^j_1, c^j_2$ and $c^j_3$ chosen so that $\tilde{m} \in \mathbb Z^3$.
\end{definition}

\begin{figure}[h]
\begin{center}
\includegraphics[width=1.4in,height=1.2in]{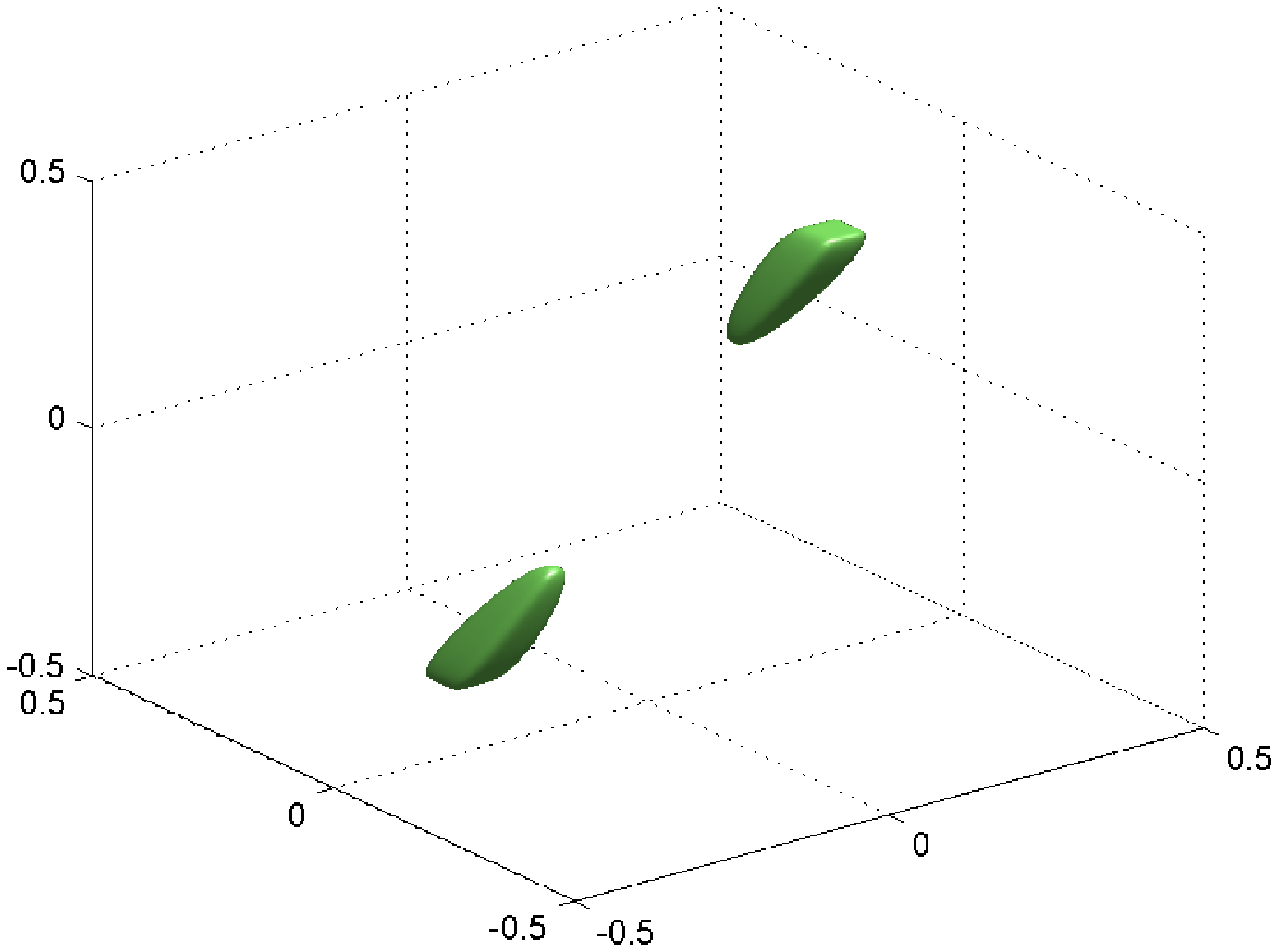}
\includegraphics[width=1.4in,height=1.2in]{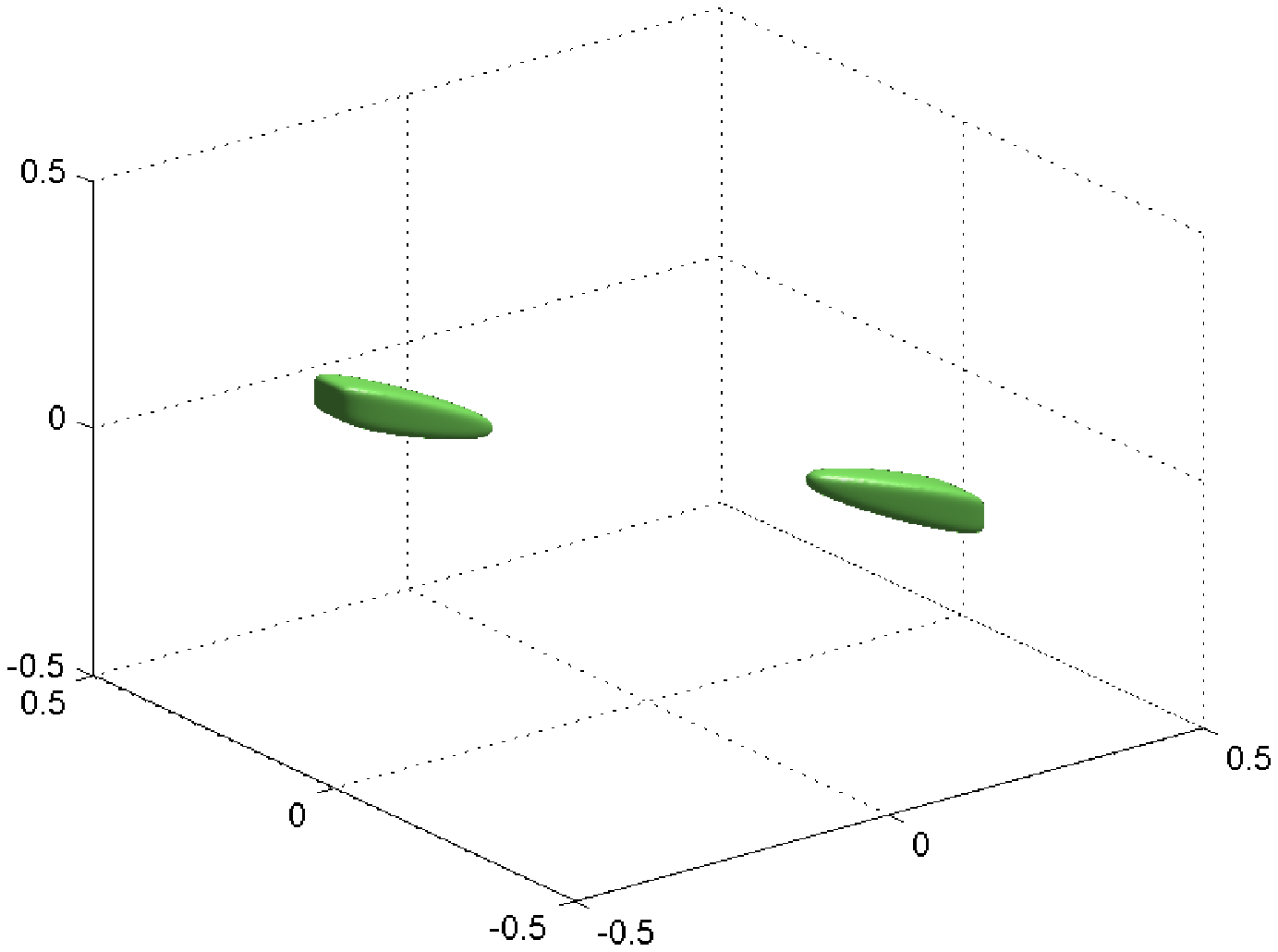}
\includegraphics[width=1.4in,height=1.2in]{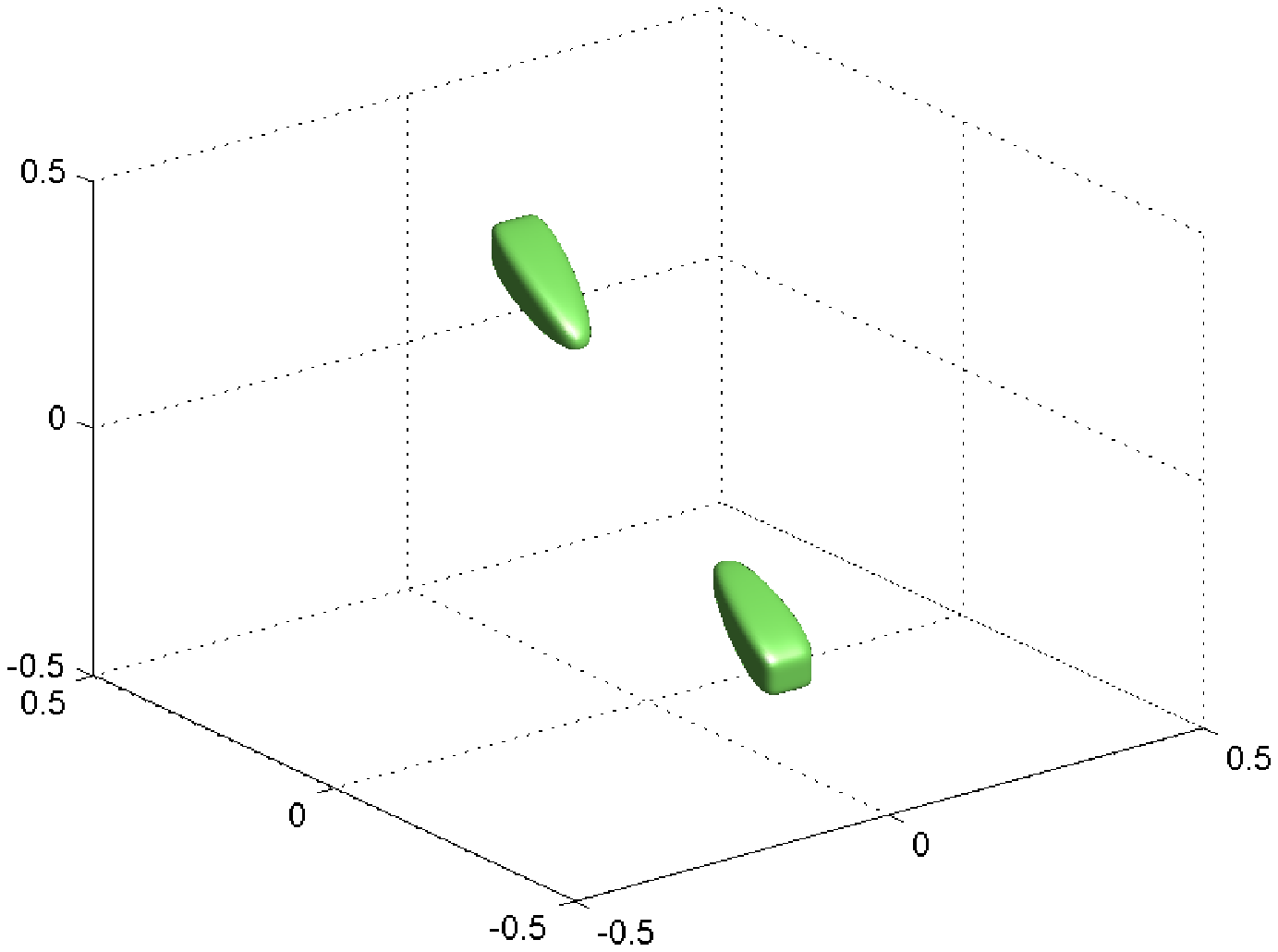}
\put(-255,-10){(a)}
\put(-155,-10){(b)}
\put(-55,-10){(c)}
\end{center}
\caption{3D digital shearlets in the frequency domain:  $\hat{\psi}^{\text{d}}_{j,k}$}\label{fig:3d_shearlets_frequency}
\end{figure}

The chosen 3D digital shearlet filters are illustrated in Figure~\ref{fig:3d_shearlets_frequency}.

\subsubsection{The General Case}

We now describe the digitalization of the 3D universal shearlet transform. Similar to the 2D situation, we first slightly
modify the 3D scaling matrix $A_{\alpha_j, 2^j}$ to consider integer matrices $A_{d_j} = \text{diag}(2^j,2^{j-d_j},2^{j-d_j})$
with $d_j = \lceil (2-\alpha_j)j/2\rceil$ for each scale $j \in \N_0$. Also the shearing parameters $k_1,k_2\in\Z$ range from
$-2^{\lceil{(2-\alpha_j)j/2}\rceil}$ to $2^{\lceil{(2-\alpha_j)j/2}\rceil}$. Again following the 2D approach, the 3D digital shearlet
filters $\psi^d_{j,k}$ defined in Definition \ref{defi:3DdigitSH} are generalized by modifying $\Phi^{\text{d}}_{j,k_1}$ and
$\Phi^{\text{d}}_{j,k_2}$ to be
\begin{equation} \label{eq:3Dfilters1_alpha}
\Phi^{\text{d}}_{j,k_1}(n_1,n_2) = \Bigl(S^{{d}}_{k_12^{-d_j}}(h_{J-(j-d_j)}*_1p_j)\Bigr)(n_1,n_2)
\end{equation}
and
\begin{equation} \label{eq:3Dfilters2_alpha}
\Phi^{\text{d}}_{j,k_2}(n_1,n_3) = \Bigl(S^{{d}}_{k_22^{-d_j}}(h_{J-(j-d_j)}*_1p_j)\Bigr)(n_1,n_3),
\end{equation}
where $p_{j}(n)$ are the Fourier coefficients of $P(2^{J-j-1}\xi_1,2^{J-(j-d_j)}\xi_2)$ with a 2D fan filter $P$, the 1D filter
$h_{J-(j-d_j)}$ is defined as in \eqref{eq:lowhigh}, and $S^{{d}}_{k_22^{-d_j}}$ is the discrete shear operator defined in
\eqref{eq:dshear2_alpha}. Thus, the 3D digital shearlet transform associated with universal shearlets is defined as described
in the following definition.

\begin{definition}\label{defi:3DdigitSH_alpha}
Let $f_J \in \ell^2(\mathbb Z^3)$ be the scaling coefficients given in \eqref{eq:code3D}, and retain the definitions and
notions of this subsection. Then the {\em digital shearlet transform} associated with $\Psi(\psi; \alpha, c)$ is defined by
\[
DST^{\text{3D}}_{j,k}(f_J)(m) = (f_J*\overline{\psi}^{\text{d}}_{j,k})(\tilde{m}) \quad \mbox{for } j = 0,\dots,J-1,
\]
where the 3D digital shearlet filters $\psi^{\text{d}}_{j,k}$ are defined using \eqref{eq:3Dfilters1_alpha} and \eqref{eq:3Dfilters2_alpha} by
\[
\hat{\psi}^{\text{d}}_{j,k}(\xi) =  \hat{g}_{J-j}(\xi_1)\hat{\Phi}^{\text{d}}_{j,k_1}(\xi_1,\xi_2)\hat{\Phi}^{\text{d}}_{j,k_2}(\xi_1,\xi_3)
\]
and
\[
\tilde{m} = (2^{J-j}c^j_1m_1,2^{J-j/2}c^j_2m_2,2^{J-j/2}c^j_3m_3)
\]
with the sampling constants $c^j_1, c^j_2$ and $c^j_3$ chosen so that $\tilde{m} \in \mathbb Z^3$.
\end{definition}

\subsection{Inverse Shearlet Transform}\label{subsec:inverse}

In this section, we define an inverse digital shearlet transform, which provides a stable reconstruction of
$f_J \in \ell^2(\mathbb Z^d)$ from the shearlet coefficients obtained by the digital shearlet transforms
from Subsections \ref{subsec:2Dtrafo} and \ref{subsec:3Dtrafo}. For this, we will consider only the 2D parabolic
case retaining notations of Subsection \ref{subsubsec:parabolic}, since this can be extended to the general 2D
case as well as the 3D case in a straightforward manner.

We first observe that in general, the forward shearlet transform defined in Definition~\ref{defi:digitSH_alpha}
can be inverted by a frame reconstruction algorithm based on the conjugate gradient method due to frame property
of shearlets -- see \cite{Mal08}.

It seems impossible to obtain a direct reconstruction formula unless we skip subsampling, which would then lead to a highly redundant transform.
However, one possibility was indeed recently discovered in \cite{Lim2013a} in the parabolic situation, which we now describe.
In this approach, we first set $c^j_1 = 2^{j-J}$ and $c^j_2 = 2^{j/2-J}$ in Definition~\ref{defi:digitSH}.
In this situation, the digital shearlet transform is merely a 2D convolution with shearlet filters, yielding a shift-invariant
linear transform. Hence, for $f_J \in \ell^2(\Z^2)$, the digital shearlet transform takes the form
\begin{equation}\label{eq:convolution}
DST^{\text{2D}}_{j,k,m}(f_J) = f_J *\overline{\psi}^{\text{d}}_{j,k}(m).
\end{equation}
As indicated before, the digital shearlet filters $\tilde{\psi}^{\text{d}}_{j,k}$ corresponding to $\tilde{\psi}_{j,k,m}$
are derived by switching the order of variables, which implies that the same convolution formula as \eqref{eq:convolution}
holds for the shearlet transform associated with $\tilde{\Psi}(\tilde{\psi}; c)$, and we can define
\begin{equation}\label{eq:convolution01}
{\widetilde{DST^{\text{2D}}_{j,k,m}}}(f_J)  = f_J *\overline{\tilde{\psi}^{\text{d}}}_{j,k}(m).
\end{equation}
Let us now select separable low-pass filter by
\[
\hat{\phi}^{\text{d}}(\xi) = \hat{h}_{J}(\xi_1)\cdot\hat{h}_{J}(\xi_2),
\]
set
\begin{equation}\label{eq:frame_filter}
\hat \Psi^{\text{d}}(\xi) := |\hat{\phi}^{\text{d}}(\xi)|^2+\sum_{j=0}^{J-1}\sum_{|k| \leq 2^{\lceil j/2 \rceil}} \Bigl(|\hat{\psi}^{\text{d}}_{j,k}(\xi)|^2
+|\hat{\tilde{\psi}}^{\text{d}}_{j,k}(\xi)|^2\Bigr),
\end{equation}
and also define dual shearlet filters by
\begin{equation}\label{eq:dual_filter}
\hat{\varphi}^{\text{d}}(\xi) = \frac{\hat{\phi}^{\text{d}}(\xi)}{\hat \Psi^{\text{d}}(\xi)}, \,\, \hat{\gamma}^{\text{d}}_{j,k}(\xi) =
\frac{\hat{\psi}^{\text{d}}_{j,k}(\xi)}{\hat \Psi^{\text{d}}(\xi)}
,\,\,\hat{\tilde{\gamma}}^{\text{d}}_{j,k}(\xi) = \frac{\hat{\tilde{\psi}}^{\text{d}}_{j,k}(\xi)}{\hat \Psi^{\text{d}}(\xi)}.
\end{equation}
Using \eqref{eq:convolution} and \eqref{eq:convolution01}, the Fourier transform of $f_J$ can be written as
\begin{eqnarray*}
\hat{f}_J &=& \frac{\hat{f}_J}{\hat \Psi^{\text{d}}}\Bigl( |\hat{\phi}^{\text{d}}|^2+\sum_{j=0}^{J-1}\sum_{|k| \leq 2^{\lceil j/2 \rceil}} \Bigl(|\hat{\psi}^{\text{d}}_{j,k}|^2
+|\hat{\tilde{\psi}}^{\text{d}}_{j,k}|^2\Bigr) \Bigr)   \\
&=& \hat{f}_J \cdot (\hat{\phi}^{\text{d}})^*\frac{\hat{\phi}^{\text{d}}}{\hat \Psi^{\text{d}}} +\sum_{j=0}^{J-1}\sum_{|k| \leq 2^{\lceil j/2 \rceil}} \Bigl(
\hat{f}_J\cdot(\hat{\psi}^{\text{d}}_{j,k})^*\frac{\hat{\psi}^{\text{d}}_{j,k}}{\hat \Psi^{\text{d}}}
+\hat{f}_J\cdot(\hat{\tilde{\psi}}^{\text{d}}_{j,k})^*\frac{\hat{\tilde{\psi}}^{\text{d}}_{j,k}}{\hat \Psi^{\text{d}}}
\Bigr) \\
&=& \hat{f}_J \cdot (\hat{\phi}^{\text{d}})^*\hat{\varphi}^{\text{d}} + \sum_{j=0}^{J-1}\sum_{|k| \leq 2^{\lceil j/2 \rceil}}
\mathcal{F}({{DST^{\text{2D}}_{j,k,m}}}(f_J))\hat{\gamma}^{\text{d}}_{j,k} +   \mathcal{F}({\widetilde{DST^{\text{2D}}_{j,k,m}}}(f_J))\hat{\tilde{\gamma}}^{\text{d}}_{j,k},
\end{eqnarray*}
where $\mathcal{F} : \ell^2(\Z^2) \rightarrow L^2([0,1]^2)$ is defined as the (discrete time) Fourier transform and $*$ as a superscript denotes the complex conjugate.
These considerations then yield the reconstruction formula given by
\begin{eqnarray}\label{eq:reconstruction}
f_J &=& (f_J*\overline{\phi}^{\text{d}})*{\varphi^{\text{d}}}+\sum_{j=0}^{J-1}\sum_{|k| \leq 2^{\lceil j/2 \rceil}} ({{DST^{\text{2D}}_{j,k,m}}}(f_J))*{\gamma^{\text{d}}}_{j,k} \nonumber
\\ &+&
\sum_{j=0}^{J-1}\sum_{|k| \leq 2^{\lceil j/2 \rceil}}
({\widetilde{DST^{\text{2D}}_{j,k,m}}}(f_J))*{\tilde\gamma^{\text{d}}}_{j,k}.
\end{eqnarray}

Now turning to the more general situation of universal shearlet systems, it can be easily observed that the definition of the dual shearlet filters from
\eqref{eq:dual_filter} can be extended using the generalized digital shearlet filters $\psi^d_{j,k}$ from Definition \ref{defi:digitSH_alpha}. In addition,
the definition of $\hat \Psi^{\text{d}}$ in \eqref{eq:frame_filter} then needs to be generalized to
\[
\hat \Psi^{\text{d}}(\xi) =: |\hat{\phi}^{\text{d}}(\xi)|^2+\sum_{j=0}^{J-1}\sum_{|k| \leq 2^{d_j}} \Bigl(|\hat{\psi}^{\text{d}}_{j,k}(\xi)|^2
+|\hat{\tilde{\psi}}^{\text{d}}_{j,k}(\xi)|^2\Bigr)
\]
with the generalized digital shearlet filters $\psi^d_{j,k}$ in this case.

\section{Implementation: ShearLab 3D} \label{sec:implementation}

An implementation of the digital transforms
\begin{itemize}
\item[$\bullet$] 2D Digital Shearlet Transform,
\item[$\bullet$] 3D Digital Shearlet Transform,
\item[$\bullet$] Forward 2D Digital Shearlet Transform, 
\item[$\bullet$] Inverse 3D Digital Shearlet Transform,
\end{itemize}
described in Section \ref{sec:digitransform} is provided in the MATLAB toolbox ShearLab 3D, which can be downloaded from
\url{www.shearlab.org}. ShearLab 3D requires the Signal Processing Toolbox and the Image Processing Toolbox of MATLAB.
If additionally the Parallel Computing Toolbox is available, CUDA-compatible NVidia graphics cards can be used to gain
a significant speed up.

The ShearLab 3D toolbox provides codes to compute the digital shearlet transform of arbitrarily sized two- and three-dimensional signals according
to the formulas in Definitions \ref{defi:digitSH_alpha} and \ref{defi:3DdigitSH} as well as the inverse shearlet transform \eqref{eq:reconstruction}.
Applying the convolution theorem, these formulas can be computed by multiplying conjugated digital shearlet filters $\overline{\psi^d_{j,k}}$,
their duals $\gamma^d_{j,k}$, and the given signal $f_J$ in the frequency domain.

We now provide more details on the forward (Subsection \ref{subsec:forwardtrafo}) and inverse transform (Subsection \ref{subsec:inversetrafo}),
provide a brief example for a potential use case (Subsection \ref{subsec:example}), discuss the possibility to avoid numerical instabilities
along the seam lines (Subsection \ref{subsec:boundary}), and compute the complexity of our algorithms in Subsection \ref{subsec:complexity}.

\subsection{Forward Transform} \label{subsec:forwardtrafo}

A schematic descriptions of the forward transform is given in Algorithm~\ref{alg:forward_transform}.

 \begin{algorithm}[h]
 \SetAlgoLined
 \KwIn{A signal $\textit{f} \in \R^{X\times Y\times Z}$ , the number of scales $\textit{nScales}\in\N$, a vector $\textit{shearLevels}\in
 \N^{\textit{nScales}}$ specifying the number of differently sheared filters on each scale, a matrix $directionalFilter$ specifying the
 directional filter $P$ (compare equation \eqref{eq:2dshear_filters}) and a vector $quadratureMirrorFilter$ describing the lowpass filter $h_1$
 of a quadrature mirror filter pair.}
 \KwOut{Coefficients $\textit{shearletCoeffs} \in \R^{X\times Y\times Z\times R}$ where $R$ denotes the redundancy of the applied shearlet system.}

 \tcp{Compute shearlet filters in the frequency domain of size $X\times Y\times Z$ according to the parameters \textit{nScales} and \textit{shearLevels}.}

  $\textit{shearletFilters}$ := computeShearletFilters($X,Y,Z,\textit{nScales},\textit{shearLevels},\textit{directionalFilters},\textit{quadratureMirrorFilter}$)\;
  \hfill

  \tcp{Compute frequency representation of the input signal.}
  $\textit{f}_{\textit{freq}}$ := FFT(\textit{f})\;
  \hfill

  \tcp{For each shearlet filter, compute a three-dimensional vector (or two-dimensional, for 2D input) of shearlet coefficients of size
  $X\times Y\times Z$ as the convolution of the time-domain representation of the shearlet filter with the input signal $\textit{f}$.
  According to the convolution theorem, this can be done by pointwise multiplication (denoted by .*) in the frequency domain.}
  \For{\textit{i} := 1 \KwTo R}{
  \textit{shearletCoeffs(i)} := IFFT$\left(\textit{shearletFilters(i)}\text{.*}\overline{\textit{f}_{\textit{freq}}}\right)$\;
  }
 \caption{ShearLab 3D forward transform}
 \label{alg:forward_transform}
\end{algorithm}

As can be seen from Algorithm~\ref{alg:forward_transform}, the computation of a shearlet decomposition of a 2D signal $f\in \ell_2(\Z^2)$
with ShearLab 3D requires the following input parameters:

\begin{itemize}
\item[$\bullet$] $\textit{nScales}$: The number of scales of the shearlet system associated with the desired decomposition.
Each scale corresponds to a ring-like passband in the frequency plane that is constructed from the quadrature mirror filter
pair defined via the lowpass filter $\textit{quadratureMirrorFilter}$. These frequency bands can then be further partitioned
into directionally sensitive elements using a 2D directional filter. Note that $nScales$ can also be viewed as the upper
bound of the parameter $j$ in Definitions \ref{defi:digitSH_alpha} and \ref{defi:3DdigitSH}. Naturally,
increasing the number of scales significantly increases the redundancy of the corresponding shearlet system.\\[-1.8ex]
\item[$\bullet$] $\textit{shearLevels}$: A vector of size  $\textit{nScales}$, specifying for each scale the fineness of the
partitioning of the corresponding ring-like passband. The larger the shear level at a specific scale, the more differently
sheared atoms will live on this scale with increasingly smaller essential support sizes in the frequency domain. To be
precise, let $d_j$ be the $j$-th component in $\textit{shearLevels}$. Then, in the 2D case this choice will generate the
shearlet filters $\psi^{\text{d}}_{j,k}$ defined in Definition \ref{defi:digitSH_alpha} and associated with the scaling
matrix $A_{d_j}$. With this choice, the range of shearing parameters $k$ is given by
$|k| \leq 2^{d_j}$ for each cone, which generates $2(2\cdot2^{d_j}+1)$ shearlet filters for each scale $j$.\\[-1.8ex]
\item[$\bullet$] $\textit{directionalFilter}$: A 2D directional filter that is used to partition the passbands of the
several scales and hence serves as the basis of the directional 'component' of the shearlets. Our default choice in
ShearLab 3D is a maximally flat 2D fan filter (see \cite{dCZD2006} and Figure \ref{fig:filters}). Other directional
filters can for instance be constructed using the $\textit{dfilters}$ method from the Nonsubsampled Contourlet Toolbox.
This value corresponds to the trigonometric polynomial $P$ from equation \eqref{eq:nonsep}.\\[-1.8ex]
\item[$\bullet$] $\textit{quadratureMirrorFilter}$: A 1D lowpass filter defining a quadrature mirror filter pair with
the corresponding highpass filter and thereby a wavelet multiresolution analysis. These filters induce the passbands
associated with the several scales, the number of which is defined by the parameter $\textit{nScales}$. The default
choice in ShearLab 3D is a symmetric maximally flat 9-tap lowpass filter (for an extensive discussion of the properties
of this filter, see Subsection \ref{subsec:parameters}); but basically, any lowpass filter can be used here. This parameter
corresponds to $h_1$ in equation \eqref{eq:lowhigh}.
\end{itemize}

Given the parameters $\textit{nScales}$, $\textit{shearLevels}$, $\textit{quadratureMirrorFilter}$ and $\textit{directionalFilter}$,
ShearLab 3D can compute a set of 2D digital shearlet filters whose inner products with a given 2D signal (and all its translates)
are the desired shearlet coefficients (see Algorithm \ref{alg:forward_transform}). As each shearlet coefficient corresponds to one
shearlet with a specific scale, a specific shearing, and a specific translation, the total number of coefficients computed by one
shearlet decomposition is $X\cdot Y\cdot R$, where $X$ and $Y$ denote the size of the given signal (and therefore the number of
different translates) and $R$ denotes the redundancy of the shearlet system, which is defined by the parameters $\textit{nScales}$
and $\textit{shearLevels}$. In fact, the redundancy $R$ including the low frequency part is given by
\[
R = 1+\sum_{j = j_0}^{j_0+\textit{nScales}-1} 2(2\cdot 2^{d_j}+1),
\]
where $j_0$ is the coarsest scale $j = j_0$ for the shearlet transform, and one can specify any nonnegative integer for $j_0$.

Let us now consider the 3D situation.
Due to the formula for ${\psi}^{\text{d}}_{j,k}$ in Definition \ref{defi:3DdigitSH}, we know that a 3D digital shearlet filter can be
constructed by combining two 2D digital shearlet filters living on the same scale but with possibly differing shearings. Therefore,
the input parameters for a three-dimensional decomposition are the same as in the 2D case but their meanings slightly differ. Where
in 2D, each scale corresponds to a ring-like passband, in the 3D case each scale is associated with a sphere-like passband in the
3D frequency domain with the parameter $nScales$ defining the number of such spheres. The parameter $shearLevels$ on the other hand
still defines the number of differently sheared atoms on one scale but there are two shearing parameters $k = (k_1,k_2)$ for a
shearlet in 3D and the 3D frequency domain is partitioned in three pyramids instead of two cones.

To be more precise, let $d_j$
be a $j$-th component in $\textit{shearLevels}$. Then the shearlet filters $\psi^{\text{d}}_{j,k}$ defined in Definition
\ref{defi:3DdigitSH} and associated with the 3D scaling matrix $A_{d_j}$ are generated. In
this case, the range of shearing parameters $k = (k_1,k_2)$ is given by $\max\{|k_1|,|k_2|\} \leq 2^{d_j}$ for each pyramid, which
gives  $3(2\cdot2^{d_j}+1)^2$ shearlet filters for each scale $j$. Thus, in the 3D case, the redundancy $R$ is given by
\[
R = 1+\sum_{j = j_0}^{j_0+\textit{nScales}-1} 3(2\cdot 2^{d_j}+1)^2.
\]

\subsection{Inverse Transform} \label{subsec:inversetrafo}

A schematic description of the inverse transform is given in Algorithm~\ref{alg:inverse_transform}.

 \begin{algorithm}[h]
 \SetAlgoLined
 \KwIn{A set of shearlet coefficients $\textit{shearletCoeffs} \in \R^{X\times Y\times Z \times R}$ and dual shearlet filters $\textit{dualFilters}
 \in \C^{X\times Y\times Z\times R}$ in the frequency domain, where $R$ denotes the redundancy of the corresponding shearlet system.}
 \KwOut{A reconstructed signal $\textit{f}_{\textit{rec}}\in \R^{X\times Y\times Z}$. }

  $\textit{f}_{\textit{rec}}$ := 0

    \tcp{The reconstructed signal $\textit{f}_{\textit{rec}}$ is computed as the sum over all convolutions of the dual shearlet filters
    with their corresponding three-dimensional vectors (two-dimensional for 2D data) of shearlet coefficients. Due to the convolution theorem,
    this can be computed using pointwise multiplication (denoted by .*) in the frequency domain.}
  \For{\textit{i} := 1 \KwTo R}{
  $\textit{f}_{\textit{rec}}$ := $\textit{f}_{\textit{rec}}$ +  FFT(\textit{shearletCoeffs(i)}).*\textit{dualFilters(i)}\;
  }
  $\textit{f}_{\textit{rec}}$ := $IFFT\left(\textit{f}_{\textit{rec}}\right)$\;
 \caption{ShearLab 3D inverse transform}
 \label{alg:inverse_transform}
\end{algorithm}

To perform an inverse transform in both 2D and 3D (see algorithm \ref{alg:inverse_transform}), ShearLab 3D requires a set of shearlet
coefficients and the corresponding digital shearlet filters from which the dual filters can be computed according to formula \eqref{eq:dual_filter}.
The reconstructed signal is then the sum over all dual filters multiplied with their corresponding coefficients.

\subsection{Example} \label{subsec:example}

\begin{figure}[h]
\begin{center}
\includegraphics[width=1.4in]{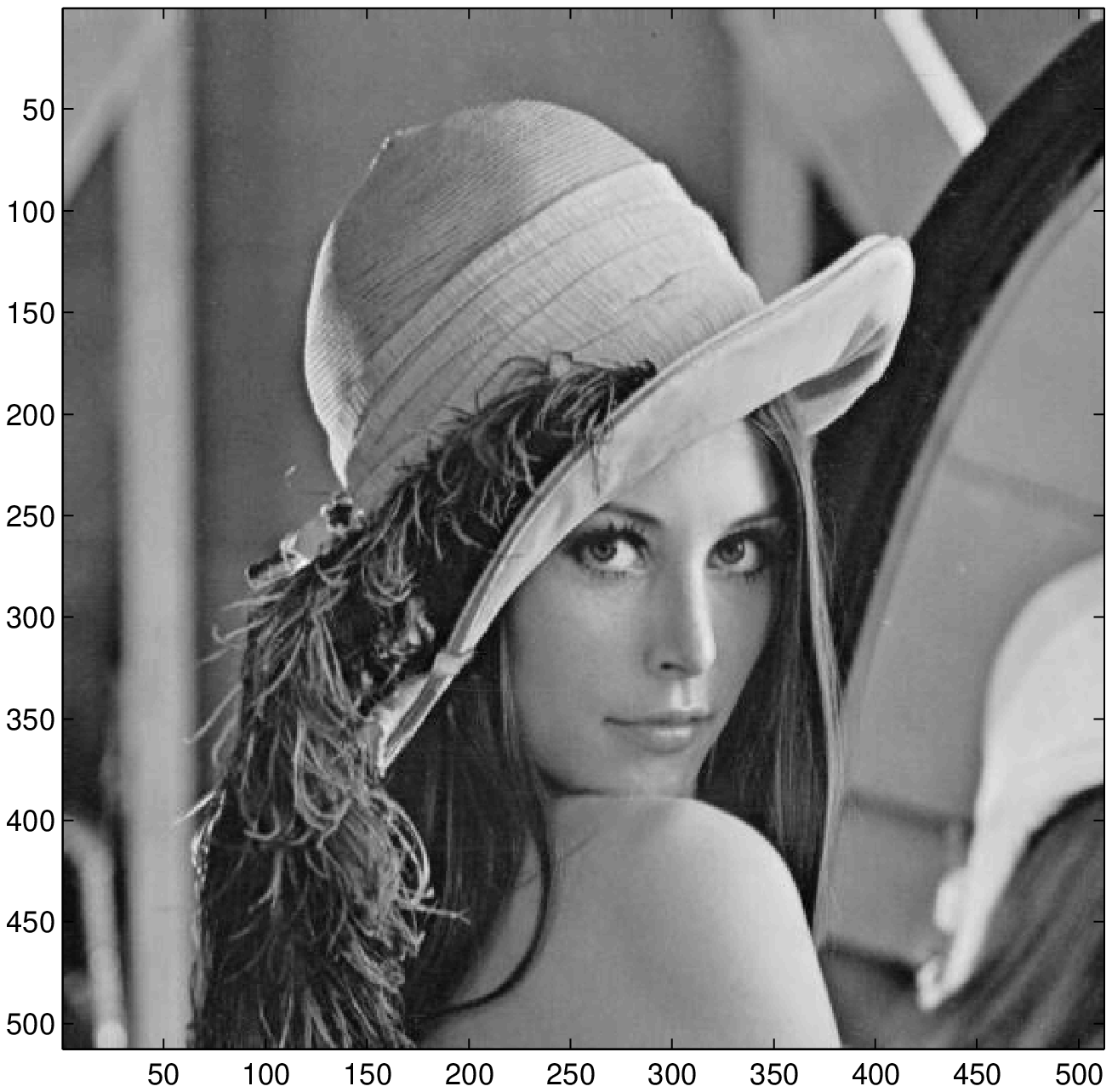}
\hspace{20pt}
\includegraphics[width=1.4in]{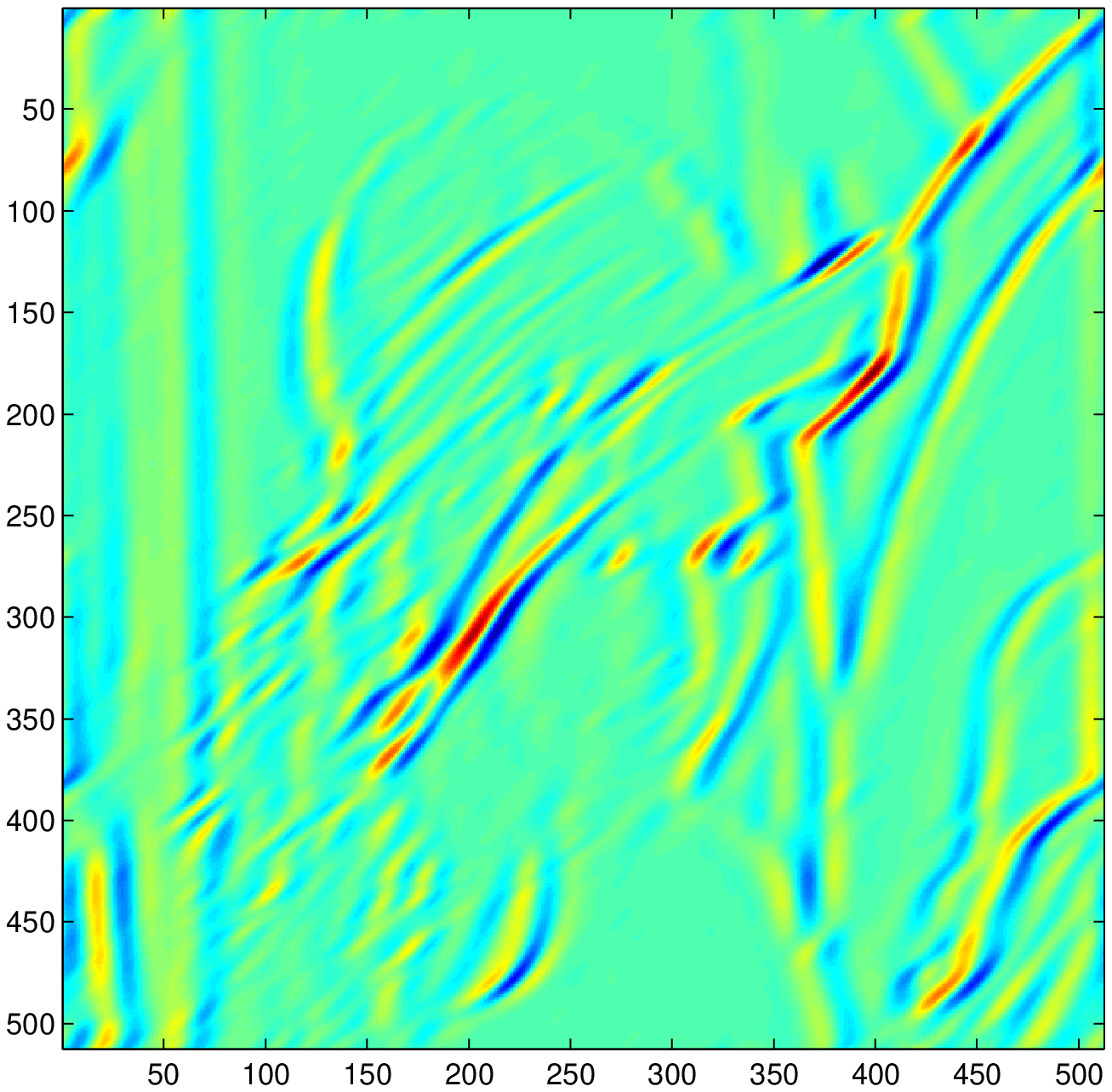}
\hspace{20pt}
\includegraphics[width=1.4in]{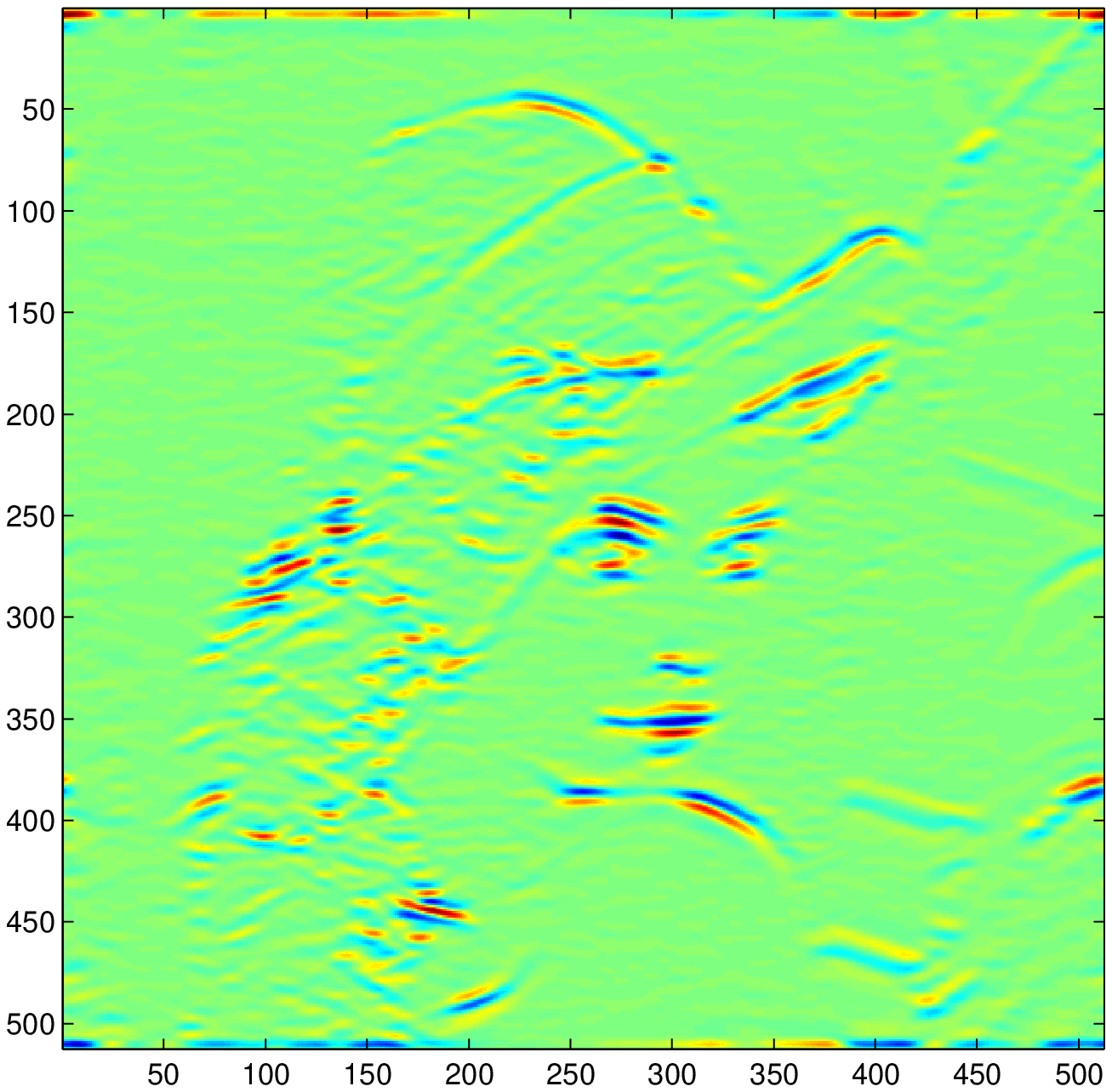}
\caption{The two images to the right show all shearlet coefficients of the translates of two different shearlets. The used system
has four scales, a redundancy of $49$ and was specified with $\textit{nScales} = 4$ and $\textit{shearLevels} = [1,1,2,2]$. The shearlet corresponding to the coefficients
in the centered picture has a scale parameter $j = 1$ a shearing parameter $k = 2$ and lives on the horizontal frequency cones. The
shearlet corresponding to the coefficients plotted in the rightmost image has a scale parameter $j = 2$, a shearing parameter $k = 0$
and lives on the vertical frequency cones.}
\label{fig:sheardec}
\end{center}
\end{figure}

We next provide a brief example of how to use ShearLab 3D to decompose and reconstruct a 2D signal in MATLAB.
Figure \ref{fig:shearlab_example} shows the MATLAB coding side, and Figure \ref{fig:sheardec} the visual outcome.
\begin{figure}[h]
\begin{center}
\begin{MatlabCode}
nScales = 4;
shearLevels = [1,1,2,2];
% this flag is used to determine whether CUDA is used or not.
useGPU = 0;
% this flag can be used to inhibit the omission of sheralet filters on the seam lines.
% it is optional with default value 0.
fullSystem = 0;

% load the lenna image.
img = double(imread('lenna.jpg'));

% construct digital shearlet filters with a default 1D lowpass and
% a default 2D directional filter.
system = SLgetShearletSystem2D(useGPU,size(img,1),size(img,2),nScales,shearLevels,fullSystem);

% compute shearlet decomposition
shearletCoefficients = SLsheardec2D(img,system);

% compute shearlet reconstruction
reconstruction = SLshearrec2D(shearletCoefficients,system);
\end{MatlabCode}
\caption{Decomposition and reconstruction of a 2D image in MATLAB using ShearLab 3D. The applied shearlet system has four scales
and the array \textit{shearLevels} induces a parabolic scaling. The array \textit{shearletCoefficients} is three-dimensional in
the 2D and four-dimensional in the 3D case. In both cases, the last dimension enumerates all digital shearlet filters within the
specified system with different shearing parameters $k$ and scaling parameters $j$ while the first two or three dimensions are
associated with the translates of one single shearlet (see Figure \ref{fig:sheardec}).}
\label{fig:shearlab_example}
\end{center}
\end{figure}

\subsection{Omission of Boundary Shearlets} \label{subsec:boundary}

\begin{figure}
 \begin{center}
 \includegraphics[width=1.5in]{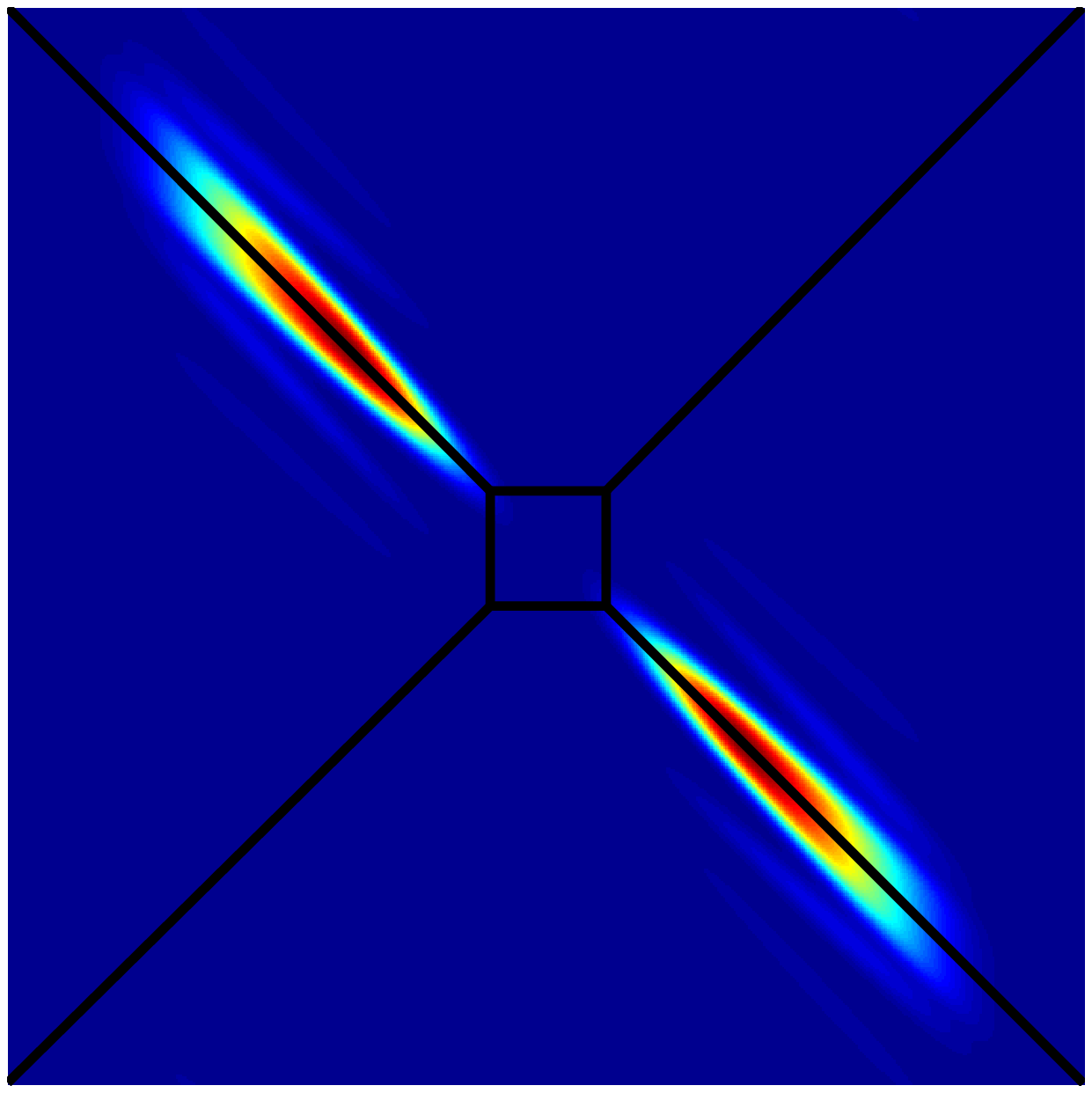}
 \hspace{50pt}
\includegraphics[width=1.5in]{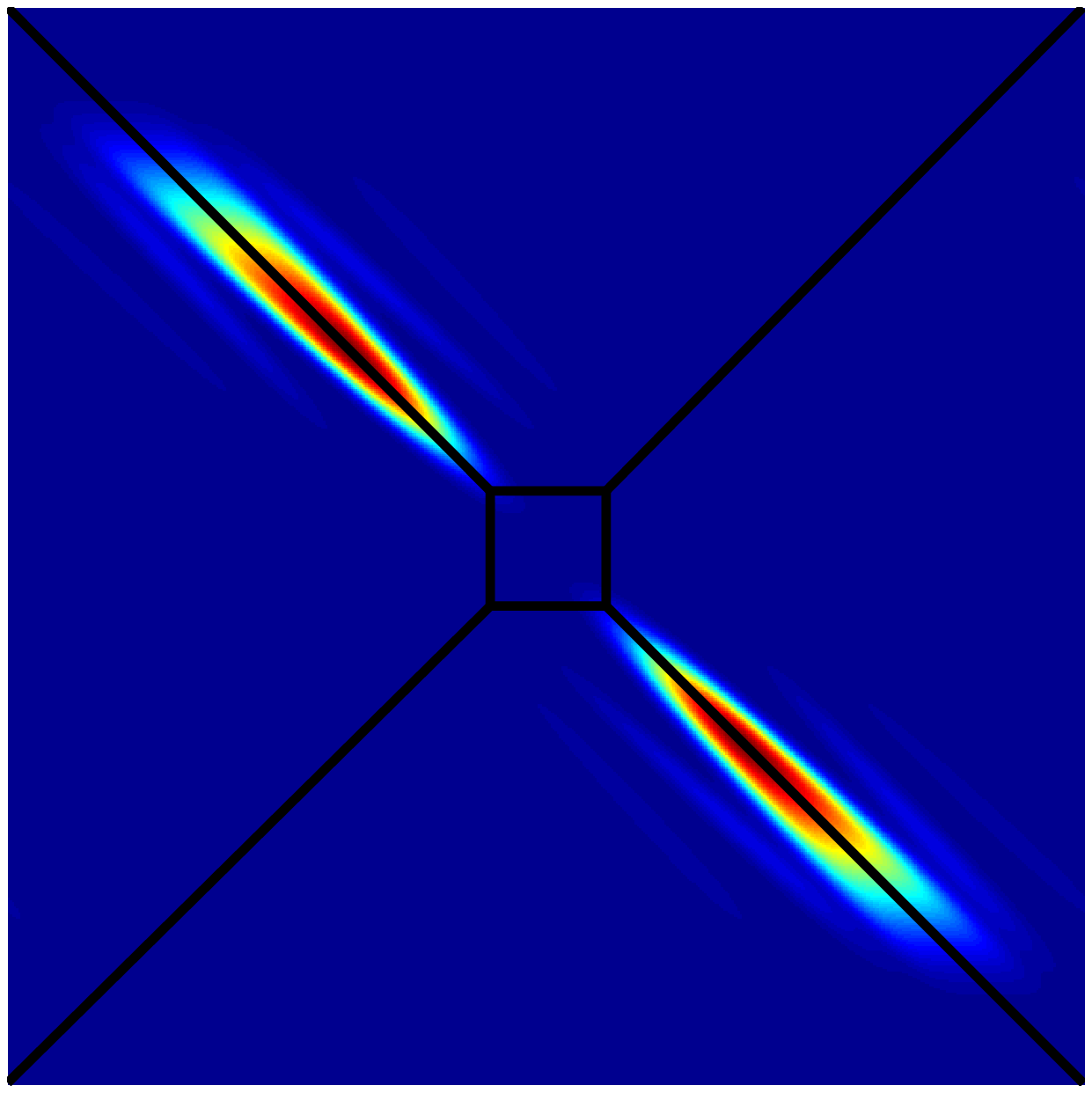}
\caption{The magnitude frequency response of the maximally sheared shearlet in the vertical cones (right image) is almost equal
to the response of the corresponding shearlet in the horizontal cones (left image). In most cases, one of these filters can be
omitted to decrease the redundancy of a shearlet system.}
\label{fig:omission}
 \end{center}
 \end{figure}

Let now $\psi^{\text{d}}_{j,k}$ and $\tilde{\psi}^{\text{d}}_{j,k}$, $k = \pm d_j$, be the four shearlet filters as
defined in Definition \ref{defi:digitSH_alpha}. We notice that the support of each of those filters is concentrated
on the boundary of the horizontal cones (or the vertical cones). Even more, the filter $\psi^{\text{d}}_{j,k}$
is almost identical to $\tilde{\psi}^{\text{d}}_{j,k}$ when $k = \pm d_j$, which is illustrated in Figure
\ref{fig:omission}.
For this reason, in ShearLab~3D the boundary shearlet filters in the vertical cones are removed for each scale $j$ in
order to improve stability as well as efficiency. This leads to $2(2\cdot 2^{d_j}+1)-2 = 2^{d_j+2}$ shearlet filters
at each scale $j$. As a example, letting  $d_j = \lceil j/2 \rceil$ yields $8,8,16,16$ shearlet filters for $j = 1,2,3,4$
(see Figure \ref{fig:systems_properties}), which corresponds to the parabolic case. In this example, the redundancy
$R$ can be computed to be $1+8+8+16+16 = 49$.

A similar strategy can be applied in the 3D case: All boundary shearlet filters whose frequency support is concentrated
on the boundary of two pyramids among three are removed, yielding
\[
3(2\cdot 2^{d_j}+1)^2-6(2\cdot2^{d_j}+1)+4
\]
3D shearlet filters for each scale $j$. Again as an example, consider $d_j = \lceil j/2\rceil$. In this case, we have
$49,49,193$ shearlet filters for $j = 1,2,3$ (compare also Figure \ref{fig:systems_properties}), which corresponds to
the parabolic case.

It should be mentioned that in ShearLab 3D the user has the option of including those boundary elements again if
needed. 

\subsection{Complexity} \label{subsec:complexity}

Since convolution can be used, both the decomposition and the reconstruction algorithm reduce to multiple computations of the Fast Fourier Transform.
Thus, their complexity is given by $\mathcal{O}\left(R\cdot N\log (N)\right)$, where $R\in\mathbb{N}$ is the redundancy of the specific digital shearlet
system, i.e., the number of digital shearlet filters filters $\psi^d_{j,k}$.

\section{Numerical Experiments} \label{sec:experiments}

This section is devoted to an extensive set of numerical experiments. The parameters in ShearLab 3D chosen for those results are specified in
Subsection \ref{subsec:parameters}, followed by a detailed description of the transforms we compare our results to (see Subsection \ref{subsec:comparison}).
We then focus on the following problems: 2D/3D denoising, 2D/3D inpainting, and also 2D decomposition of point and curvelike structures,
which are the contents of Subsections \ref{subsec:denoising} to \ref{subsec:v_inpainting}. The results of the numerical experiments are
discussed in Subsection \ref{subsec:discussion}.

All experiments have been performed with MATLAB 2013a on an Intel(R) Core(TM)2 Duo CPU E6750 processor with 2.66GHz and a NVIDIAGeForce GTX 650 Ti
graphics card with 2GB RAM. The scripts and input data for all experiments are available at \url{www.shearlab.org} in support of the idea of
reproducible research.

\subsection{Selection of Parameters} \label{subsec:parameters}

As discussed in Subsection \ref{subsec:2Dtrafo}, the construction of a 2D digital shearlet filter $\psi_{j,k}^d$, requieres a 1D lowpass
filter $h_1$ and a 2D directional filter $P$, compare equation \eqref{eq:2dshear_filters}. The 1D filter $h_1$ defines a wavelet
multiresolution analysis -- and thereby the highpass filter $g_1$ --, whereas the trigonometric polynomial $P$ is used to ensure a wedge
shape of the essential frequency support of $\psi_{j,k}^d$. The choice of these filters certainly significantly impacts crucial
properties of the generated digital shearlet system such as frame bounds and directional selectivity.

\begin{figure}[h]
\begin{center}
\includegraphics[width=1.4in]{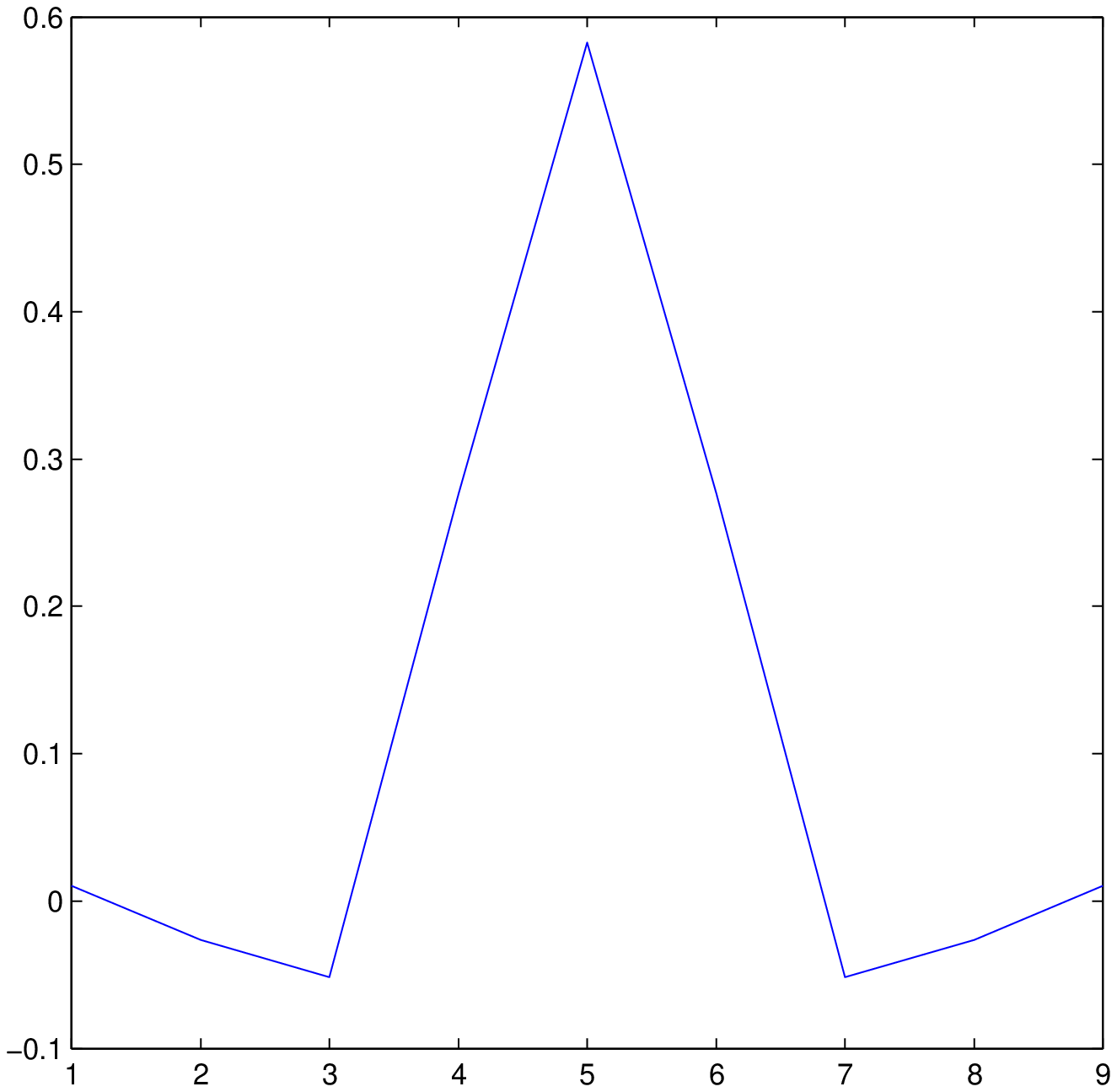}
\hspace{20pt}
\includegraphics[width=1.4in]{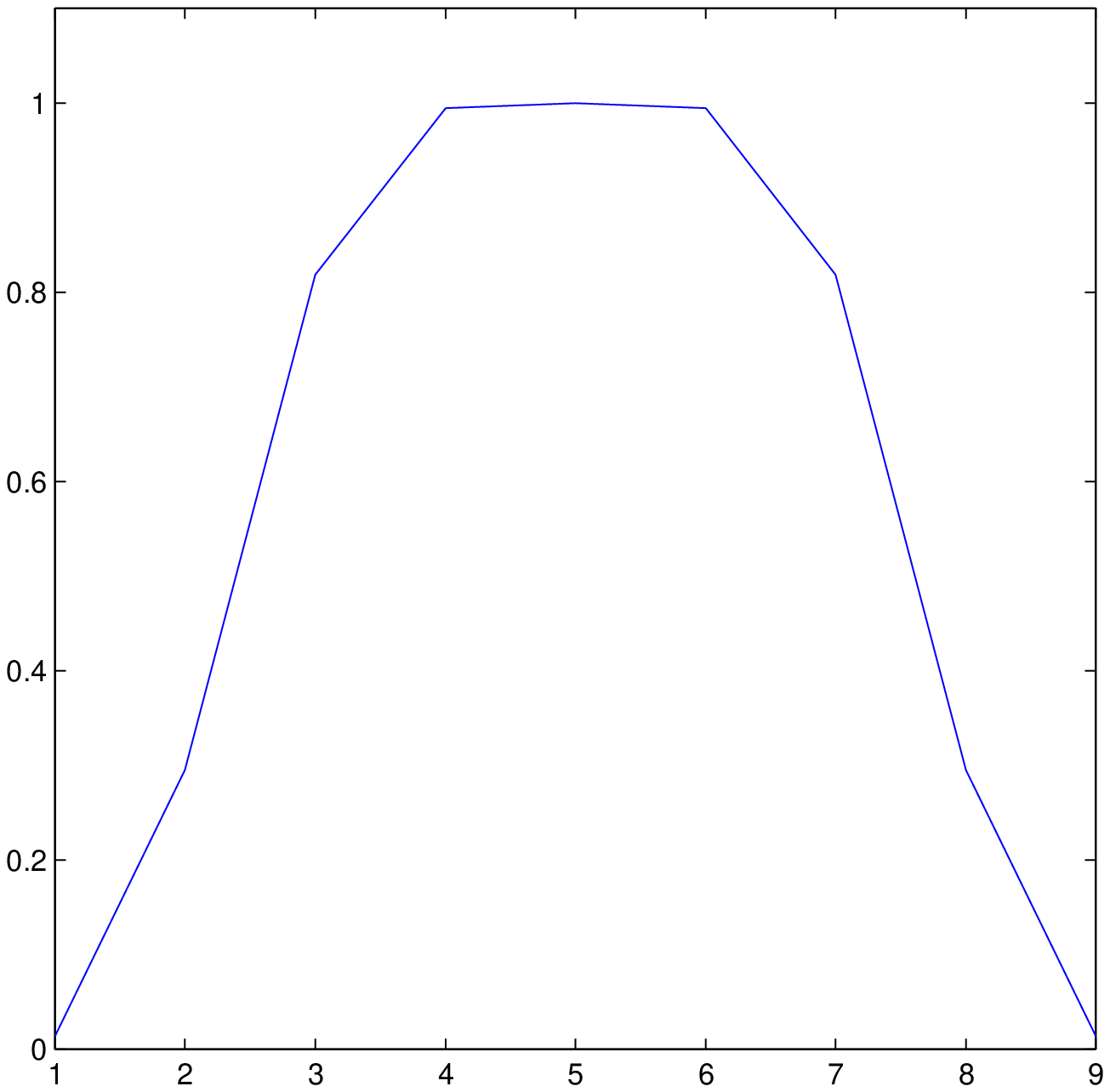}
\hspace{20pt}
\includegraphics[width=1.4in]{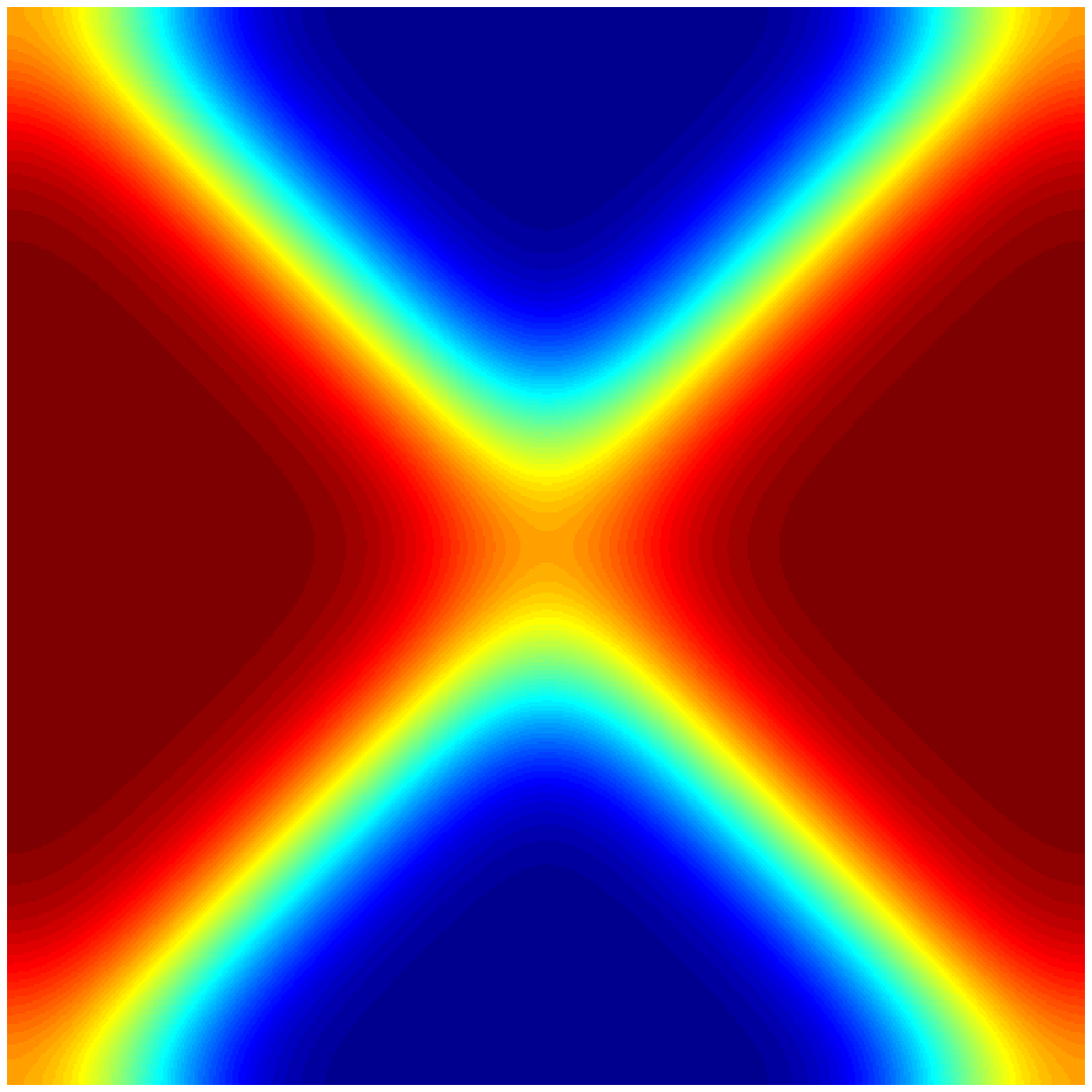}
\put(-308,-10){(a)}
\put(-182,-10){(b)}
\put(-55,-10){(c)}
\caption{(a) The coefficients of the 1D lowpass filter $h_{ShearLab}$. (b) Magnitude frequency response of $h_{ShearLab}$. (c) Magnitude response of the 2D fan filter $P_{ShearLab}$.}
\label{fig:filters}
\end{center}
\end{figure}

Our choice for $h_1$ -- from now on denoted by $h_{ShearLab}$ -- is a maximally flat, i.e. a maximum number of derivatives of the magnitude
frequency response at $0$ and $\pi$ vanish, and symmetric 9-tap lowpass filter\footnote{The MATLAB command \texttt{design(fdesign.lowpass('N,F3dB',8,0.5),'maxflat')}
generates the 9-tap filter $h_{ShearLab}$, whose approximate values are $h_{ShearLab} = (0.01049, -0.02635, -0.05178, 0.27635, 0.58257, ...)$.}
which is normalized such that $\sum_n h_{ShearLab}(n) = 1$. For an illustration, we refer to Figure \ref{fig:filters}(a) and (b). This filter has two vanishing
moments, i.e. $\int\phi(x)x^k = 0$ for $k\in\{0,1\}$. While there is no symmetric, compactly supported, and orthogonal wavelet besides the Haar wavelet,
the renormalized filter $\sqrt{2}h_{ShearLab}$ at least approximately fulfills the orthonormality condition, which is
\[
\abs{2\sum_{n}h_{ShearLab}(n)h_{ShearLab}(n+2l) - \delta_{l0}}\leq 0.0018 \quad \mbox{for all } \l\in\Z,
\]
with $\delta$ denoting Kronecker's delta. We remark that by choosing $h_{ShearLab}$ to be maximally flat, the amount of ripples in the digital filter $\psi_{j,k}^d$
is significantly reduced. This leads to an improved localization of the associated digital shearlets in the frequency domain.
The highpass filter $g_1$, hereafter denoted by $g_{ShearLab}$, is certainly chosen to be the associated mirror filter, that is
\[
 g_{ShearLab}(n) = (-1)^{n}\cdot h_{ShearLab}(n).
\]

We would like to mention that the filter coefficients $h_{ShearLab}$ are quite similar to those of the Cohen-Daubechies-Feauveau (CDF)
9/7 wavelet \cite{CDF1992}, which is used in the JPEG 2000 standard. While the CDF 9/7 wavelet has four vanishing moments and higher
degrees of regularity both in the H\"older and Sobolev sense, trading these advantageous properties for maximal flatnass seems to be
the optimal choice for most applications.

For the trigonometric polynomial $P_{ShearLab}$, we use the maximally flat 2D fan filter\footnote{The 2D fan filter $P_{ShearLab}$ can be obtained
in MATLAB using the Nonsubsampled Contourlet Toolbox by the statement \texttt{fftshift(fft2(modulate2(dfilters('dmaxflat4','d')./sqrt(2),'c')))}.}
described in \cite{dCZD2006}. This filter is illustrated in Figure \ref{fig:filters}(c).

For the numerical experiments, we used two different digital shearlet systems in both the 2D and 3D case, the reason being that
this allows us to demonstrate how different degrees of redundancy influence the performance of ShearLab 3D. The two considered
2D systems named $SL2D_1$ and $SL2D_2$, similarly $SL3D_1$ and $SL3D_2$ in 3D, are constructed as follows. Notice that they all
correspond to the case $\alpha = 1$, hence the parabolic case.

The system $SL2D_1$ has four scales with four differently oriented digital shearlet filters on scales one and two, and eight
directions in each of the higher scales. Including the 2D lowpass filter, the total redundancy of this system is $25$. We remark
that a maximally sheared filter within the horizontal cones has always an almost identical counterpart contained in the vertical
cones that can be omitted without affecting the performance in most applications. The system $SL2D_2$ also consists of four scales,
but in contrast to $SL2D_2$ has a redundancy of $49$ with $8, 8, 16$ and $16$ differently oriented shearlets on the respective
scales.

In the 3D experiments, we used the three-scale digital systems $SL3D_1$ with $13, 13$ and $49$ directions on scales one, two and
three as well as the system $SL3D_2$ with $49, 49$ and $193$ differently oriented shearlet filters on the corresponding scales.
The total redundancy of $SL3D_1$ is $76$, while $SL3D_2$ contains $292$ different digital 3D shearlet filters. These numbers along
with other properties of these systems are compiled in Figure \ref{fig:systems_properties}.

\begin{figure}[h]
\begin{center}
\begin{tabular}{|c|c|c|c|c|c|c|}
\hline
&Scales&Directions&Redundancy&A&B&B/A\\
\hline
$SL2D_1$&4&(4,4,8,8)&25&0.0893&1.0000&11.19\\
\hline
$SL2D_2$&4&(8,8,16,16)&49&0.0669&1.0000&14.94\\
\hline
$SL3D_1$&3&(13,13,49)&76&0.0075&1.0000&133.39\\
\hline
$SL3D_2$&3&(49,49,193)&292&0.0045&1.0000&220.84\\
\hline
\end{tabular}
\caption{Properties of the digital shearlet systems $SL2D_1$, $SL2D_2$, $SL3D_1$, and $SL3D_2$ used in our numerical
experiments. Columns A and B show approximates of the lower and upper frame bounds.}
\label{fig:systems_properties}
\end{center}
\end{figure}

\subsection{Systems for Comparison} \label{subsec:comparison}

In a total of five different experiments, we compare the transforms associated with the shearlet systems $SL2D_1$, $SL2D_2$,
$SL3D_1$, and $SL3D_2$ to various transforms similarly associated with a specific representation system. In each of these experiments, an algorithm
based on a sparse representation of the input data is used to complete a certain task like image denoising or image inpainting.
In order to get a meaningful comparison, we simply run the same algorithm with the same input several times for each of the
transforms -- i.e., with the associated representation system -- for computing the sparse representation at each execution. To
assess the performance of a sparse representation scheme, for each of the tasks we introduce a performance measure for the quality
of the output and also measure the overall running time of the algorithm.

The transforms considered in our experiments besides the digital shearlet transform implemented in ShearLab 3D are:
\begin{itemize}
 \item[$\bullet$] {\em Nonsubsampled Shearlet Transform ($NSST$, 2D \& 3D)}.\\[0.5ex]
 The NSST was introduced in 2006 by Labate, Easley, and Lim \cite{ELL08} and later extended to a 3D transform by Negi and Labate \cite{NL2013}.
 It is also based on the theory of shearlets and uses the nonsubsampled Laplacian pyramid transform with specially designed bandpass filters \cite{dCZD2006} to decompose input data into several high-frequency layers and a low-frequency part, while directional filters are constructed on the pseudo-polar
 grid from a certain window function, e.g. the Meyer wavelet window. The main conceptional difference to ShearLab 3D is that these directional filters are not
 compactly supported in the time domain. Still, by representing the directional filters with small matrices, the NSST also manages to construct digital shearlet
 filters that are highly localized in the time domain. An implementation is publicly available at \url{http://www.math.uh.edu/~dlabate/software.html}.\\[-1.2ex]
 \item[$\bullet$] {\em Nonsubsampled Contourlet Transform ($NSCT$, 2D)}.\\[0.5ex]
  The NSCT was developed by da Cunha, Zhou, and Do \cite{dCZD2006} and uses a nonsubsampled pyramid decomposition as well as a directional filter bank based on
  two-channel fan filter banks to construct directionally sensitive digital filters on several scales. It was shown that the NSCT can be applied to construct
  frames for $\ell_2(\Z^2)$ and that (band-limited) contourlets can achieve the optimal approximation rate for cartoon-like images \cite{DV2005}.
 The Nonsubsampled Contourlet Toolbox can be downloaded from \url{http://www.mathworks.com/matlabcentral/fileexchange/10049-nonsubsampled-contourlet-toolbox}.\\[-1.2ex]
 \item[$\bullet$] {\em Fast Discrete Curvelet Transform ($FDCT$, 2D)}.\\[0.5ex]
 Curvelets were first introduced in 1999 by Cand\'{e}s and Donoho \cite{CD1999} with the goal of constructing a non-adaptive frame of representing functions
 providing optimal approximation rates for cartoon-like images. Indeed, the curvelet transform was the first non-adaptive method published to achieve this and
 can be viewed as a precursor to the theory of shearlets. The most significant conceptional differences are that the shearlet transform is associated with a
 single (or finite set of) generating function(s) that can be subject to anistropic scaling and shearing, while curvelet atoms are constructed by rotating
 'mother' curvelets that exist on each scale. The FDCT used in our experiments is described in \cite{CDDY2006} and can be downloaded from
 \url{http://www.curvelet.org/software.html}.\\[-1.2ex]
 \item[$\bullet$] {\em Surfacelet Transform ($SURF$, 3D)}.\\[0.5ex]
 In the surfacelet transform, first published by Do and Lu in 2007 \cite{DL2007}, the two-dimensional Bamberger and Smith directional filter bank \cite{BS1992}
 (which is also used in the NSCT) is extended to higher dimensions. Together with a pyramid transform similar to one developed by Simoncelli and Freeman
 \cite{SFAH1992}, a directional selective three-dimensional multiscale transform can be constructed. An implementation is available at
 \url{http://www.mathworks.com/matlabcentral/fileexchange/14485-surfacelet-toolbox}.\\[-1.2ex]
 \item[$\bullet$] {\em Stationary Wavelet Transform ($SWT$, 2D)}.\\[0.5ex]
 The SWT, also known as algorithme \`{a} trous, is a redundant and translation invariant version of the discrete wavelet transform. Instead of dyadically
 downsampling the signal at each transition from one scale to another, the filter coefficients are dyadically upsampled. In our experiments, we used the
 method SWT2 from the MATLAB Wavelet Toolbox.
\end{itemize}

\subsection{Image Denoising} \label{subsec:denoising}

As input, we consider several grayscale images of size 512x512 that are distorted with Gaussian white noise. In order to denoise these images,
we use hard thresholding on the coefficients of a sparse representation scheme before computing the reconstruction. That is, for an image $f\in \ell^2(\Z^2)$ and
\[
f_{noisy}(i,j) = f(i,j) + e(i,j)
\]
where $e(i,j)  \sim \mathcal{N} (0,\sigma^2)$, we compute
\[
f_{denoised} = \mathcal{T}^{-1}T_\delta\mathcal{T}f_{noisy},
\]
where $\mathcal{T}$ is the forward and $\mathcal{T}^{-1}$ is the inverse transform associated with a certain sparse representation scheme
and $T_\delta$ is the hard thresholding operator given by
\begin{equation}
\label{eq:hard_tresholding}
(T_\delta x)(n) = \begin{cases}x(n) &\text{if }\abs{ x(n)} \geq \delta,\\0 &\text{else.}\end{cases}
\end{equation}
In order to increase the performance, we will use different thresholds $\delta_j$ on different scales $j$, that are of the form
\[
\delta_j = K_j\sigma
\]
where for four scales, we typically have $K = [K_j]_j = [2.5, 2.5, 2.5, 3.8]$.

The quality of the reconstruction is measured using the peak signal-to-noise ratio (PSNR), defined as
\begin{equation}
\label{eq:PSNR}
PSNR = 20\log_{10}\frac{255\sqrt{N}}{\|f - f_{denoised}\|_F}
\end{equation}
where $N$ is the number of pixels and $\|\cdot\|_F$ denotes the Frobenius norm. Note that $255$ is the maximum value a pixel can attain in a grayscale image.

We compare the performance of $SL2D_1$, $SL2D_2$, the nonsubsampled shearlet transform $NSST$, the nonsubsampled contourlet transform $NSCT$, the fast discrete
curvelet transform $FDCT$, and the stationary wavelet transform $SWT$. All quantitative results for a total of four grayscale images and several levels of noise
as well as the running times and redundancies associated with the considered transforms are displayed in Figure \ref{fig:results_denoising_2d}. For a visual
comparison, we refer to Figure \ref{fig:results_denoising_2d_images} in the Appendix.
\begin{figure}[h]
\begin{center}
\begin{tabular}{|c|c|c|c|}
\hline
\multicolumn{4}{|c|}{\textbf{Image Denoising: Running Times}}\\
\hline
&{Redundancy}&{Running Time}&{Running Time with CUDA}\\
\hline
$SL2D_1$&25&\unit{4.8}{\second}&\unit{1.7}{\second}\\
\hline
$SL2D_2$&49&\unit{10.9}{\second}&\unit{5.7}{\second}\\
\hline
$NSST$&49&\unit{9.4}{\second}&-\\
\hline
$NSCT$&49&\unit{452}{\second}&-\\
\hline
$FDCT$&2.85&\unit{0.8}{\second}&-\\
\hline
$SWT$&13&\unit{0.7}{\second}&-\\
\hline
\end{tabular}

\vspace{0.5cm}

\begin{tabular}{|l|c|c|c|c|c|c|c|c|c|c|}
\hline
\multicolumn{11}{|c|}{\textbf{Image Denoising: Quantitative Results in PSNR}}\\
\hline
&\multicolumn{5}{|c|}{{Lenna}}&\multicolumn{5}{|c|}{{Barbara}}\\
\hline
&$\sigma = ${$10$}&{$20$}&{$30$}&{$40$}&{$50$}&{$10$}&{$20$}&{$30$}&{$40$}&{$50$}\\\hline
\textbf{$SL2D_1$}&35.79&32.51&30.52&29.01&27.79&33.38&29.42&27.03&25.40&24.37\\\hline
\textbf{$SL2D_2$}&\textbf{35.90}&32.75&30.87&29.43&28.26&\textbf{33.63}&\textbf{29.98}&\textbf{27.83}&\textbf{26.28}&\textbf{25.17}\\\hline
\textbf{$NSST$}&35.85&\textbf{32.83}&\textbf{31.06}&\textbf{29.69}&\textbf{28.58}&33.56&29.91&27.75&26.20&25.15\\\hline
\textbf{$NSCT$}&35.67&32.53&30.69&29.30&28.14&33.43&29.58&27.36&25.89&24.88\\\hline
\textbf{$FDCT$}&34.01&31.41&29.71&28.38&27.33&28.96&25.47&24.48&23.87&23.43\\\hline
\textbf{$SWT$}&34.19&30.79&28.91&27.64&26.62&31.08&26.55&24.55&23.55&23.00\\\hline
\end{tabular}\\
\vspace{0.1cm}
\begin{tabular}{|l|c|c|c|c|c|c|c|c|c|c|}
\hline
&\multicolumn{5}{|c|}{{Boat}}&\multicolumn{5}{|c|}{{Peppers}}\\
\hline
&$\sigma = $\textbf{$10$}&\textbf{$20$}&\textbf{$30$}&\textbf{$40$}&\textbf{$50$}&\textbf{$10$}&\textbf{$20$}&\textbf{$30$}&\textbf{$40$}&\textbf{$50$}\\\hline
\textbf{$SL2D_1$}&33.06&30.00&28.16&26.87&25.86&34.10&31.78&30.04&28.67&27.49\\\hline
\textbf{$SL2D_2$}&\textbf{33.14}&\textbf{30.18}&\textbf{28.42}&\textbf{27.17}&26.18&\textbf{34.12}&\textbf{31.92}&\textbf{30.32}&\textbf{29.06}&27.97\\\hline
\textbf{$NSST$}&33.05&30.09&28.34&27.14&\textbf{26.20}&34.05&31.84&30.26&29.05&\textbf{28.02}\\\hline
\textbf{$NSCT$}&32.87&29.82&28.05&26.83&25.90&33.87&31.51&29.83&28.53&27.51\\\hline
\textbf{$FDCT$}&30.73&28.36&27.00&25.99&25.16&32.37&29.67&28.17&27.11&26.25\\\hline
\textbf{$SWT$}&31.54&28.23&26.45&25.31&24.47&33.09&30.39&28.44&26.99&25.92\\\hline
\end{tabular}
\caption{Numerical results of the image denoising experiment.}
\label{fig:results_denoising_2d}
\end{center}
\end{figure}

\subsection{Image Inpainting} \label{subsec:inpainting}

We consider a grayscale image $f\in\ell^2(\Z^2)$ to be partially occluded by a binary mask $M\in\{0,1\}^{\Z\times\Z}$, i.e.,
\[
f_{masked}(i,j) = f(i,j)M(i,j).
\]
The algorithm for inpainting the missing parts is based on an iterative thresholding scheme published in 2005 by Starck et al.
\cite{SED2005} (see also \cite{FSED2010} and Algorithm \ref{alg:image_inpainting}).
 \begin{algorithm}[h]
 \SetAlgoLined
 \KwData{$f_{masked}$, $M$, $\delta_{init}$, $\delta_{min}$, $iterations$}
 \KwResult{$f_{inpainted}$ }

 $f_{inpainted} := 0$\;
 $\delta := \delta_{init}$\;
 $\lambda := (\delta_{min})^{1/(iterations-1)}$\;

  \For{i := 1 \KwTo $iterations$}{
  $f_{res} := M.*(f_{masked} - f_{inpainted})$\tcp*{$.*$ denotes the pointwise multiplication}
  $f_{inpainted} := \mathcal{T}^{-1}T_\delta\mathcal{T}(f_{res} + f_{inpainted})$\tcp*[r]{Forward transform, thresholding, synthesis}
  $\delta := \lambda\delta$\;
  }

 \caption{Image inpainting via iterative thresholding }
 \label{alg:image_inpainting}
\end{algorithm}
In each step, a forward transform is performed on the unoccluded parts of the image combined with everything already inpainted in the missing areas.
The resulting coefficients are then subject to hard thresholding before an inverse transform is carried out. By gradually decreasing the thresholding
constant, this algorithm approximates a sparse set of coefficients whose synthesis is very close to the original image on the unoccluded parts.

\begin{figure}[h]
\begin{center}
\includegraphics[width=1.4in]{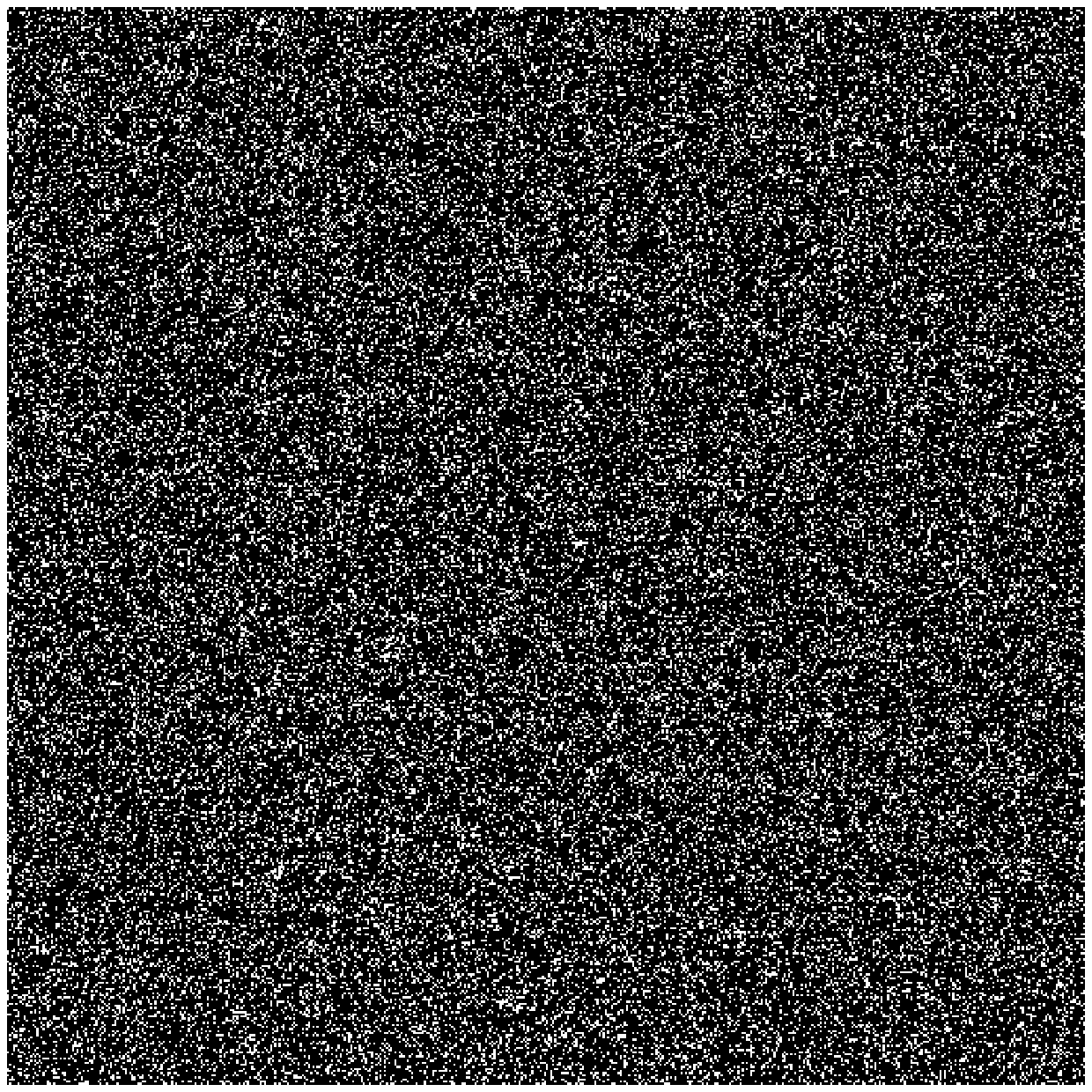}
\hspace{20pt}
\includegraphics[width=1.4in]{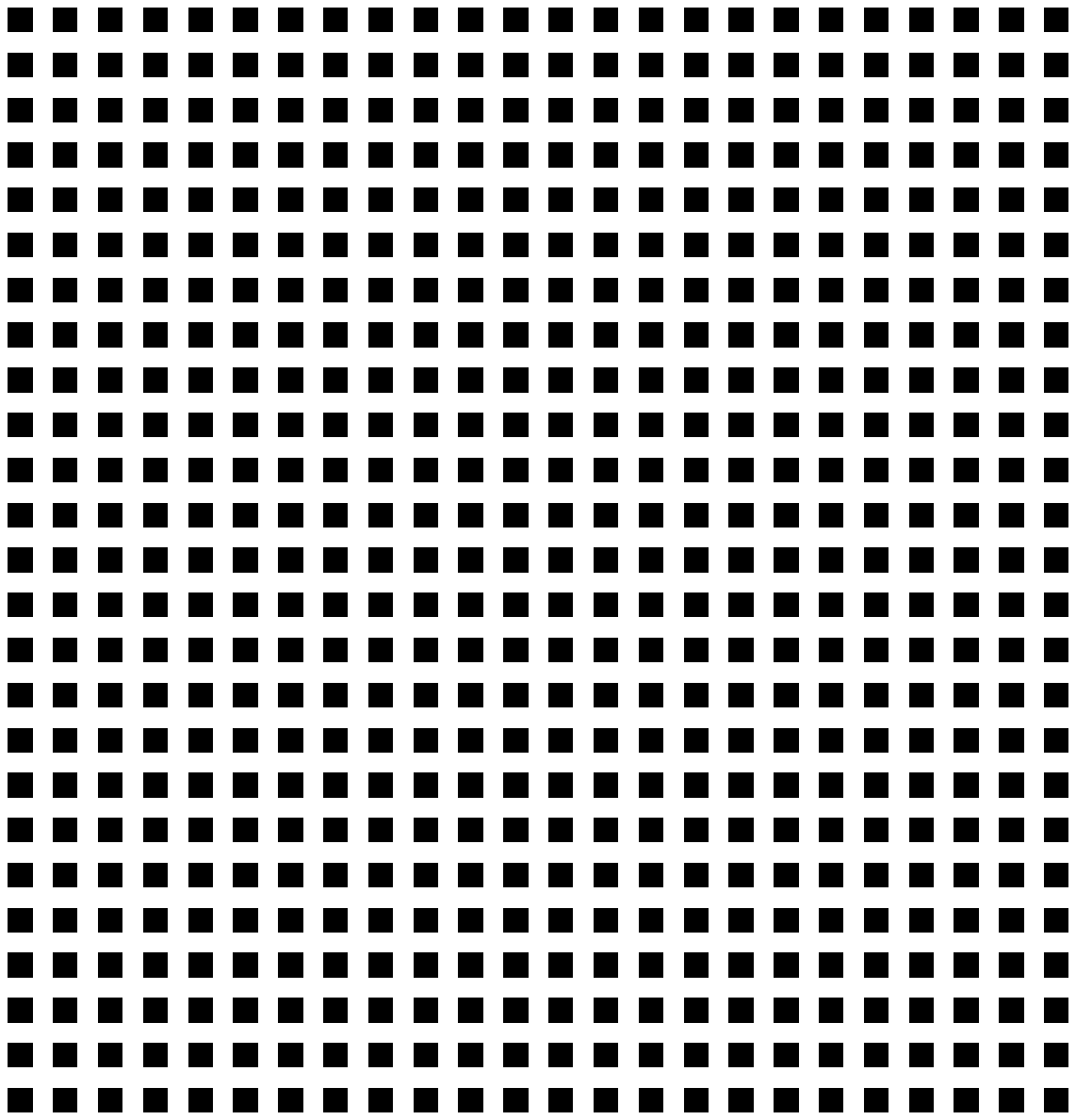}
\hspace{20pt}
\includegraphics[width=1.4in]{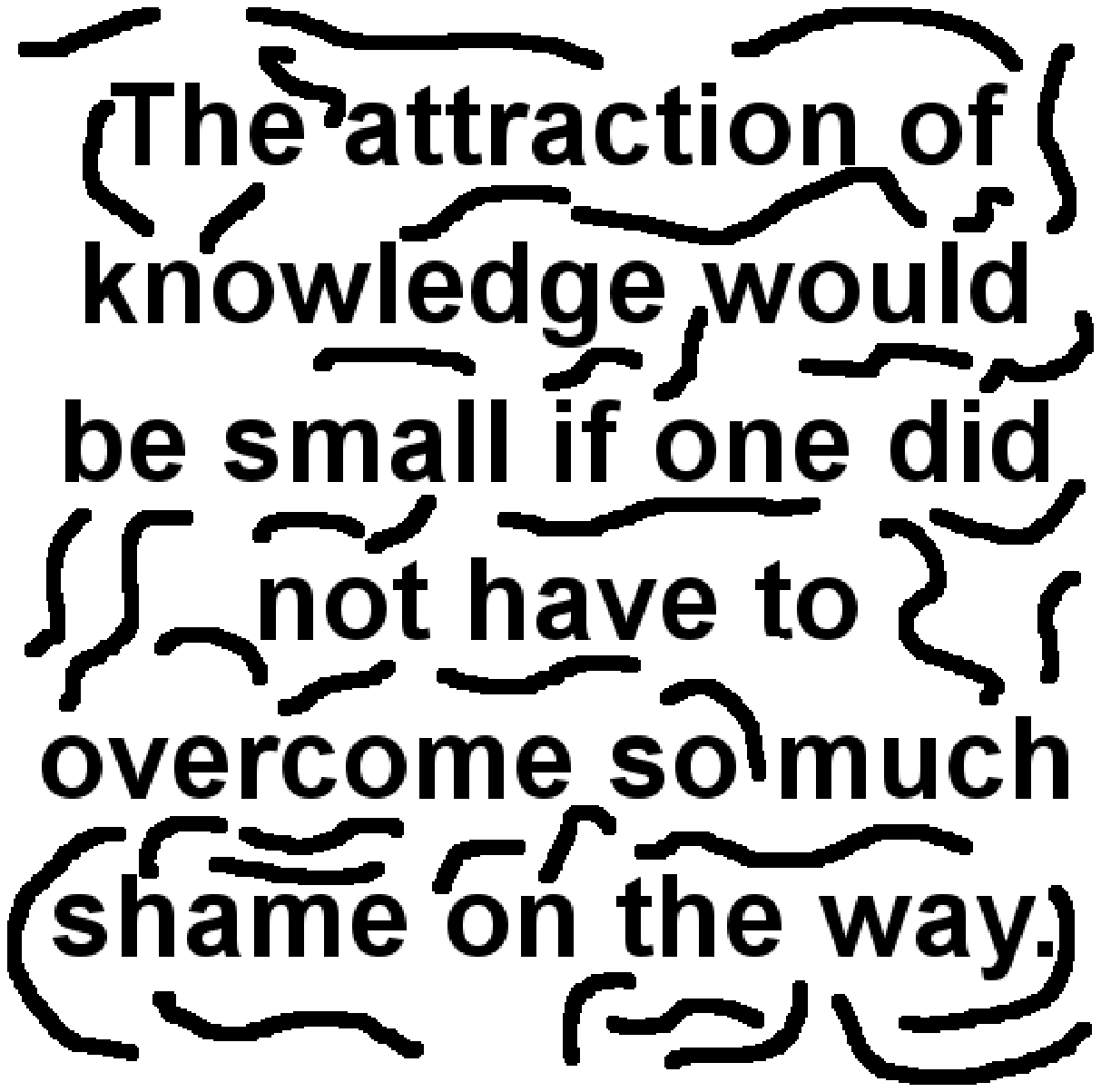}
\put(-360,-10){Random (\unit{80}{\%} occlusion)}
\put(-235,-10){Squares (\unit{28}{\%} occlusion)}
\put(-98,-10){Text (\unit{25}{\%} occlusion)}
\caption{Three different binary masks.}
\label{fig:inpainting_2d_masks}
\end{center}
\end{figure}

Again, we use the PSNR defined in \eqref{eq:PSNR} to measure the qualitative performance and to compare the systems $SL2D_1$ and $SL2D_2$ with the
nonsubsampled shearlet transform ($NSST$), the fast discrete curvelet transform ($FDCT$) and the stationary wavelet transform ($SWT$). To maximize the
quality of the output, we perform $300$ iterations during each inpainting task. Due to this heavy computational workload, we do not consider the
nonsubsampled contourlet transform ($NSCT$), as one iteration would already take more than $7$ minutes. Furthermore, we use three different types of
masks to properly emulate several typical inpainting problems, which are displayed in Figure~\ref{fig:inpainting_2d_masks}. For the quantitative results
and a compilation of running times and redundancies see Figure~\ref{fig:results_inpainting_2d}. A visual comparison is provided in Figure~\ref{fig:results_inpainting_2d_images} in the Appendix.

\begin{figure}[h]
\begin{center}
\begin{tabular}{|c|c|c|c|}
\hline
\multicolumn{4}{|c|}{\textbf{Image Inpainting: Running Times}}\\
\hline
&{Redundancy}&{Running Time}&{Running Time with CUDA}\\
\hline
$SL2D_1$&25&\unit{463}{\second}&\unit{45.3}{\second}\\
\hline
$SL2D_2$&49&\unit{931}{\second}&\unit{78.7}{\second}\\
\hline
$NSST$&49&\unit{2495}{\second}&-\\
\hline
$FDCT$&2.85&\unit{193}{\second}&-\\
\hline
$SWT$&13&\unit{162}{\second}&-\\
\hline
\end{tabular}

\vspace{0.5cm}

\begin{tabular}{|l|c|c|c|c|c|c|c|c|c|}
\hline
\multicolumn{10}{|c|}{\textbf{Image Inpainting: Quantitative Results in PSNR}}\\
\hline
&\multicolumn{3}{|c|}{{Lenna}} & \multicolumn{3}{|c|}{{Barbara}}&\multicolumn{3}{|c|}{{Flintstones}}
\\
\hline
&{Rand}&{Squares}&{Text}&{Rand}&{Squares}&{Text}&{Rand}&{Squares}&{Text}\\\hline
\textbf{$SL2D_1$}&32.17&33.11&31.44&26.59&30.08&29.70&23.84&24.16&22.98\\\hline
\textbf{$SL2D_2$}&32.08&\textbf{33.92}&32.57&\textbf{27.82}&\textbf{31.53}&\textbf{31.05}&{23.59}&\textbf{25.14}&\textbf{23.72}\\\hline
\textbf{$NSST$}&\textbf{32.31}&33.79&\textbf{32.78}&27.62&31.23&30.93&\textbf{24.12}&24.89&23.57\\\hline
\textbf{$FDCT$}&30.40&32.78&32.06&24.07&27.78&28.00&22.08&23.63&22.78\\\hline
\textbf{$SWT$}&30.46&31.69&30.95&23.84&28.46&28.75&21.86&22.91&22.18\\\hline
\end{tabular}
\caption{Numerical results of the image inpainting experiment.}
\label{fig:results_inpainting_2d}
\end{center}
\end{figure}

\subsection{Image Decomposition} \label{subsec:decomposition}

\begin{figure}[h]
\begin{center}
\includegraphics[width=1.4in]{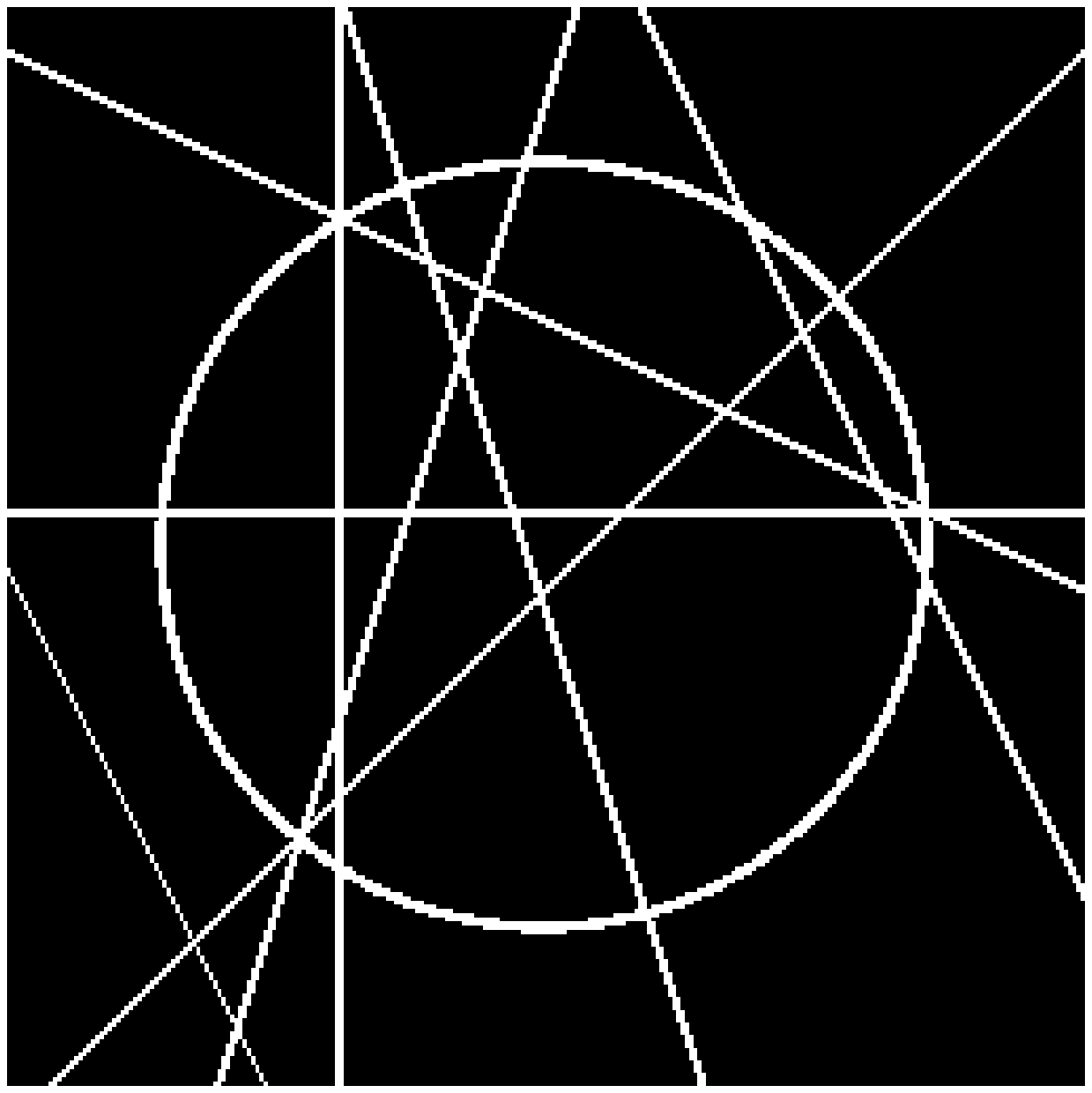}
\hspace{20pt}
\includegraphics[width=1.4in]{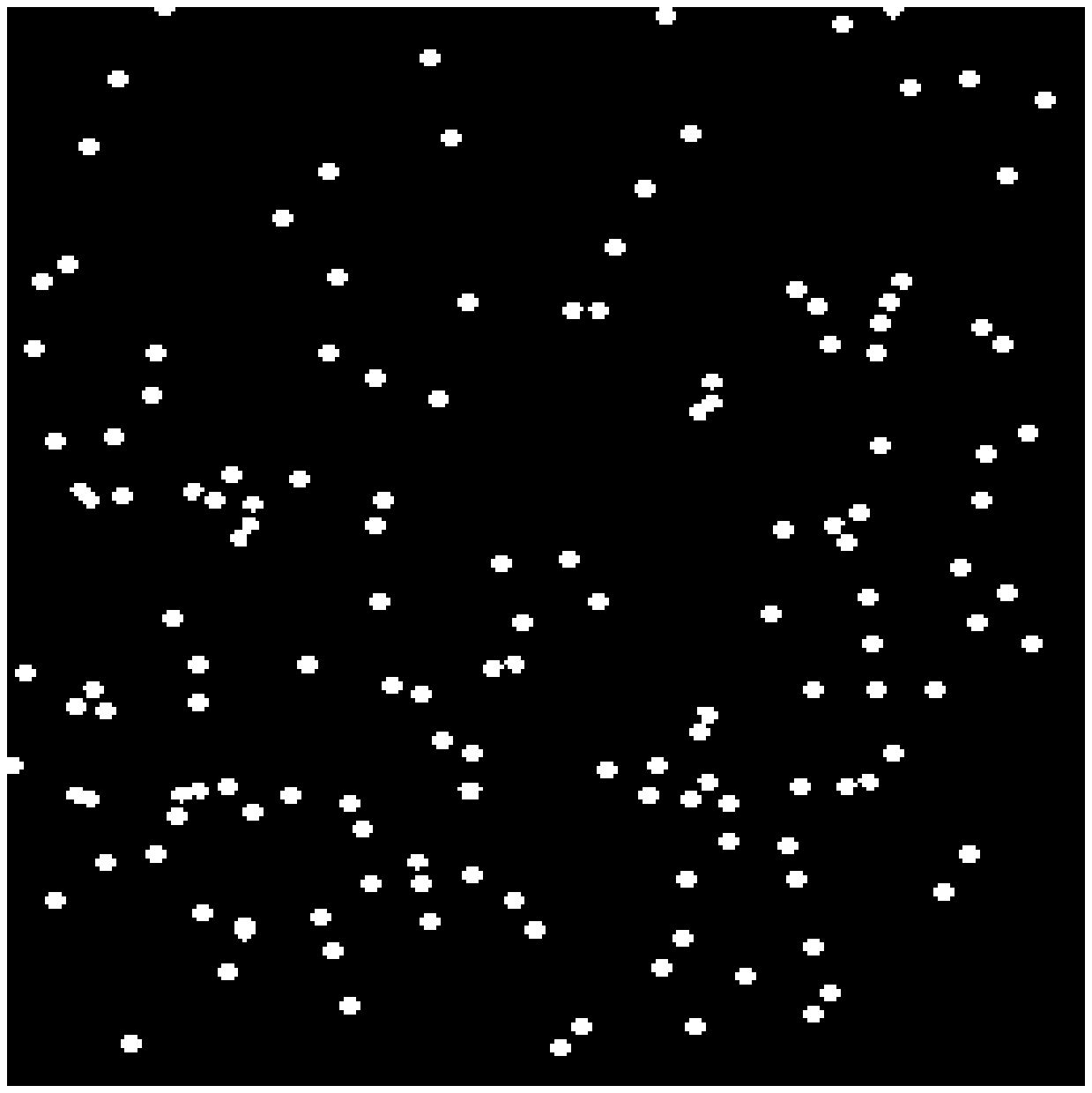}
\hspace{20pt}
\includegraphics[width=1.4in]{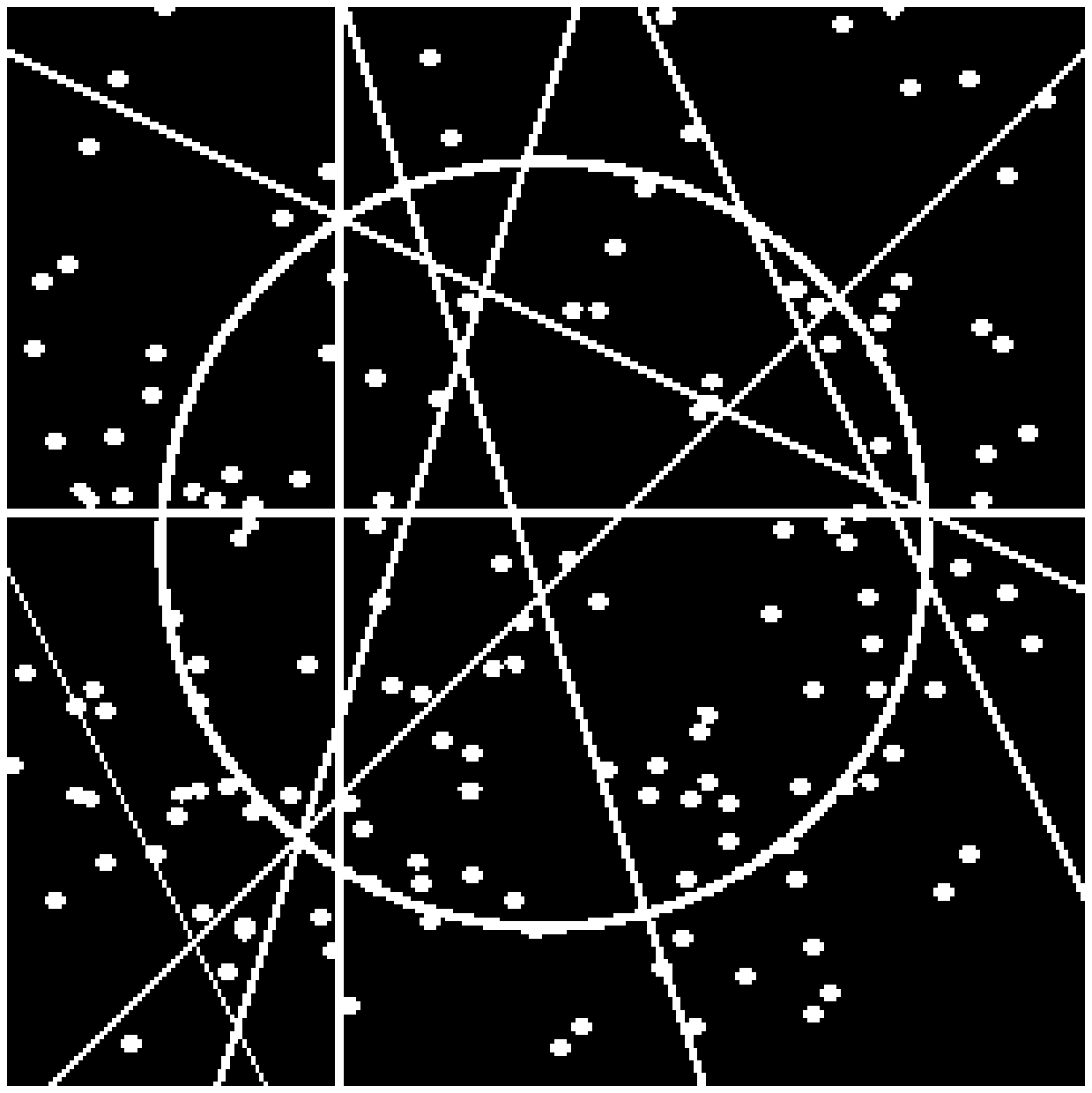}
\put(-305,-10){$f_0$}
\put(-183,-10){$f_1$}
\put(-70,-10){$f = f_0 + f_1$}
\caption{In the image decomposition experiment, we try to recover the curvilinear part $f_0$ and the isotropic part $f_1$ from $f = f_0 + f_1$.}
\label{fig:decomposition_input}
\end{center}
\end{figure}

We assume a given image $f\in\ell^2(\Z^2)$ can be split into a curvilinear part $f_0$ and a part  $f_1$ containing isotropic
blob-like structures in the sense of
\begin{equation}
f = f_0 + f_1.
\end{equation}
To compute a meaningful decomposition, we use a directional transform nicely adapted to curvilinear structures (e.g., shearlets or
curvelets) together with the stationary wavelet transform, and apply the iterative thresholding algorithm from \cite{SED2005}
(see also \cite{FSED2010} and Algorithm \ref{alg:image_decomposition}).
\begin{algorithm}[h]
 \SetAlgoLined
 \KwData{$f$, $\delta_{init}$, $\delta_{min}$, $iterations$}
 \KwResult{$\tilde{f_0},\tilde{f_1}$ }

 $\tilde{f_0} := \tilde{f_1} := 0$\;
 $\delta := \delta_{init}$\;
 $\lambda := (\delta_{min})^{1/(iterations-1)}$\;

  \For{i := 1 \KwTo $iterations$}{
  $f_{res} := f - (f_0 + f_1)$\;
  $\tilde{f_0} := \mathcal{T}_0^{-1}T_\delta\mathcal{T}_0(f_{res} + \tilde{f_0})$\tcp*[r]{Hard thresholding with anistropic transform}
  $\tilde{f_1} := \mathcal{T}_1^{-1}T_\delta\mathcal{T}_1(f_{res} + \tilde{f_1})$\tcp*[r]{Hard thresholding with istropic transform}
  \tcp{To improve the quality of the result, you can add a total variation correction for the directional part}
  $\delta := \lambda\delta$\;
  }

 \caption{Image decomposition via iterative thresholding}
 \label{alg:image_decomposition}
\end{algorithm}
During each iteration, the difference between the input and the current blob-like image is subject to a directional transform, while
the difference between the input and the current curvilinear image is subject to a stationary wavelet transform. Before computing the
inverse transforms, a hard thresholding operator is applied to both sets of coefficients. By gradually decreasing the thresholding
constant, this algorithm iteratively approximates an image close to the curvilinear part of the original data whose coefficients are
sparse in the directional dictionary and an image close to the blob-like part of the original input whose coefficients are sparse with
respect to the stationary wavelet transform. An example can be seen in Figure \ref{fig:decomposition_input}.

To quantitatively measure the performance of a sparse representation scheme, we restrict our experiments to binary images as in Figure
\ref{fig:decomposition_input}. Furthermore, we introduce the operator $B_\delta$, mapping images $g\in \ell^2(\Z^2)$ to binary images with
\[
(B_\delta g)(i,j) = \begin{cases}1 &\text{if }\abs{ g(i,j)} \geq \delta,\\0 &\text{else, }\end{cases}
\]
and the measures (see \cite{KL2012a})
\begin{align*}
Q(\tilde{f_0},\delta) &=  \frac{\|g*f_0-g*(B_\delta\tilde{f_0})\|_2}{\|g*f_0\|_2}, \\
Q_{opt}(\tilde f_0)&= \min_{\delta \in \{0,\ldots,255\}} Q(\tilde{f_0},\delta),
\end{align*}
where $\tilde{f_0}$ is the curvilinear part computed by algorithm \ref{alg:image_inpainting}, $g$ is a discrete two-dimensional Gaussian filter
and $*$ denotes the convolution. This definition can naturally also be used for the blob-like part $\tilde{f_1}$.

\begin{figure}[t]
\begin{center}
\includegraphics[width=2in]{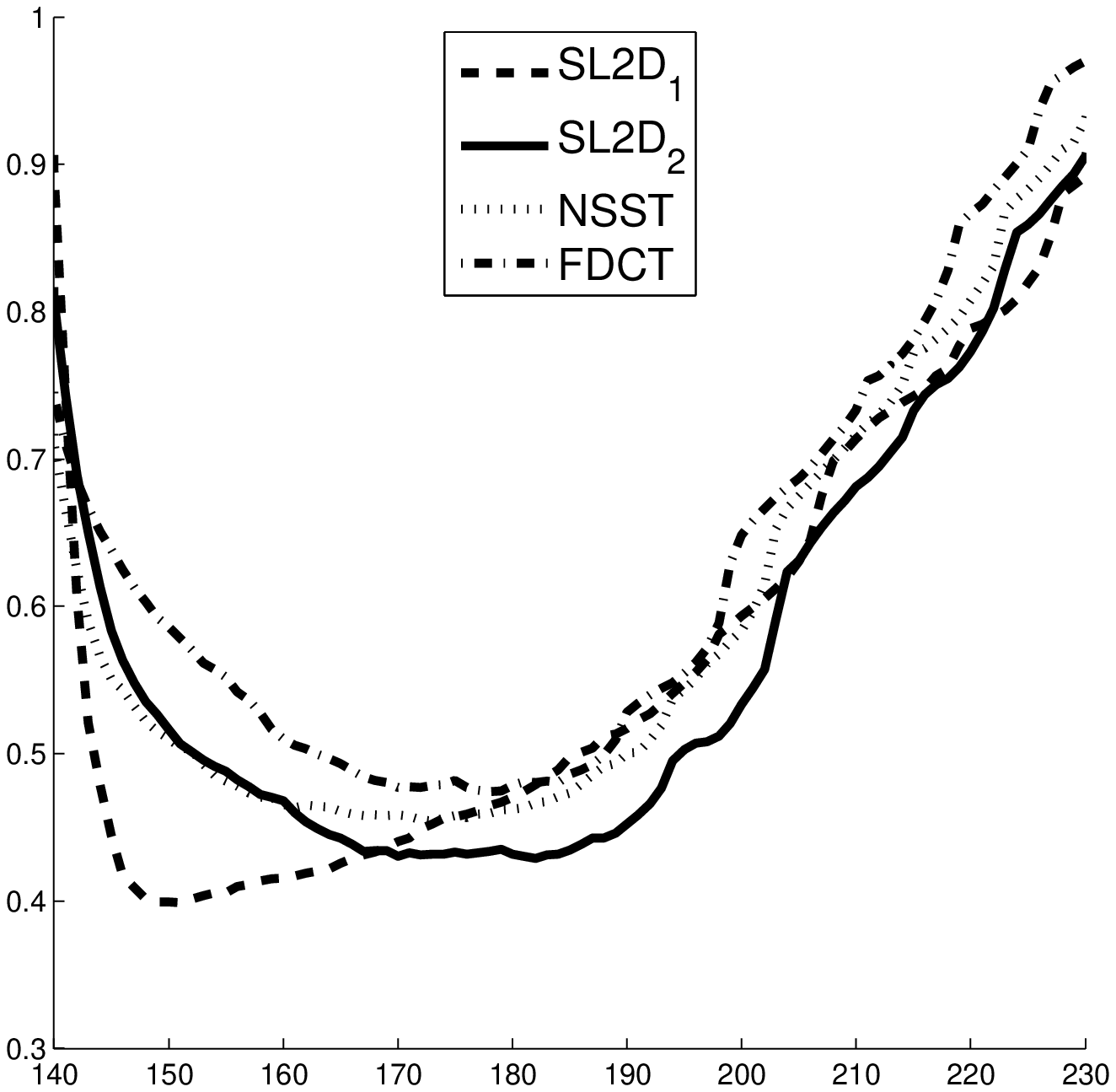}
\hspace{50pt}
\includegraphics[width=2in]{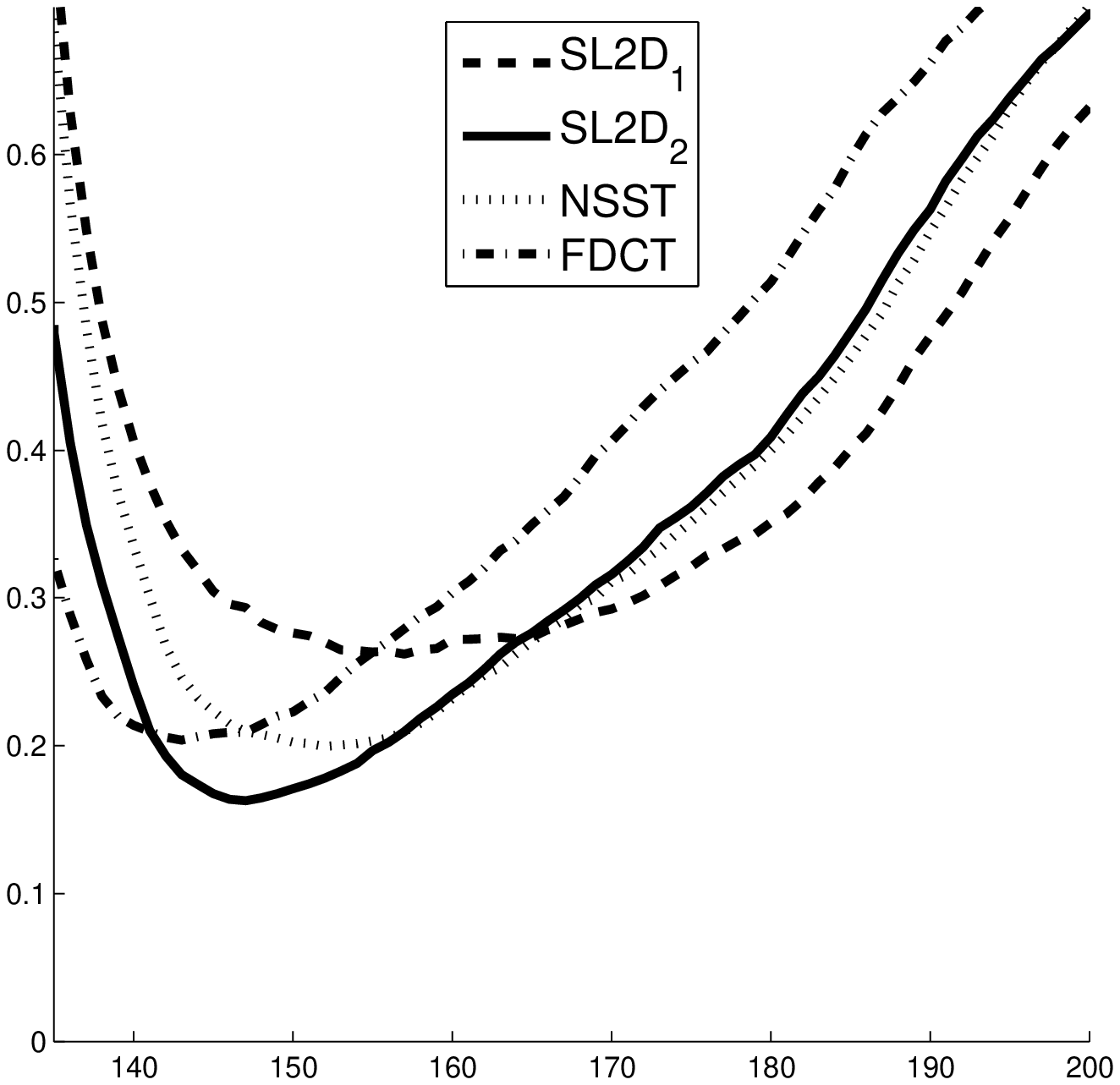}
\put(-290,-15){curvilinear}
\put(-90,-15){blob-like}
\caption{The values $Q(\tilde{f_0},\delta)$ and $Q(\tilde{f_1},\delta)$ as functions of $\delta$.}
\label{fig:results_separation_graphs}
\end{center}
\end{figure}

As input for our experiment, we used the image depicted in Figure \ref{fig:decomposition_input}. We compared the systems $SL2D_1$ and $SL2D_2$
to the directional transforms $NSST$ and $FDCT$. The numerical results together with redundancies and running times are compiled in Figures
\ref{fig:results_separation_values} and \ref{fig:results_separation_graphs}. For a visual comparison of the results, see Figure
\ref{fig:results_decomposition_images} in the Appendix.

\begin{figure}[h]
\begin{center}

\begin{tabular}{|l|c|c|c|c|}
\hline
\multicolumn{5}{|c|}{\textbf{Image Decomposition: Quantitative Results and Running Times}}\\
\hline
&{Points}&{Curves}&Redundancy&Running Time\\\hline
{$SL2D_1$}&$Q_{opt}$ = \textbf{0.3991}&$Q_{opt}$ = 0.2620&25&\unit{57.1}{\second}\\\hline
{$SL2D_2$}&$Q_{opt}$ = 0.4288&$Q_{opt}$ = \textbf{0.1627}&49&\unit{87.8}{\second}\\\hline
{$NSST$}&$Q_{opt}$ = 0.4541&$Q_{opt}$ = 0.1996&49&\unit{203.7}{\second}\\\hline
{$FDCT$}&$Q_{opt}$ = 0.4743&$Q_{opt}$ = 0.2040&2.85&\unit{84.9}{\second}\\\hline
\end{tabular}
\caption{Numerical results of the image decomposition experiment.}
\label{fig:results_separation_values}
\end{center}
\end{figure}

\subsection{Video Denoising} \label{subsec:v_denoising}

Similar to the two-dimensional case, we consider grayscale videos of size 192x192x192 distorted with Gaussian white noise.
For a video $f\in\ell^2(\Z^3)$, we have
\[
f_{noisy}(i,j,k) = f(i,j,k) + e(i,j,k)
\]
where $e(i,j,k) \sim \mathcal{N}(0,\sigma^2)$, and compute
\[
f_{denoised} = \mathcal{T}^{-1}T_\delta \mathcal{T} f_{noisy}
\]
where $\mathcal{T}$ denotes the forward and $\mathcal{T}^{-1}$ the inverse transform associated with a sparse approximation
scheme and $T_\delta$ is the hard thresholding operator already defined in \eqref{eq:hard_tresholding}. Again, our thresholds
will be of the form
\[
 \delta_j = K_j \sigma,
\]
where $j$ iterates the scales of the digital transform. Typically, we will use systems with three scales and choose $K = [K_j]_j
= [3,3,4]$.

In total, we run our experiment with three different videos and noise levels ranging from $\sigma = 10$ to $50$. To quantitatively
compare the performance of ShearLab 3D with the performance of the three-dimensional nonsubsampled shearlet transform ($NSST$) and
the surfacelet transform, we again calculate the peak signal to noise ratio (PSNR) defined in \eqref{eq:PSNR}. For a complete
listing of our numerical results, see Figure \ref{fig:results_video_denoising}. A visual comparison is provided in Figure
\ref{fig:results_video_denoising_images} in the Appendix.

\begin{figure}[t]
\begin{center}
\begin{tabular}{|l|c|c|c|c|c|c|c|c|c|c|}
\hline
\multicolumn{11}{|c|}{\textbf{Video Denoising: Quantitative Results in PSNR and Running Times}}\\
\hline
&\multicolumn{5}{|c|}{{Mobile}} & \multicolumn{5}{|c|}{{Coastguard}}\\
\hline
&$\sigma = 10$&{$20$}&{$30$}&{$40$}&{$50$}&{$10$}&{$20$}&{$30$}&{$40$}&{$50$}\\\hline
\textbf{$SL3D_1$}&35.27&31.32&29.01&27.38&26.14&33.13&29.46&27.51&26.17&25.18\\\hline
\textbf{$SL3D_2$}&\textbf{35.91}&\textbf{32.19}&\textbf{29.99}&\textbf{28.43}&\textbf{27.23}&\textbf{33.81}&\textbf{30.28}&\textbf{28.41}&\textbf{27.14}&\textbf{26.17}\\\hline
\textbf{$NSST$}&34.61&30.83&28.56&26.96&25.73&32.59&29.00&27.05&25.68&24.63\\\hline
\textbf{$SURF$}&32.79&29.96&28.25&27.04&26.11&30.86&28.26&26.87&25.91&25.18\\\hline
\textbf{$SL2D_2$}&32.49&28.58&26.42&24.97&23.90&31.24&27.67&25.79&24.52&23.59\\\hline
\end{tabular}\\
\vspace{0.1cm}
\begin{tabular}{|l|c|c|c|c|c|c|c|c|}
\hline
&\multicolumn{5}{|c|}{{Tennis }}&\multirow{2}{*}{Redundancy}&Running&Running Time\\
\cline{1-6}
&{$10$}&{$20$}&{$30$}&{$40$}&{$50$}&&Time&with CUDA\\\hline
\textbf{$SL3D_1$}&33.76&30.19&28.28&26.94&25.90&76&\unit{146}{\second}&\unit{14.1}{\second}\\\hline
\textbf{$SL3D_2$}&\textbf{34.15}&\textbf{30.67}&\textbf{28.86}&\textbf{27.62}&\textbf{26.67}&292&\unit{610}{\second}&\unit{53.1}{\second}\\\hline
\textbf{$NSST$}&33.02&29.48&27.55&26.20&25.14&208&\unit{512}{\second}&-\\\hline
\textbf{$SURF$}&29.95&27.34&26.02&25.11&24.43&6.4&\unit{33}{\second}&-\\\hline
\textbf{$SL2D_2$}&30.91&27.48&25.55&24.31&23.42&49&\unit{43}{\second}&\unit{26.4}{\second}\\\hline
\end{tabular}

\caption{Numerical results of the video denoising experiment.}
\label{fig:results_video_denoising}
\end{center}
\end{figure}

\subsection{Video Inpainting} \label{subsec:v_inpainting}

Analogous to the two-dimensional case, we consider grayscale videos $f\in\ell^2(\Z^3)$ of size $192\times192\times192$ to be
partially occluded by a binary mask $M \in\{0,1\}^{\Z\times\Z\times\Z}$, i.e.,
\[
 f_{masked}(i,j,k) = f(i,j,k)M(i,j,k)
\]
To fill in the missing gaps, we run Algorithm \ref{alg:image_inpainting} (see \cite{SED2005} and \cite{FSED2010}) with $50$
iterations.

In our video inpainting experiments, we use a random mask with $\unit{80}{\%}$ occlusion and a mask consisting of cubes of random
size with $\unit{5}{\%}$ occlusion. We compare the systems $SL3D_1$ and $SL3D_2$ to the surfacelet transform $SURF$. All numerical
results are compiled in Figure \ref{fig:results_video_inpainting}, and for a visual comparsion we refer to Figure
\ref{fig:results_video_inpainting_images} in the Appendix.

\begin{figure}[h]
\begin{tabular}{|l|c|c|c|c|c|c|c|}
\hline
\multicolumn{8}{|c|}{\textbf{Video Inpainting: Quantitative Results in PSNR and Running Times}}\\
\hline
&\multicolumn{2}{|c|}{{Mobile}}&\multicolumn{2}{|c|}{{Coastguard}}&\multirow{2}{*}{Redundancy}&Running&Running Time\\
\cline{1-5}
&Rand&Squares&Rand&Squares&&Time&with CUDA\\\hline
\textbf{$SL3D_1$}&28.15&31.65&26.72&29.05&76&-&\unit{837}{\second}\\\hline
\textbf{$SL3D_2$}&\textbf{29.98}&\textbf{32.87}&\textbf{28.09}&{30.22}&292&-&\unit{3101}{\second}\\\hline
\textbf{$SURF$}&22.40&28.24&21.50&27.1&6.4&\unit{852}{\second}&-\\\hline
\textbf{$SL2D_2$}&23.63&31.61&23.42&\textbf{30.80}&49&-&\unit{1350}{\second}\\\hline
\end{tabular}
\caption{Numerical results of the video inpainting experiment.}
\label{fig:results_video_inpainting}
\end{figure}

\subsection{Discussion} \label{subsec:discussion}

In all experiments, the sparse approximations provided by ShearLab 3D yield the best results with respect to the applied quantitative
measures, except for some cases where our algorithm is slightly outperformed by the nonsubsampled shearlet transform ($NSST$) (see,
for example, Figure \ref{fig:results_denoising_2d}). However, the $NSST$ has significantly worse running times in most experiments,
in particular, if a CUDA capable GPU is available (cf, e.g., Figure \ref{fig:results_inpainting_2d}). In computationally heavy tasks
like image or video inpainting, applying CUDA can lead to a significantly increased speed (up to a factor $10$ in Figures
\ref{fig:results_inpainting_2d} and \ref{fig:results_video_denoising}).

The main goal of our experiments was not to argue that the digital shearlet transform implemented in ShearLab 3D is specifically adapted
to a certain task like image denoising, but to compare its applicability to other, similar transforms. That being said, we would like to
mention that our video denoising results are only marginally beaten by the BM3D algorithm \cite{MBFE2012} which represents -- to our
knowledge -- the current state of the art (PSNR values for the coastguard sequence, $144\times176\times300$, $\sigma = 30$, BM3D: $29.69$, $SL3D_2$: $29.54$).


\bibliographystyle{ACM-Reference-Format-Journals}
\bibliography{ShearLab3D}

\clearpage
\begin{appendix}

 \begin{figure}[b!]
\begin{center}
\includegraphics[width=1.8in]{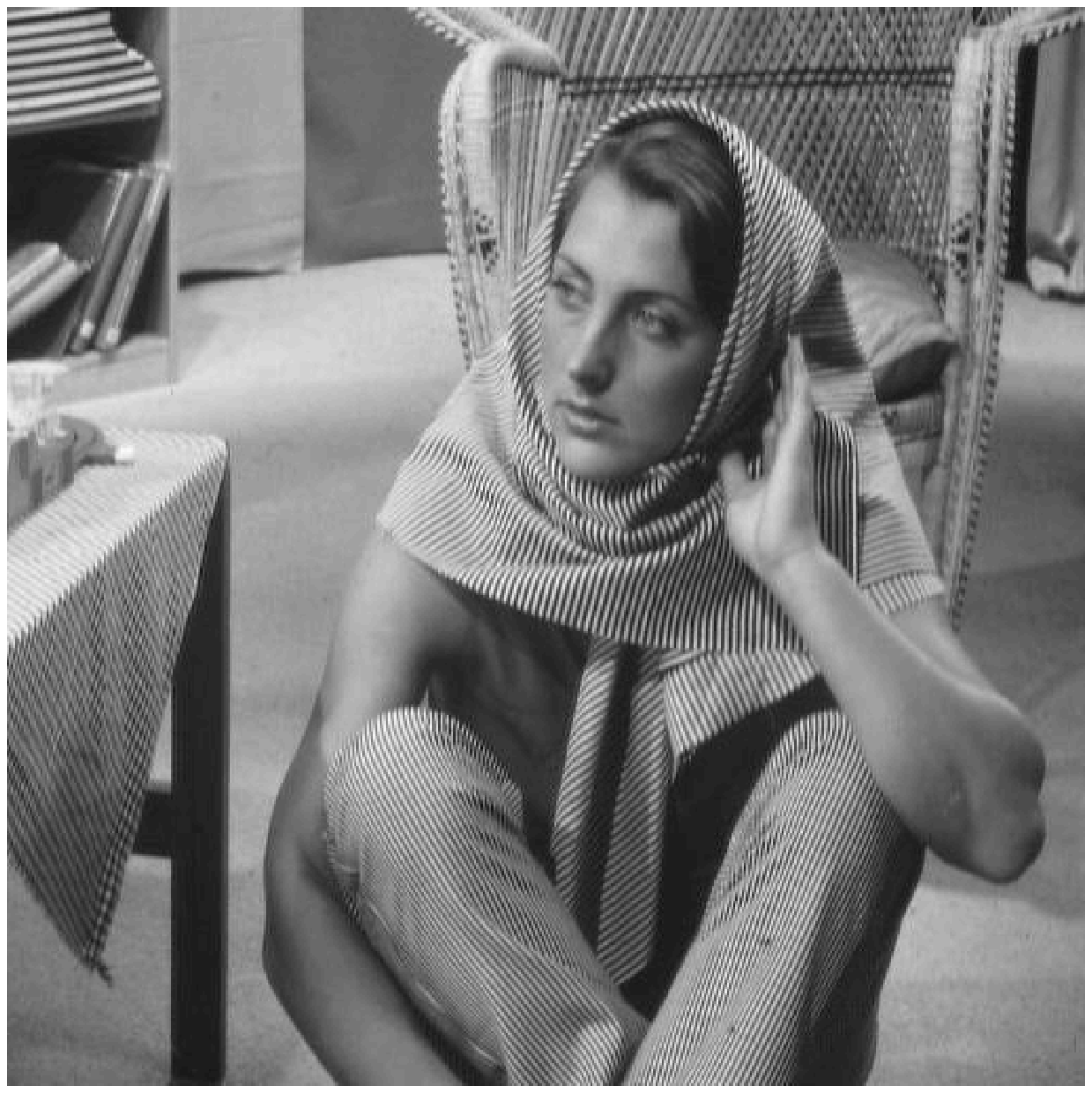}
\includegraphics[width=1.8in]{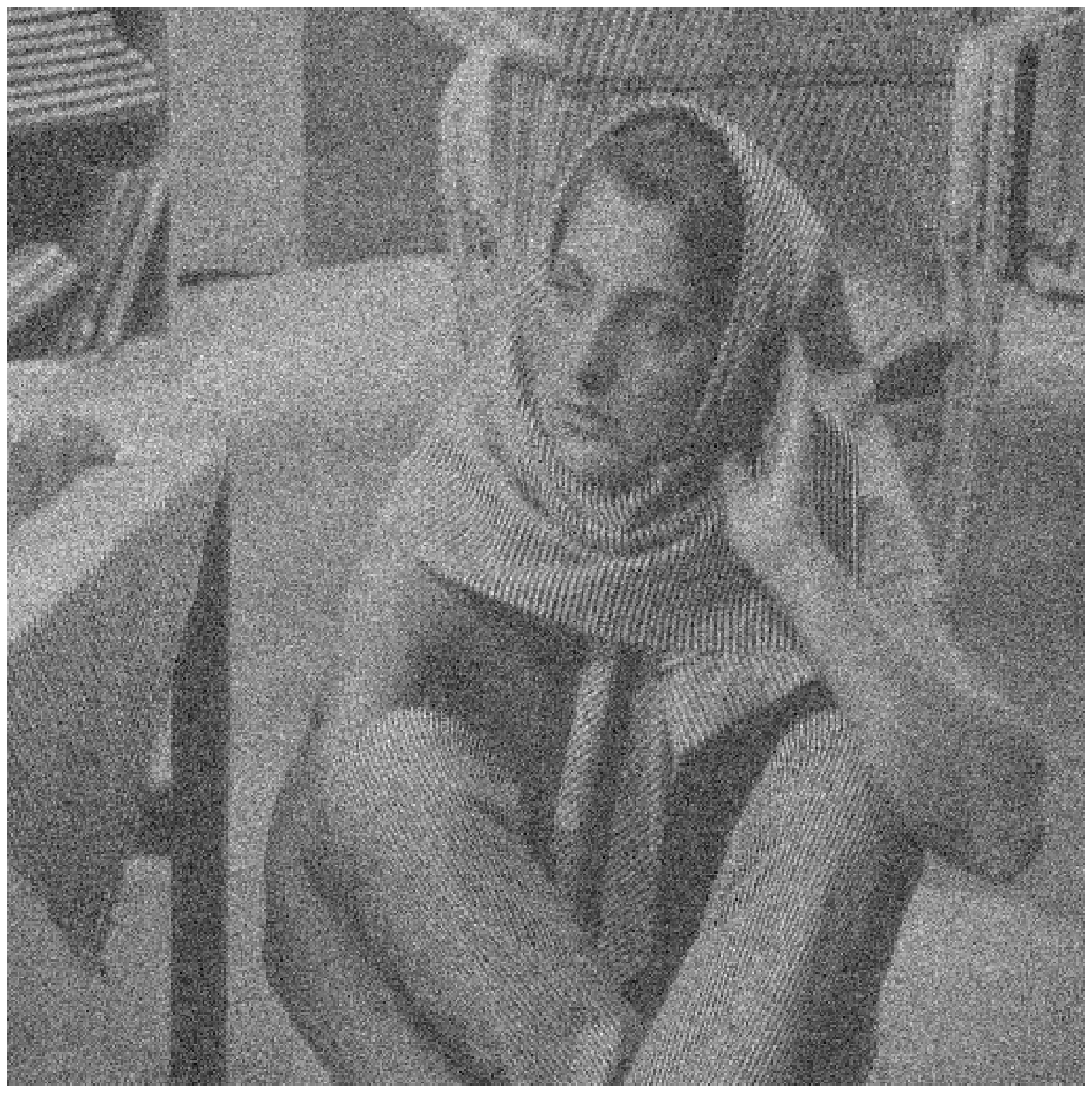}
\includegraphics[width=1.8in]{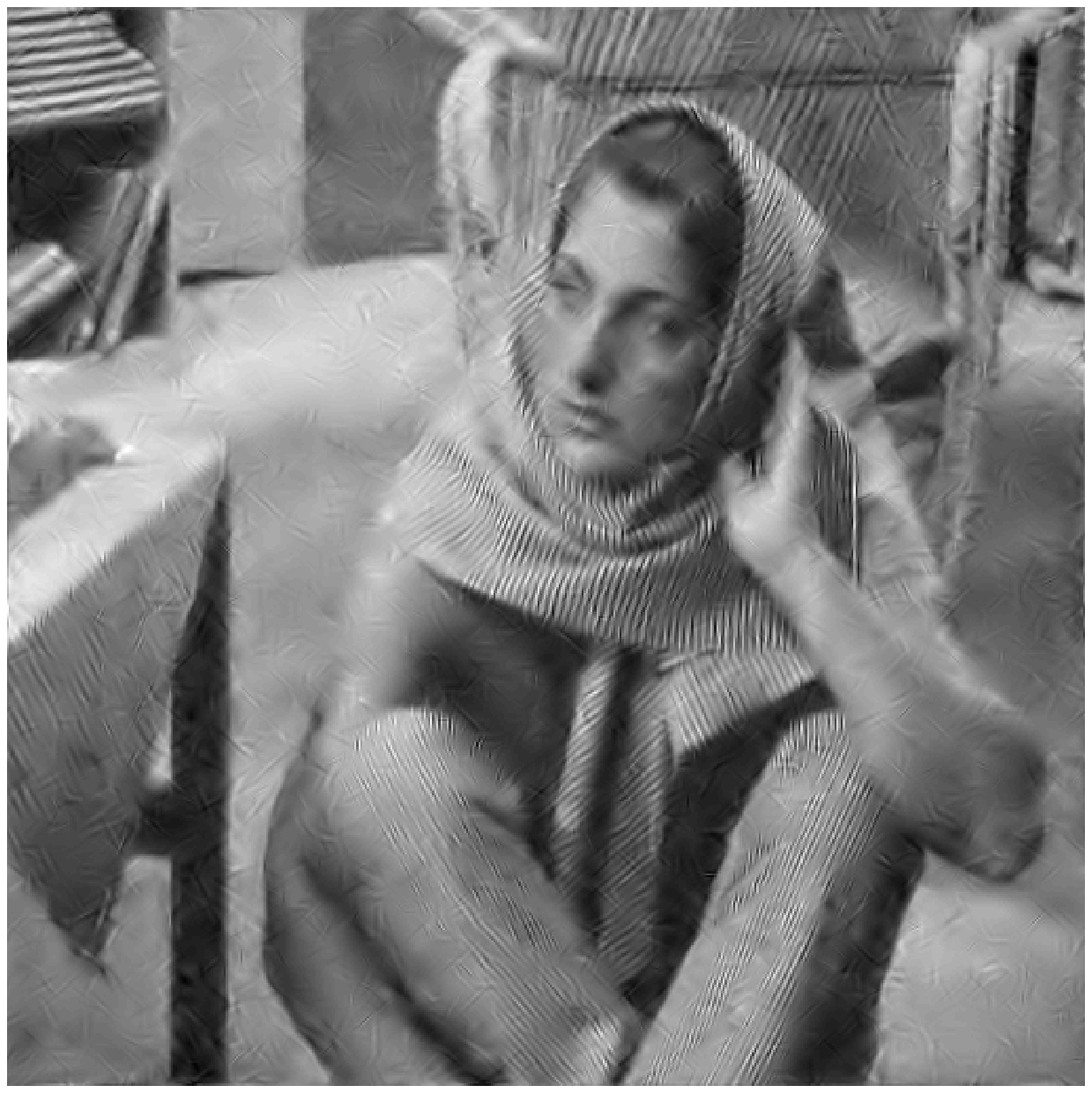}
\put(-350,-15){original}
\put(-265,-15){noisy ($\sigma = 40, \text{PSNR} = 16.06$)}
\put(-115,-15){$SL2D_1$ ($\text{PSNR} = 25.40$)}
\\
\vspace{20pt}
\includegraphics[width=1.8in]{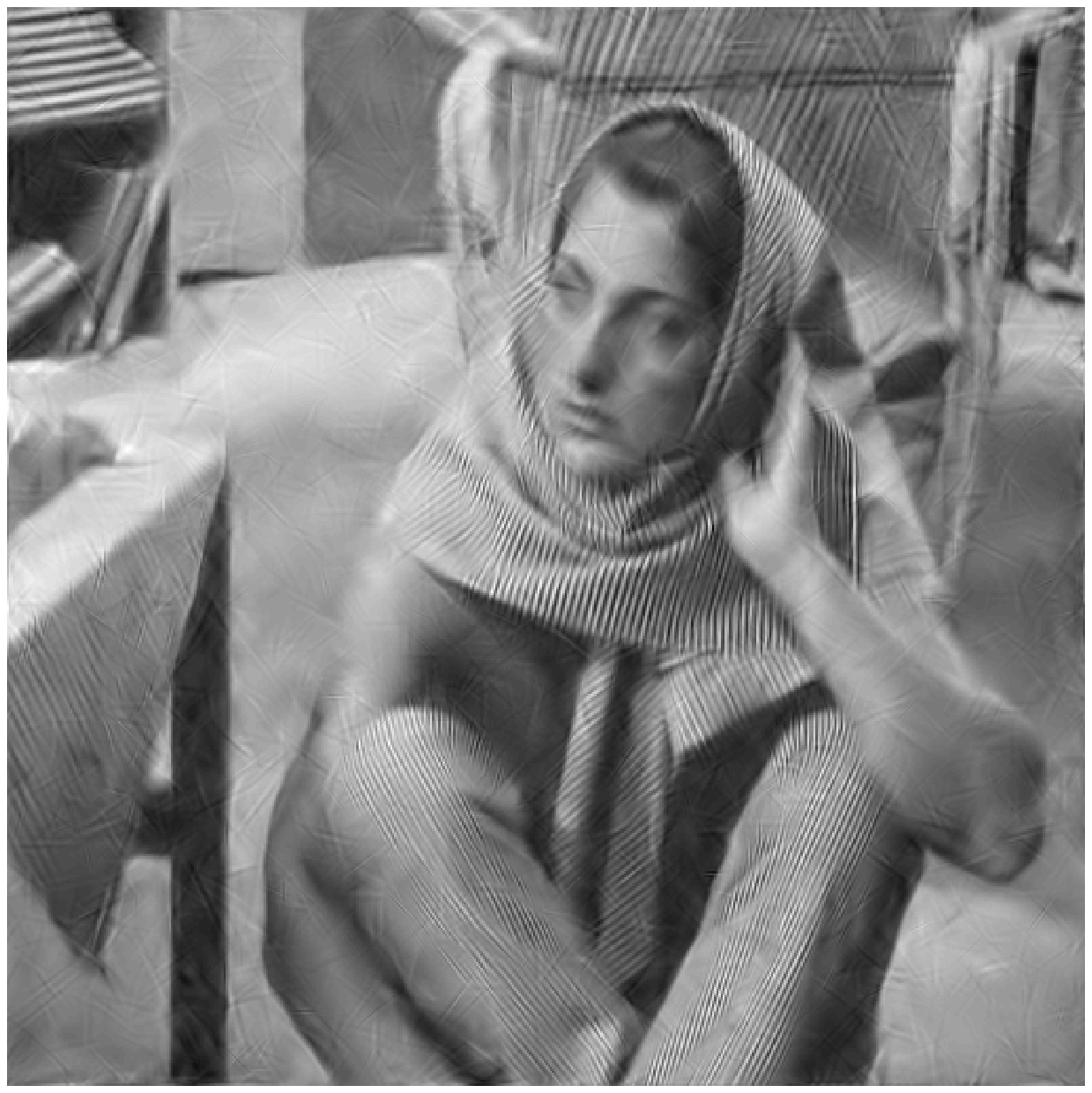}
\includegraphics[width=1.8in]{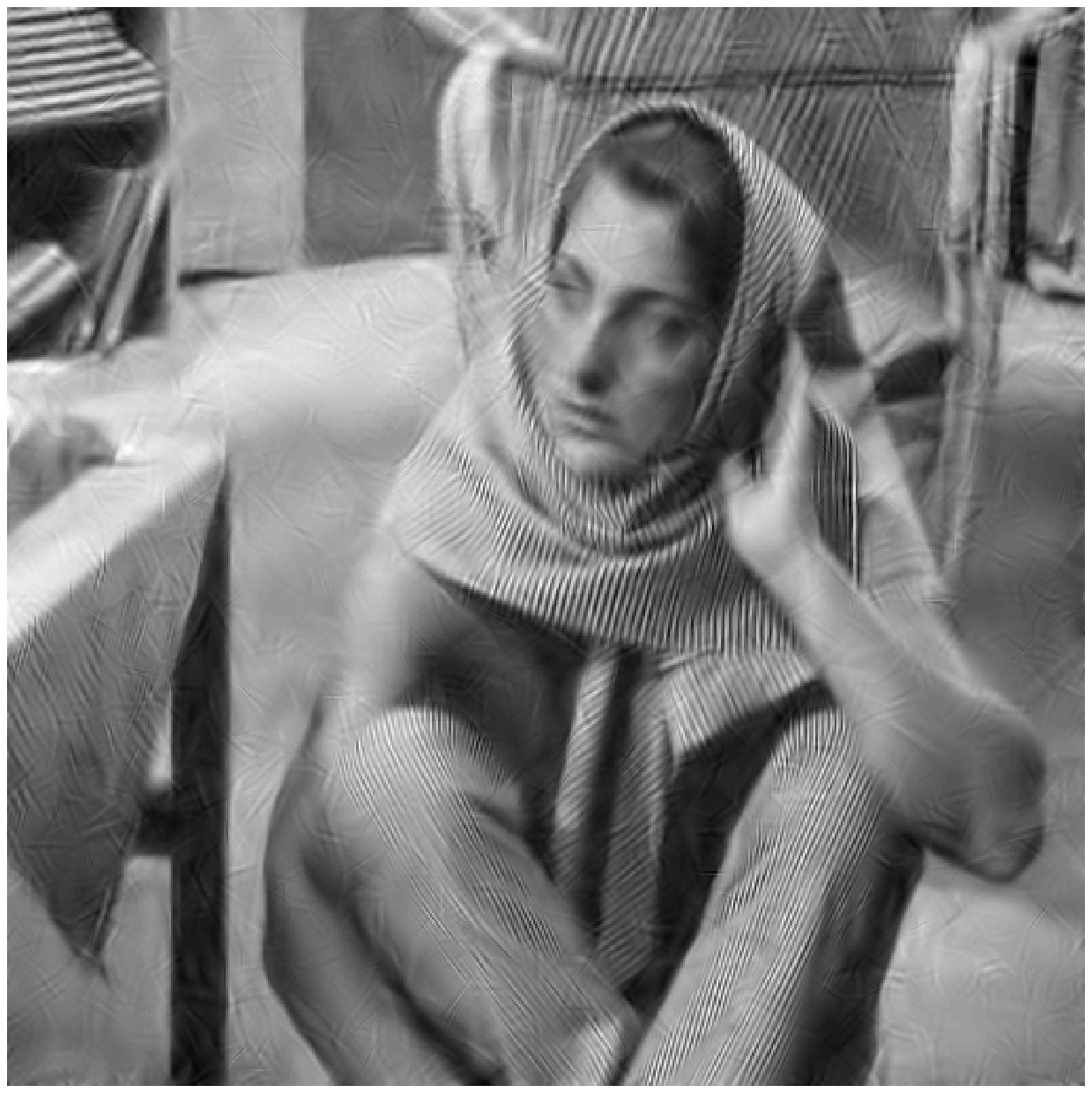}
\includegraphics[width=1.8in]{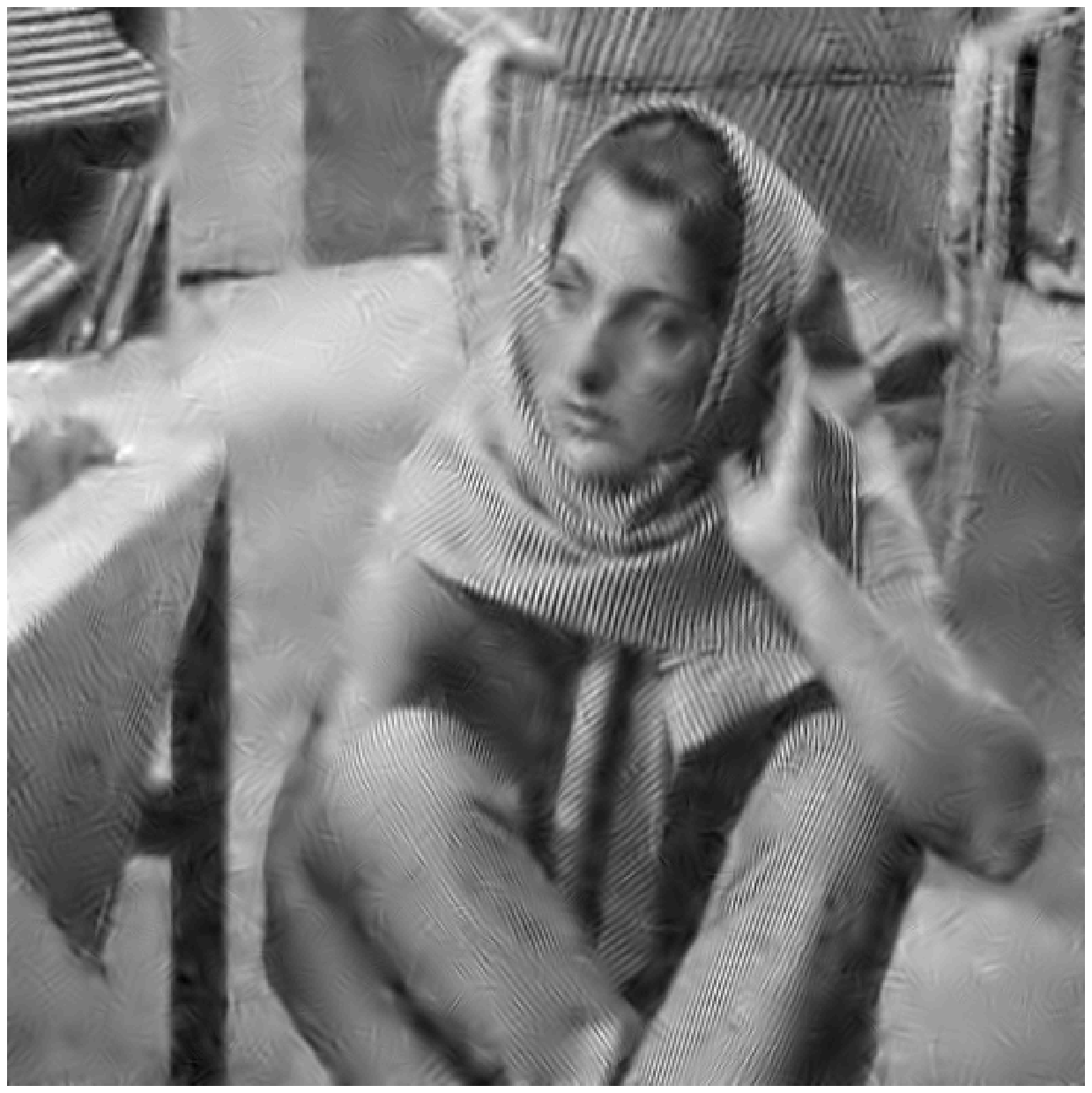}
\put(-385,-15){$SL2D_2$ ($\text{PSNR} = 26.28$)}
\put(-247,-15){$NSST$ ($\text{PSNR} = 26.20$)}
\put(-115,-15){$NSCT$ ($\text{PSNR} = 25.89$)}\\
\vspace{20pt}
\includegraphics[width=1.8in]{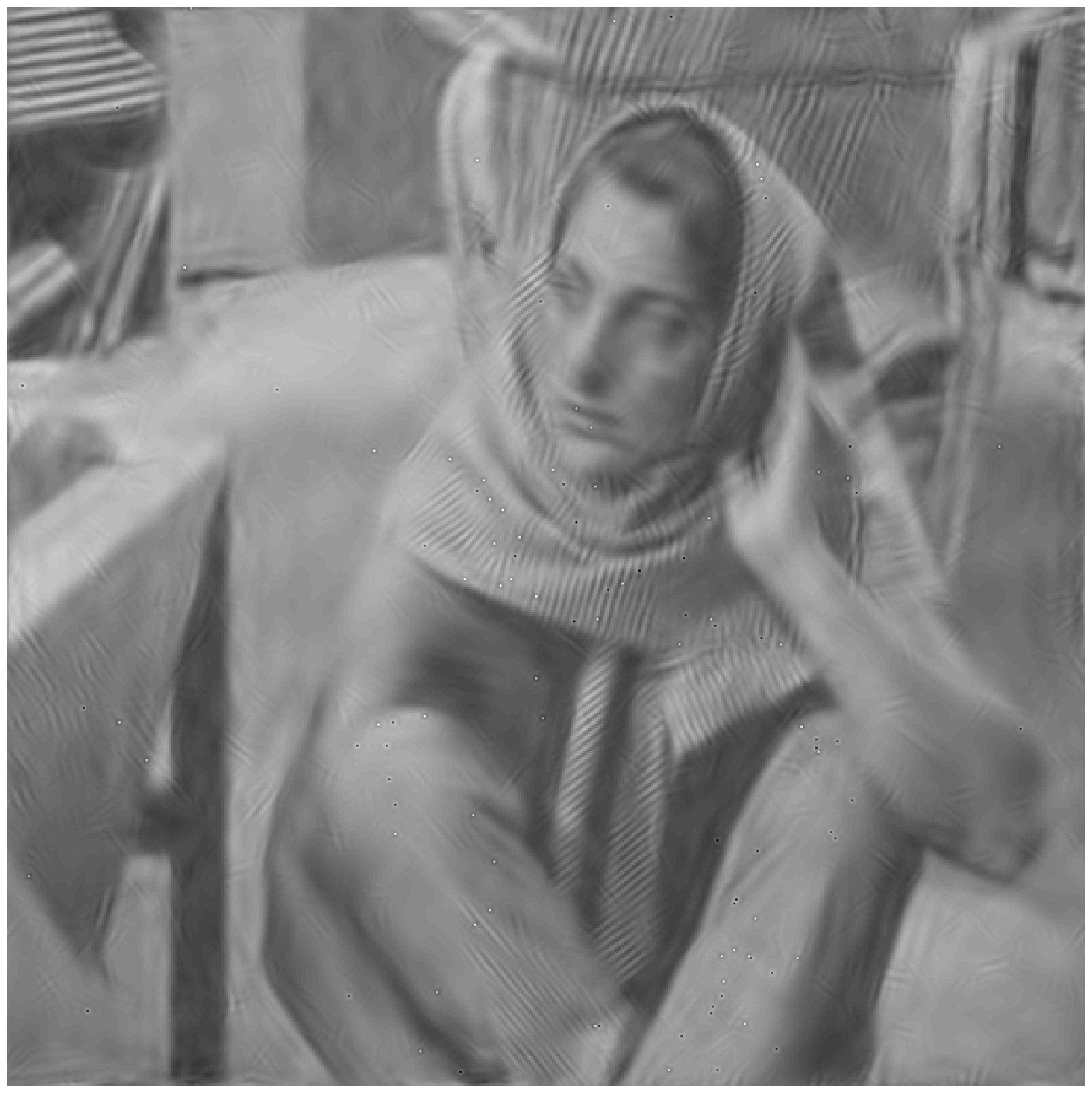}
\includegraphics[width=1.8in]{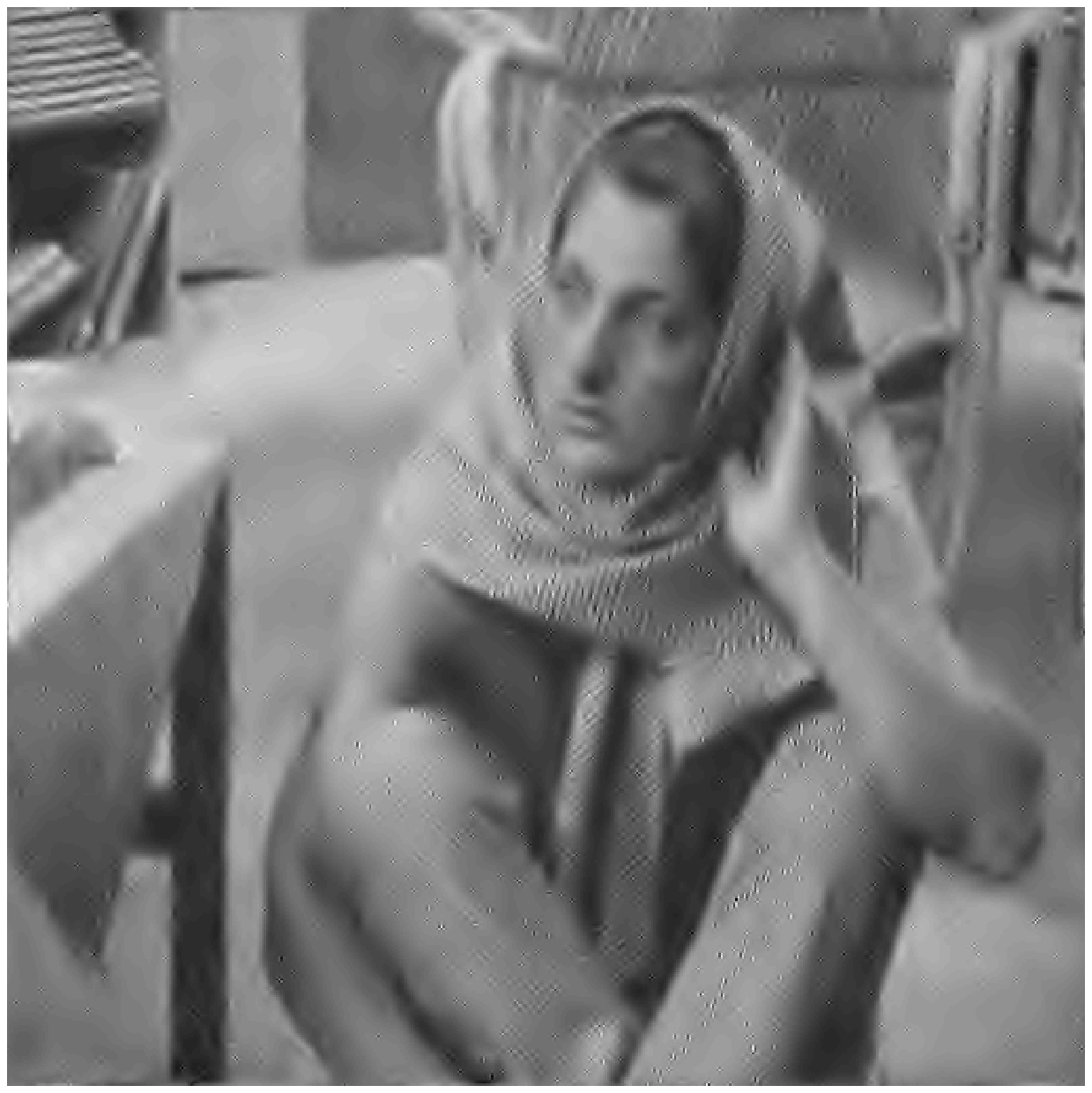}
\put(-247,-15){$FDCT$ ($\text{PSNR} = 23.87$)}
\put(-115,-15){$SWT$ ($\text{PSNR} = 23.55$)}
\end{center}
\caption{The Barbara image is distorted with Gaussian white noise with $\sigma = 40$ and denoised using various sparse approximation schemes.}
\label{fig:results_denoising_2d_images}
\end{figure}

\begin{figure}[t!]
\begin{center}
\includegraphics[width=1.8in]{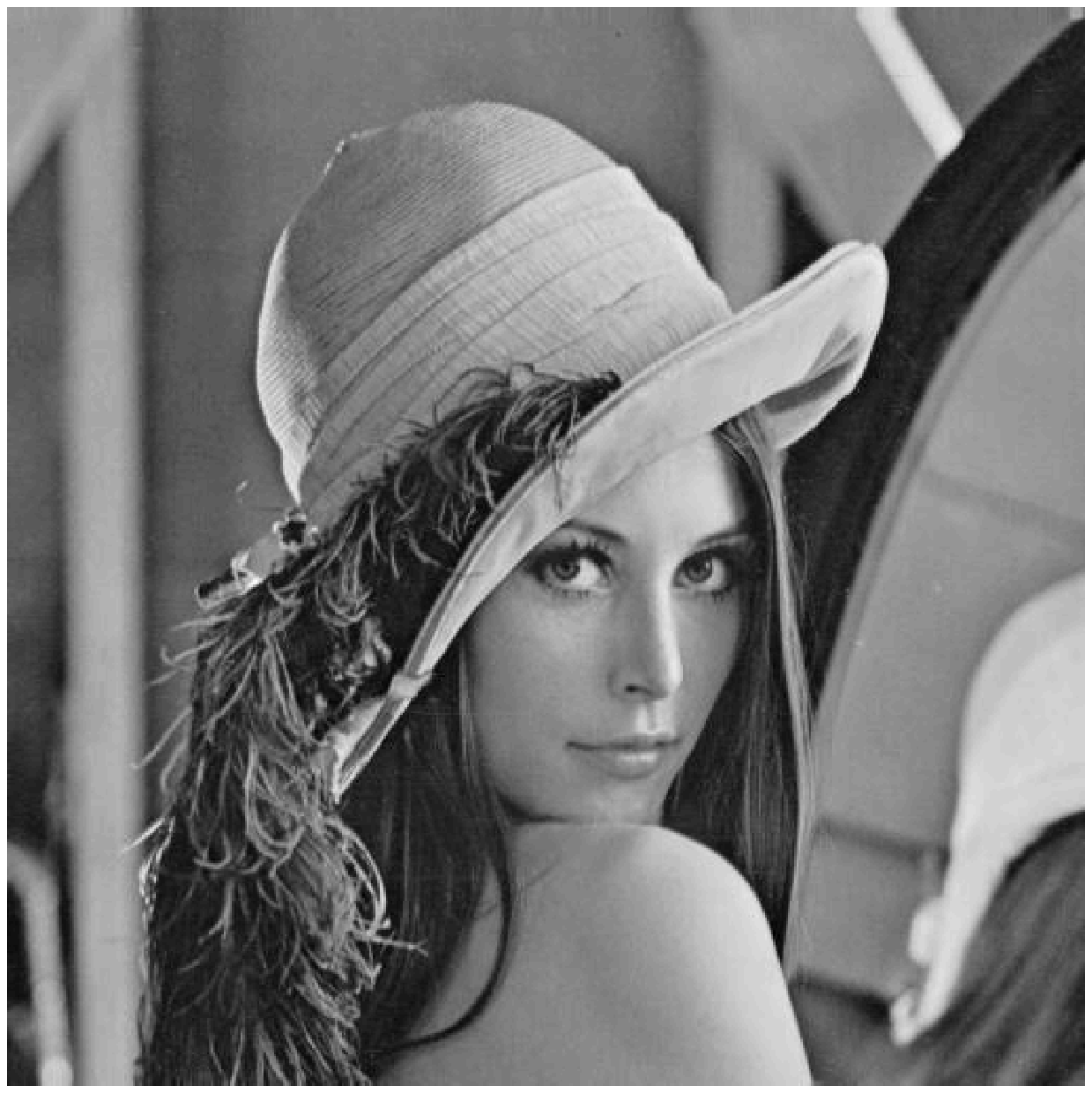}
\includegraphics[width=1.8in]{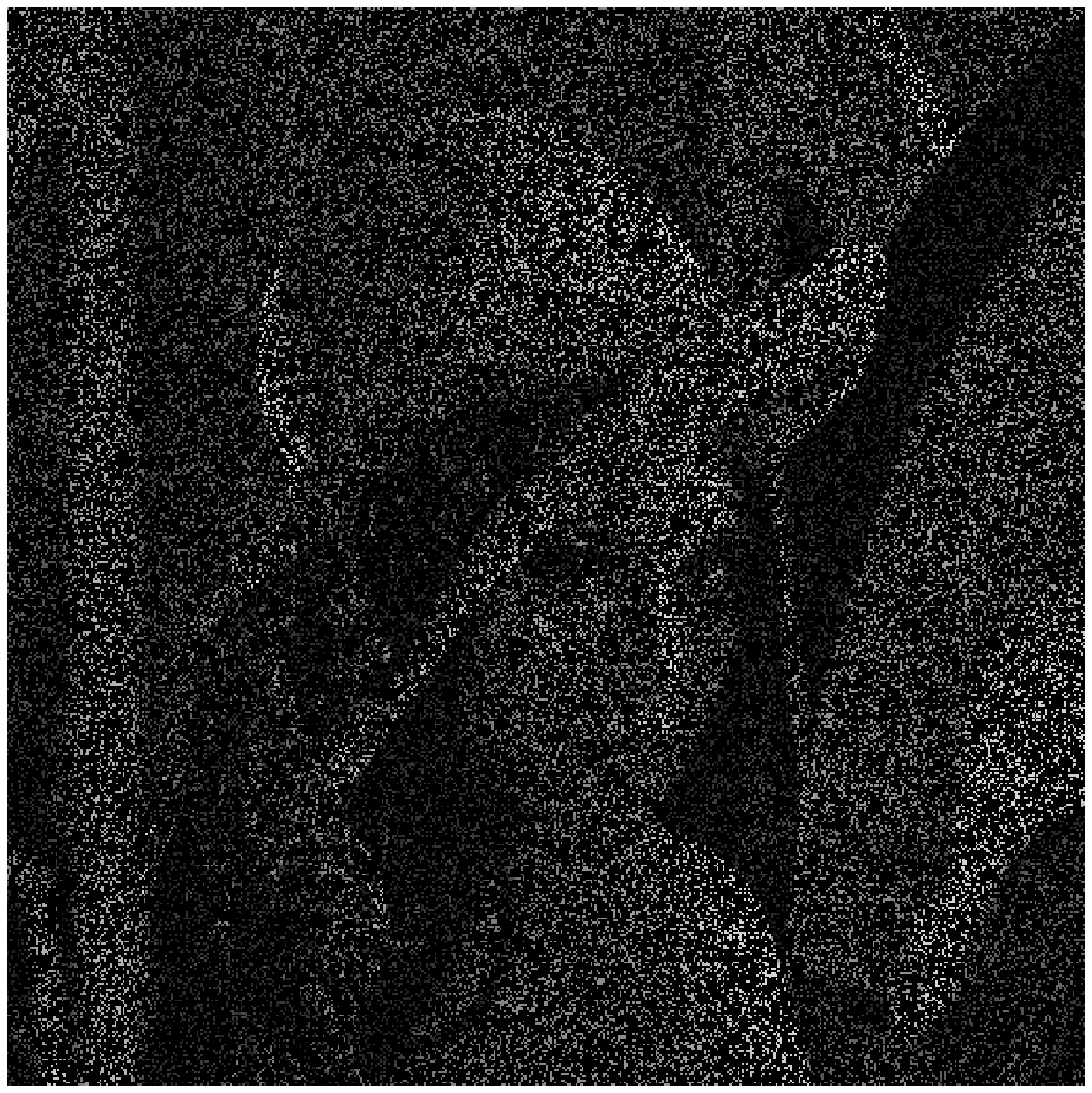}
\includegraphics[width=1.8in]{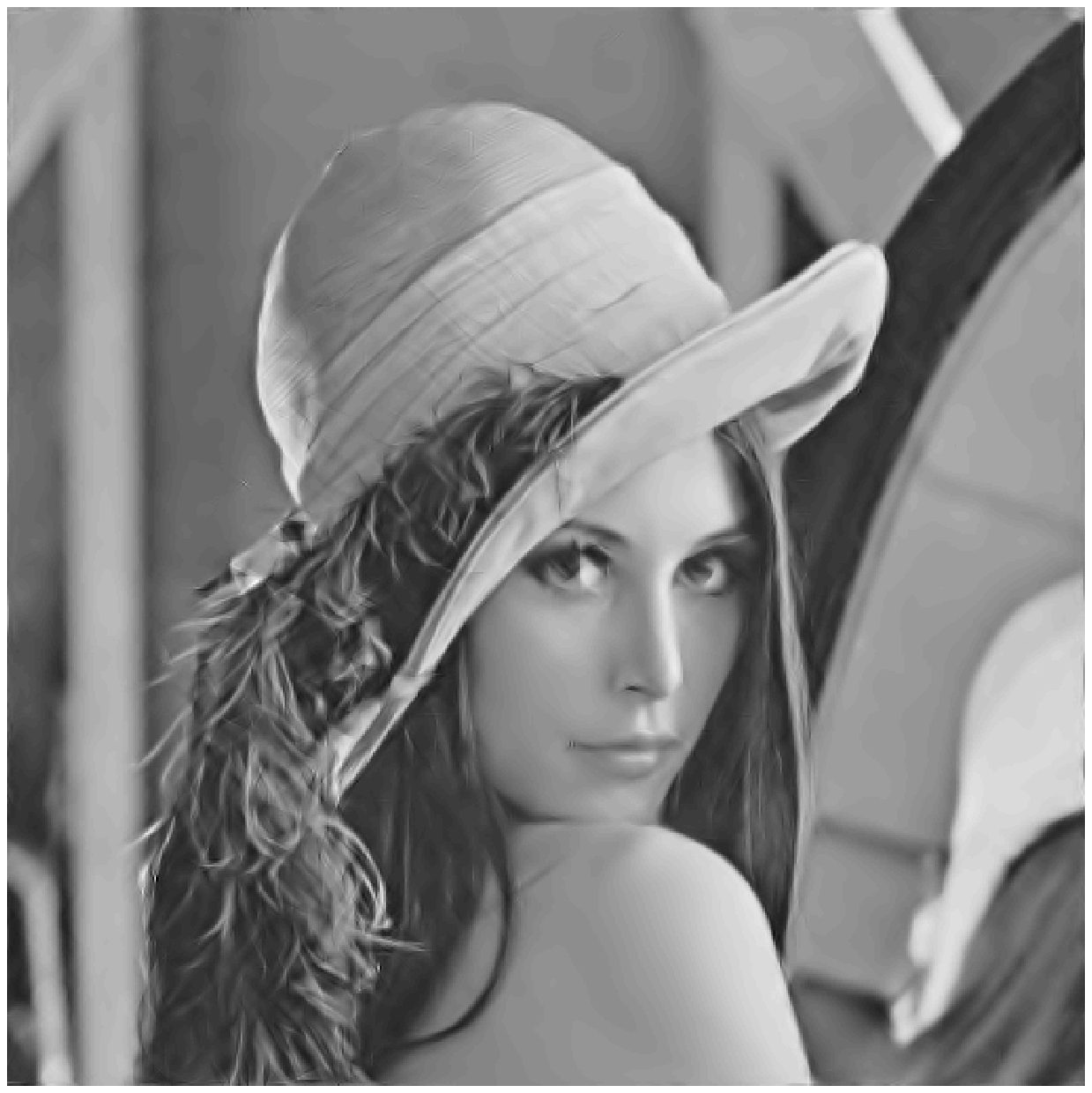}
\put(-350,-15){original}
\put(-255,-15){occluded (\unit{80}{\%} missing)}
\put(-115,-15){$SL2D_1$ ($\text{PSNR} = 32.17$)}
\\
\vspace{20pt}
\includegraphics[width=1.8in]{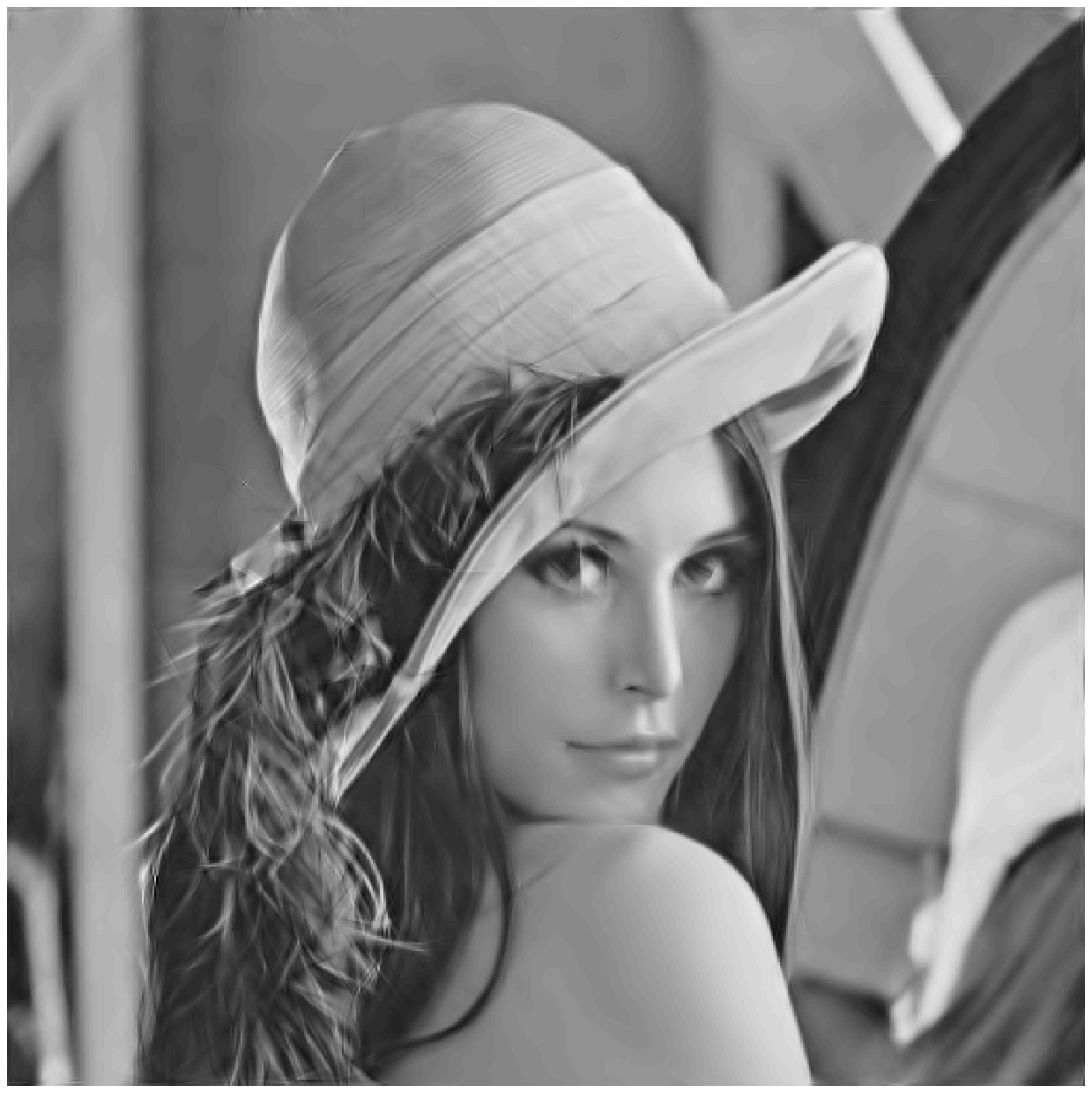}
\includegraphics[width=1.8in]{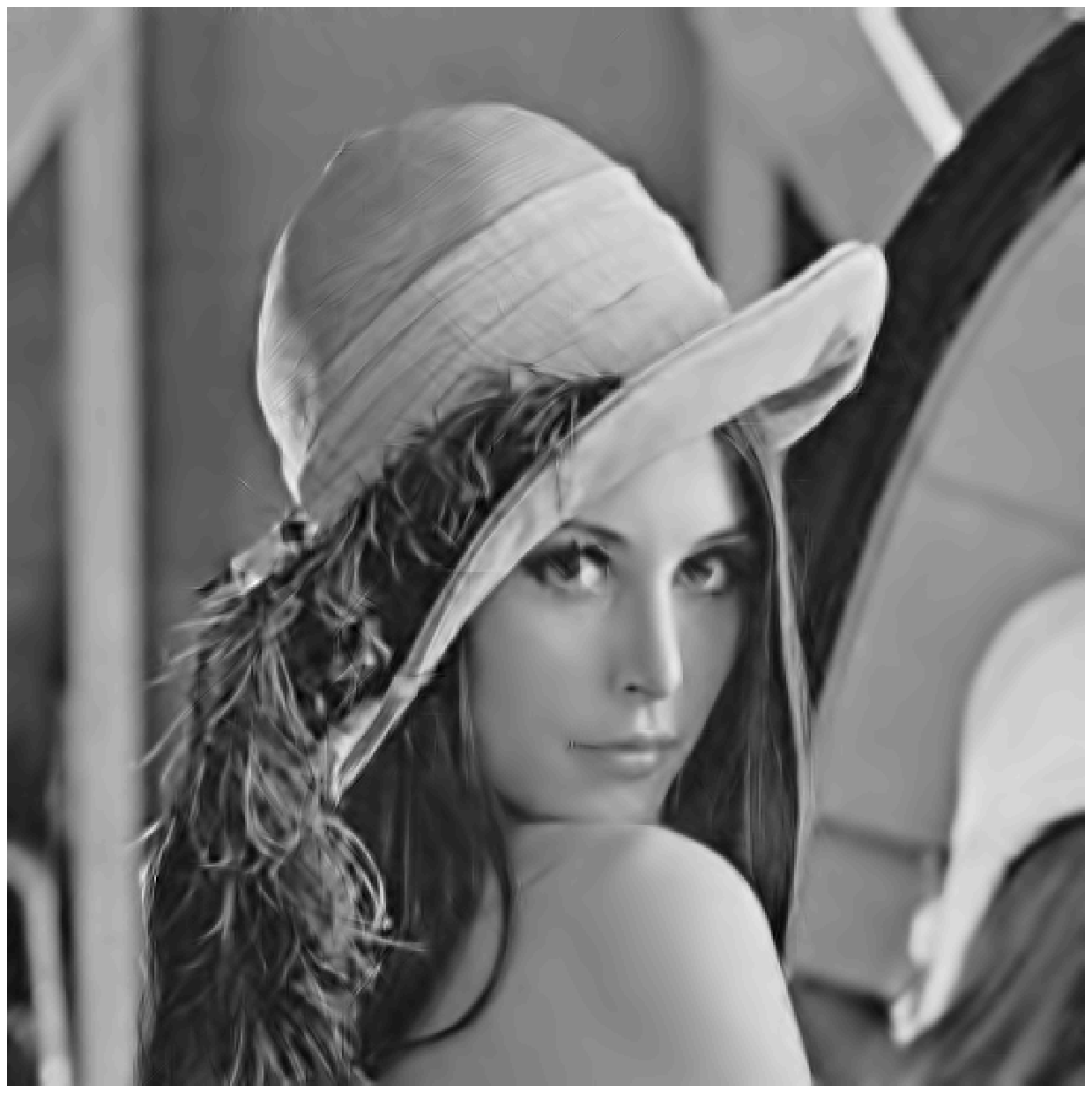}
\includegraphics[width=1.8in]{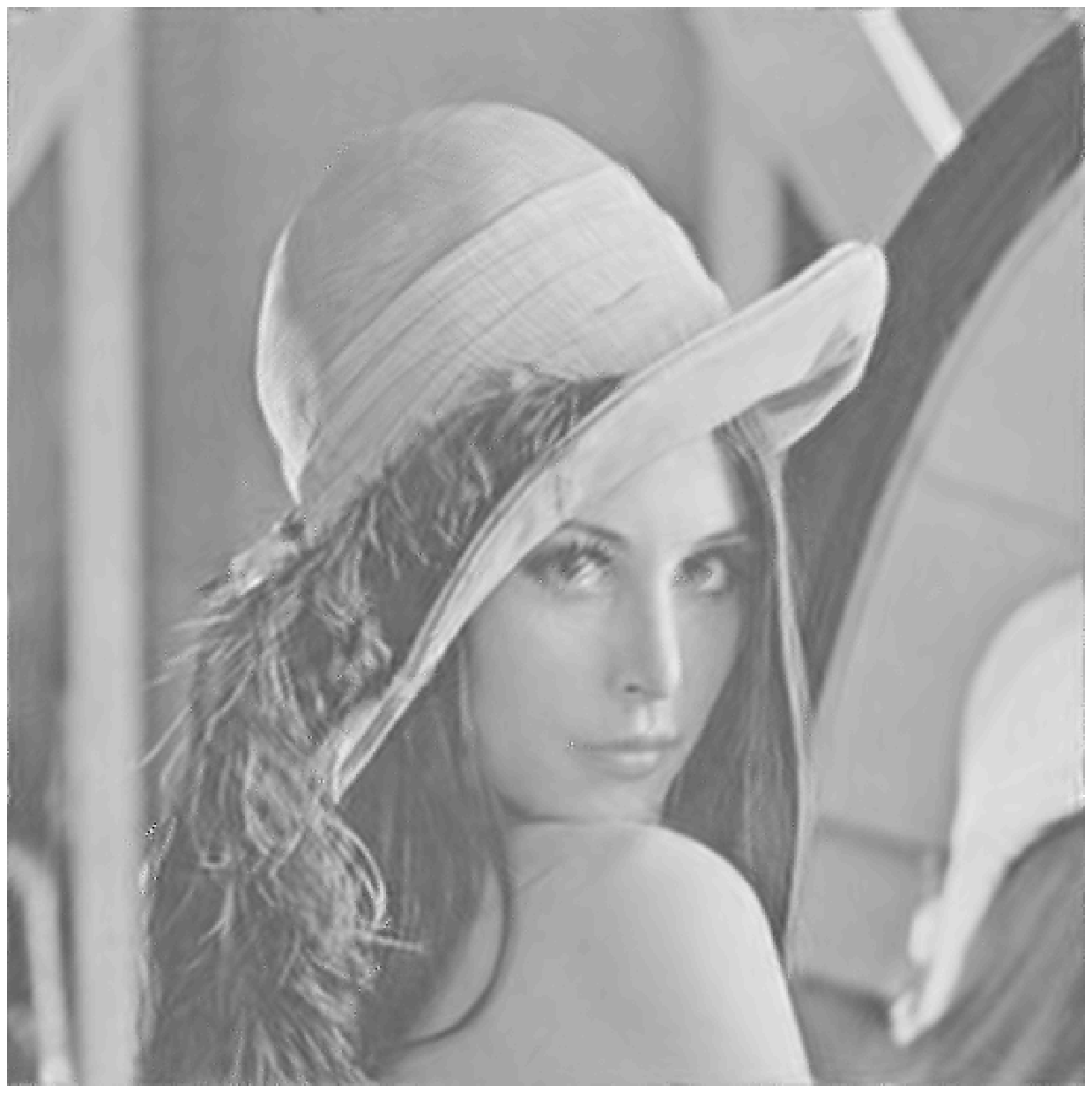}
\put(-385,-15){$SL2D_2$ ($\text{PSNR} = 32.08$)}
\put(-247,-15){$NSST$ ($\text{PSNR} = 32.31$)}
\put(-115,-15){$FDCT$ ($\text{PSNR} = 30.40$)}\\
\vspace{20pt}
\includegraphics[width=1.8in]{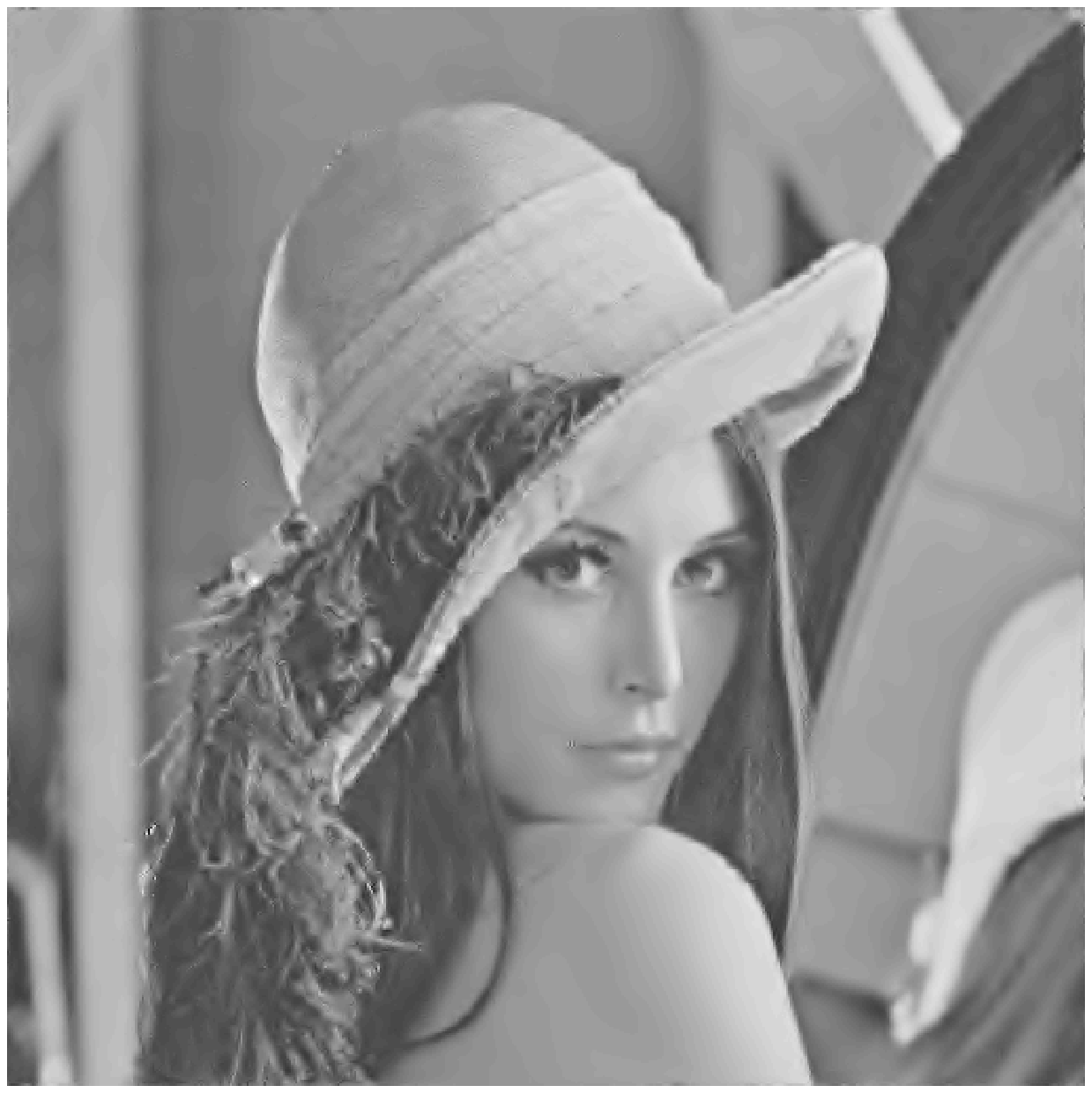}
\put(-115,-15){$SWT$ ($\text{PSNR} = 30.46$)}
\end{center}
\caption{The Lenna image is occluded with a random binary mask and denoised using various sparse approximation schemes.}
\label{fig:results_inpainting_2d_images}
\end{figure}

\begin{figure}[t!]
\begin{center}
\includegraphics[width=1.8in]{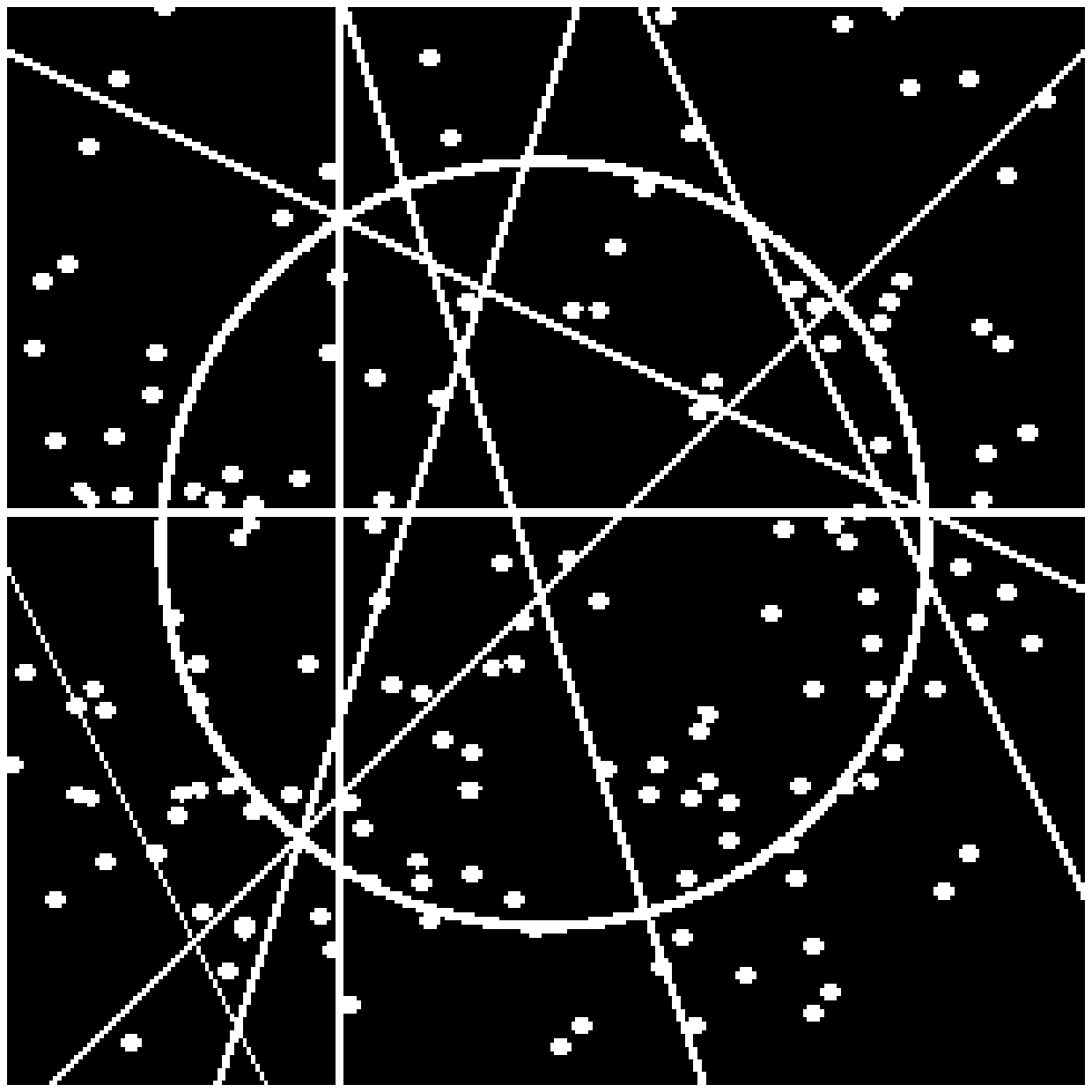}
\includegraphics[width=1.8in]{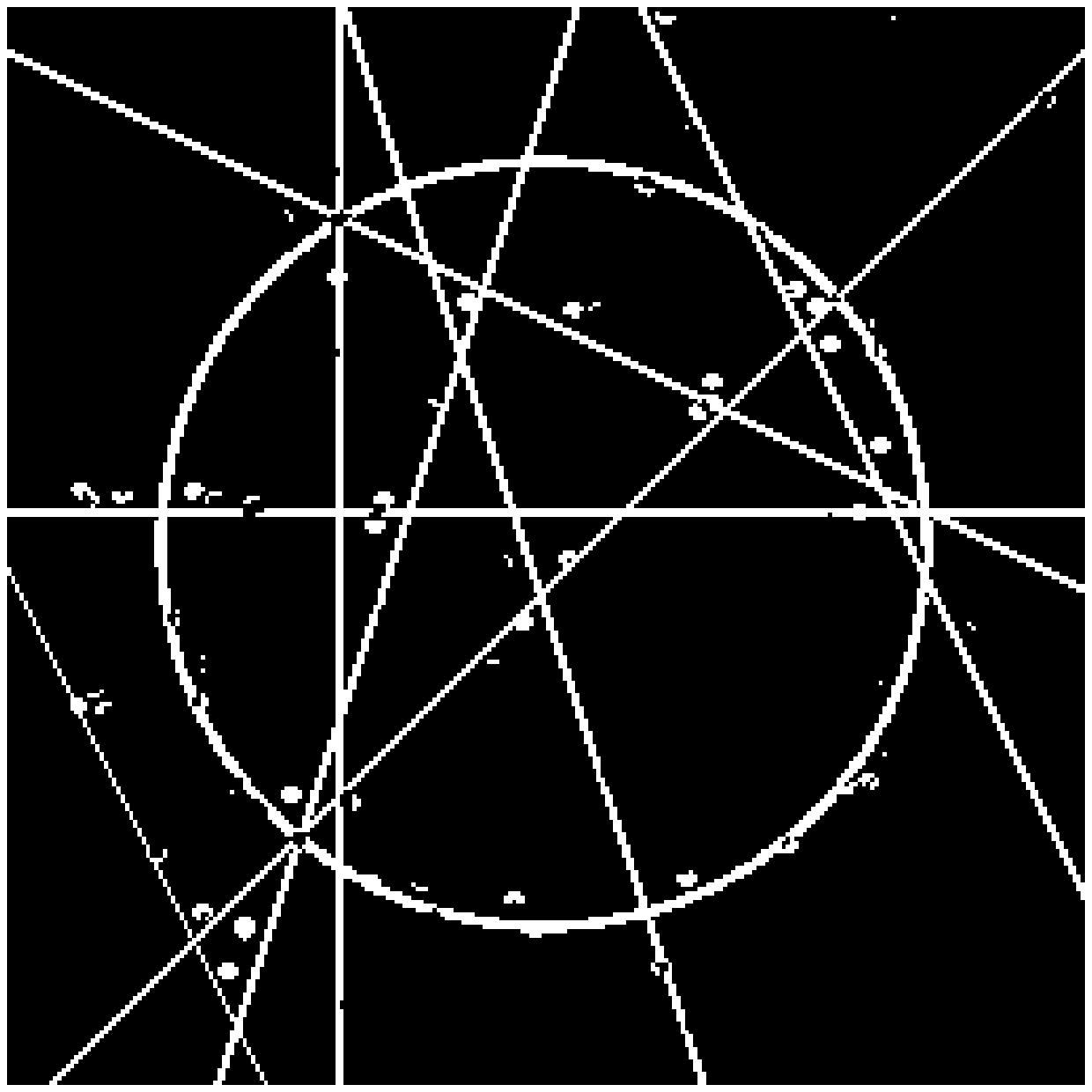}
\includegraphics[width=1.8in]{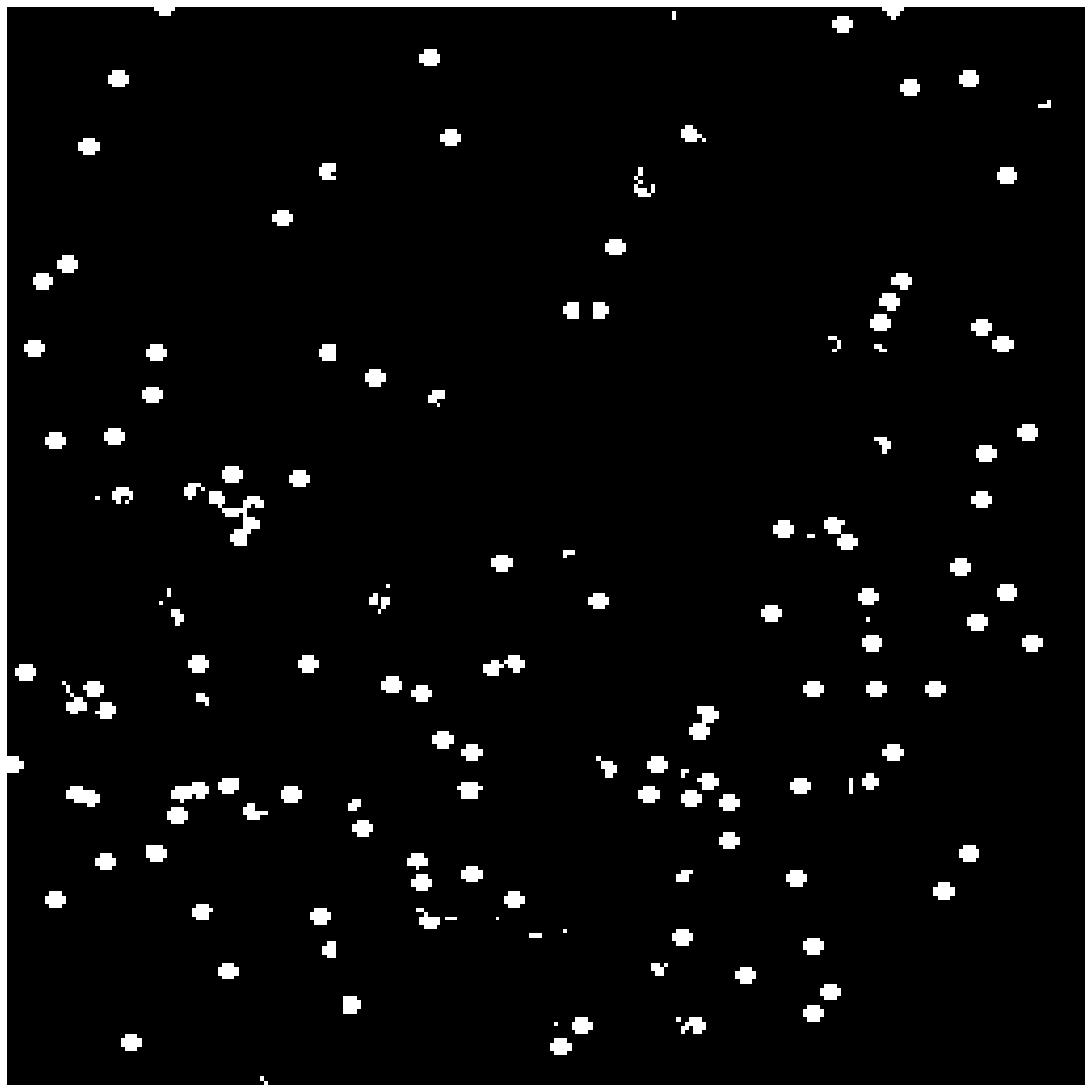}
\put(-350,-15){original}
\put(-263,-15){$SL2D_1$ (curves, $Q_{opt} = 0.26$)}
\put(-130,-15){$SL2D_1$ (points, $Q_{opt} = 0.40$)}
\\
\vspace{20pt}
\includegraphics[width=1.8in]{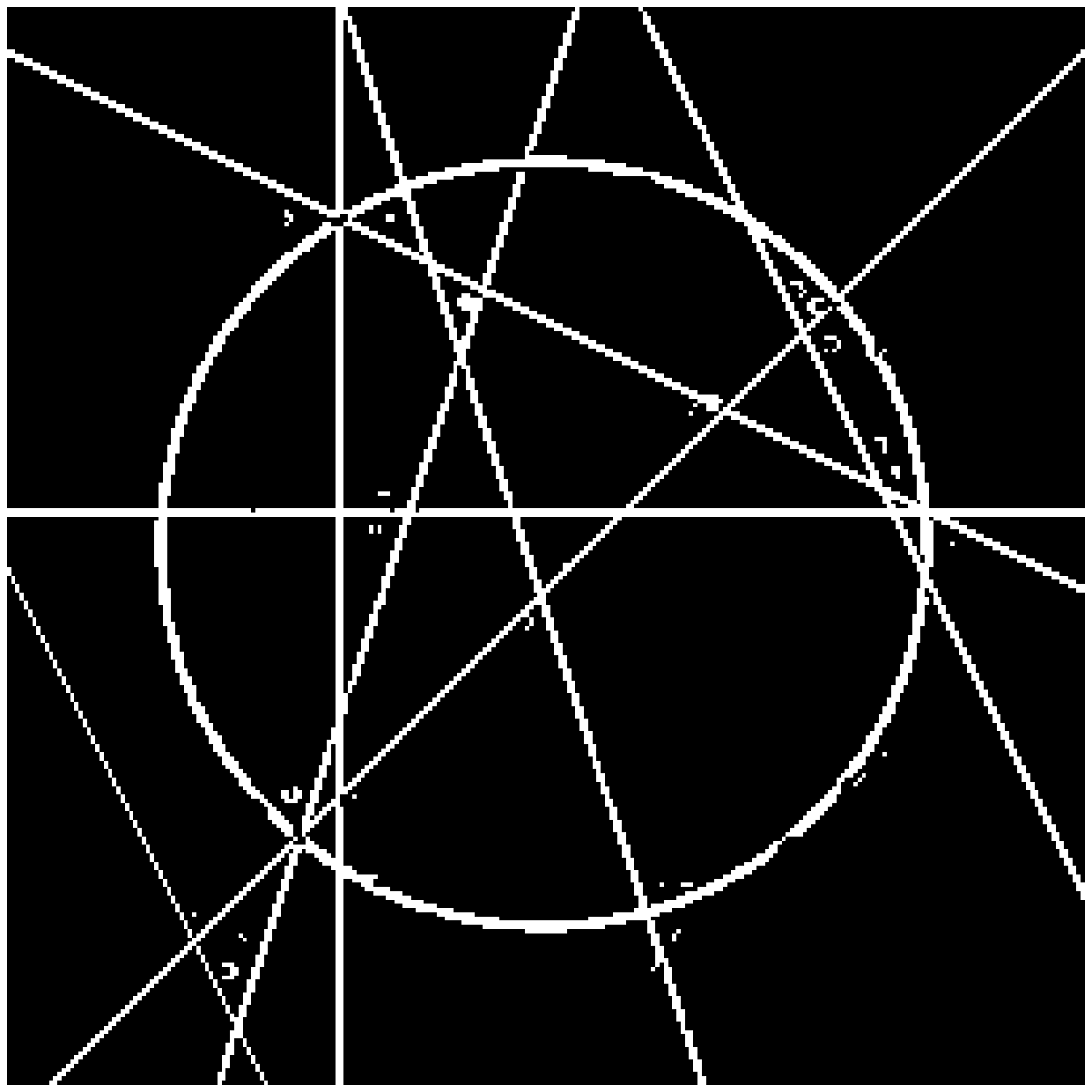}
\includegraphics[width=1.8in]{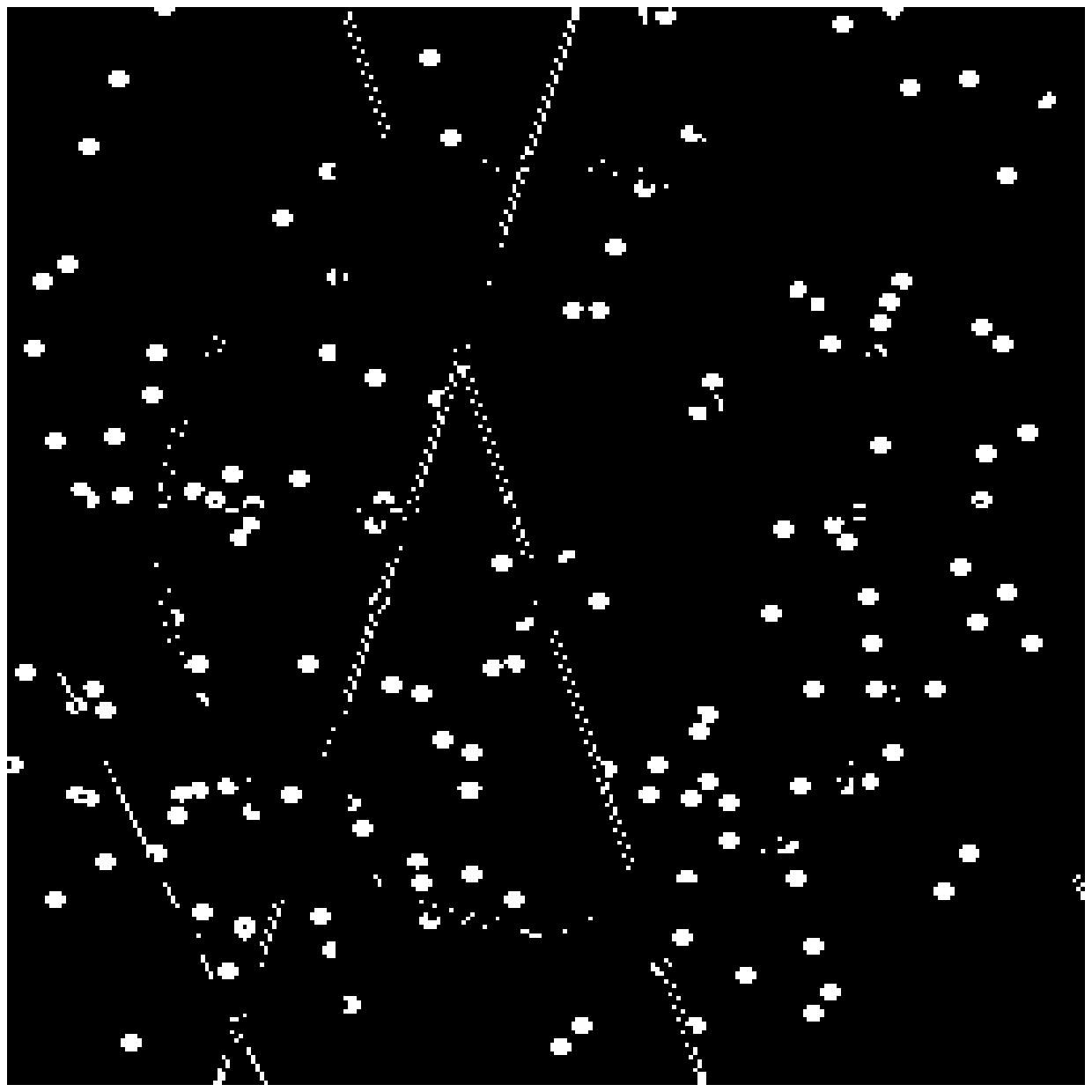}
\includegraphics[width=1.8in]{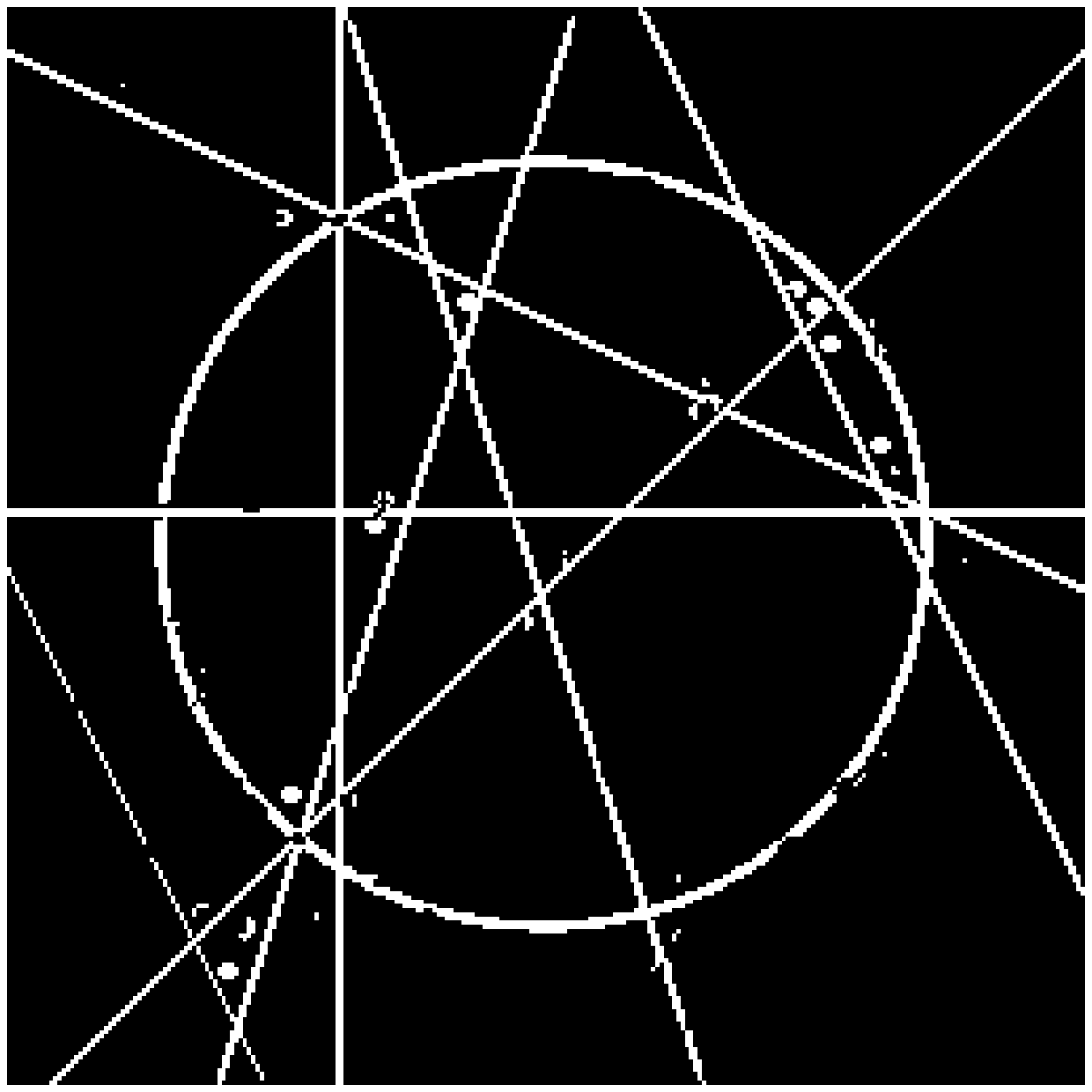}
\put(-395,-15){$SL2D_2$ (curves, $Q_{opt} = 0.16$)}
\put(-260,-15){$SL2D_2$ (points, $Q_{opt} = 0.44$)}
\put(-130,-15){$NSST$ (curves, $Q_{opt} = 0.20$)}\\
\vspace{20pt}
\includegraphics[width=1.8in]{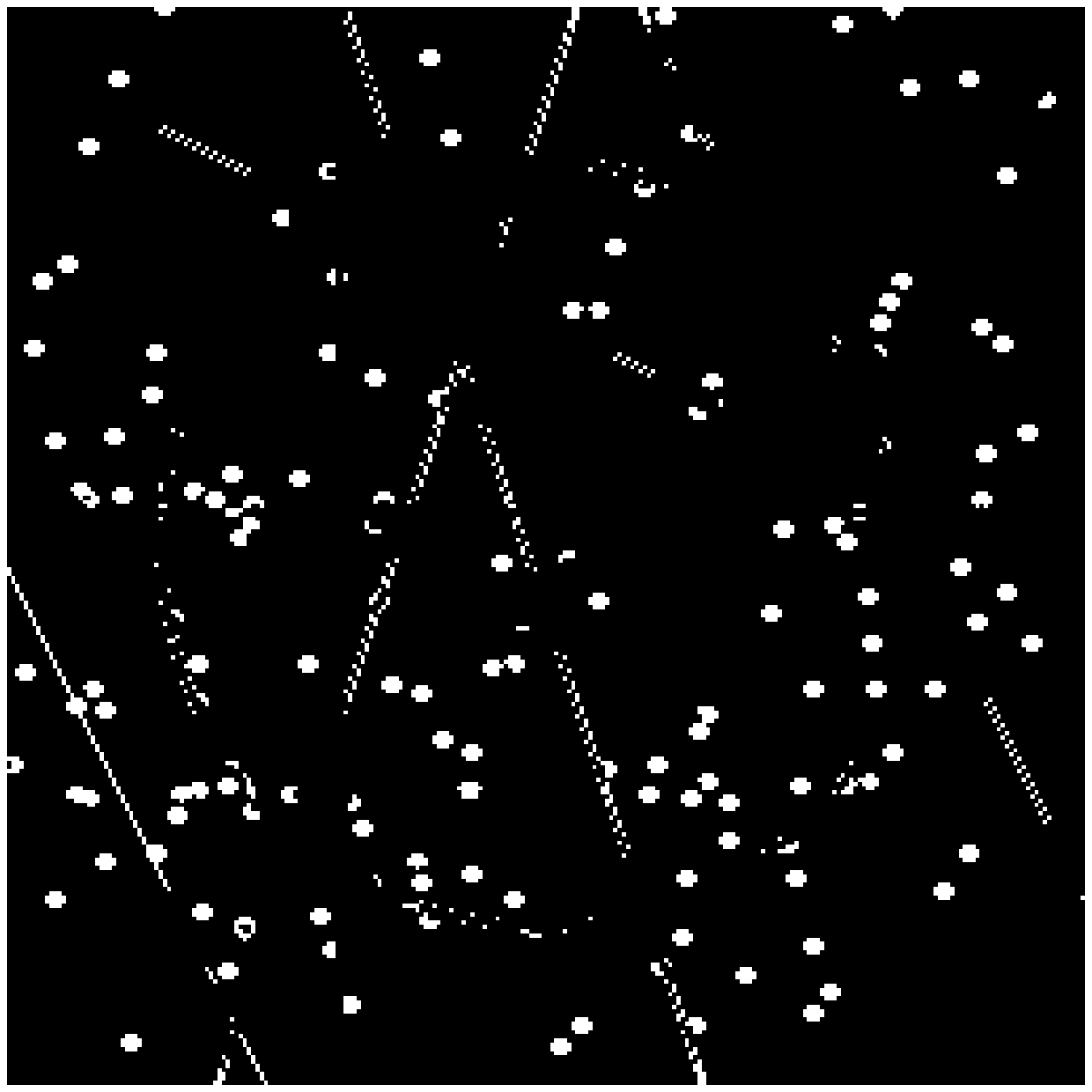}
\includegraphics[width=1.8in]{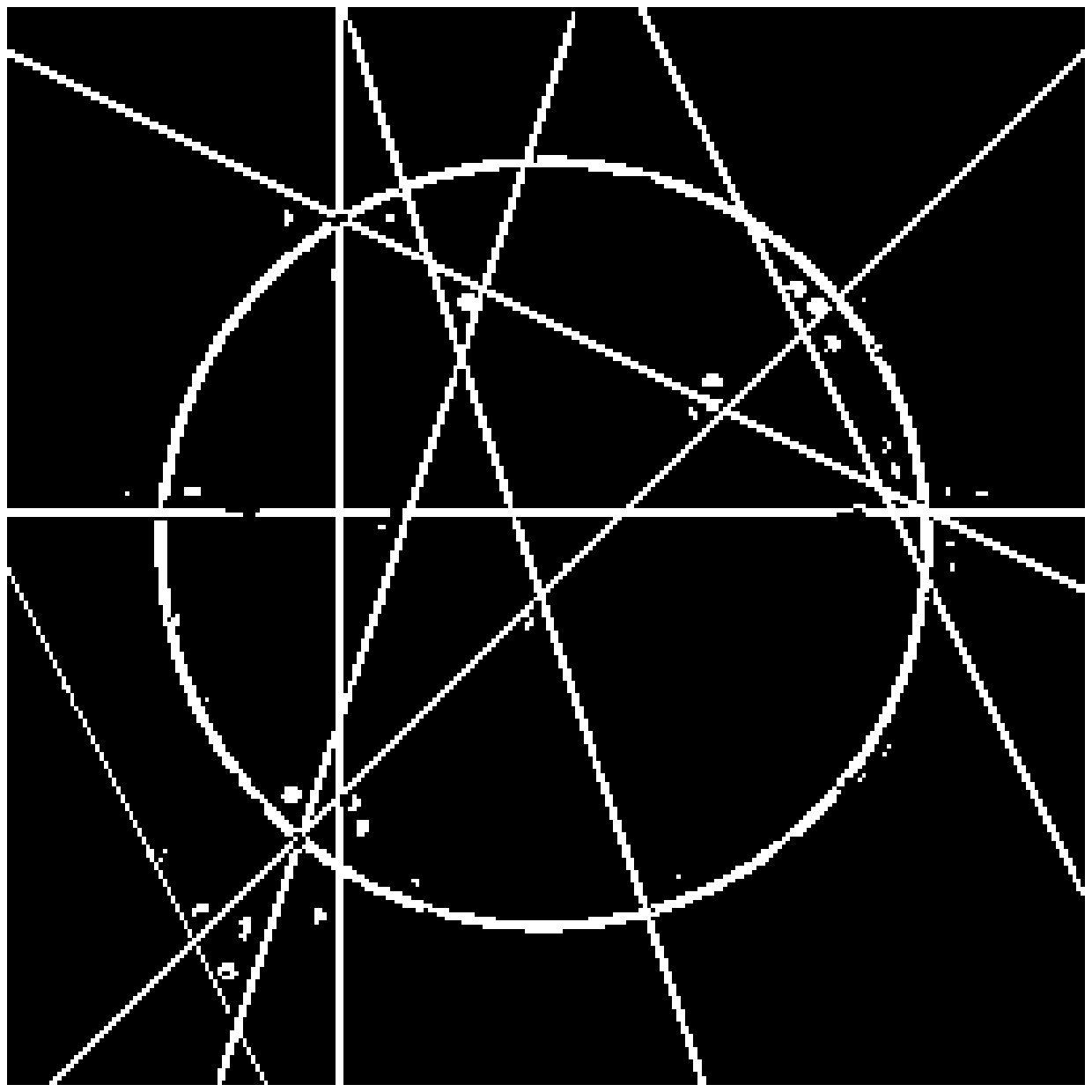}
\includegraphics[width=1.8in]{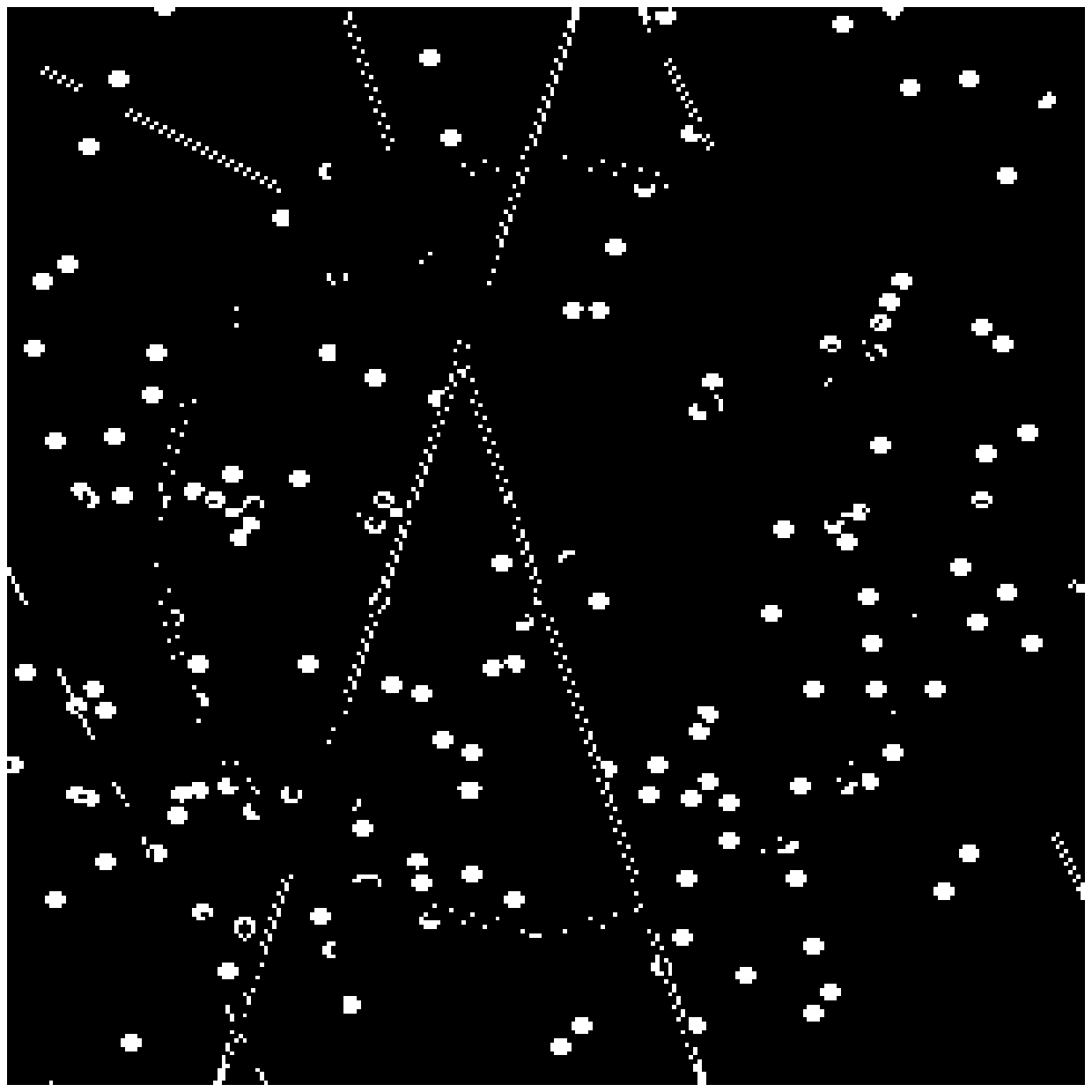}
\put(-395,-15){$NSST$ (points, $Q_{opt} = 0.45$)}
\put(-263,-15){$FDCT$ (curves, $Q_{opt} = 0.20$)}
\put(-130,-15){$FDCT$ (points, $Q_{opt} = 0.47$)}
\end{center}
\caption{The optimal decompositions of a binary image obtained from applying algorithm \ref{alg:image_decomposition} using four different directional transforms.}
\label{fig:results_decomposition_images}
\end{figure}

\begin{figure}[t!]
\begin{center}
\includegraphics[width=1.8in]{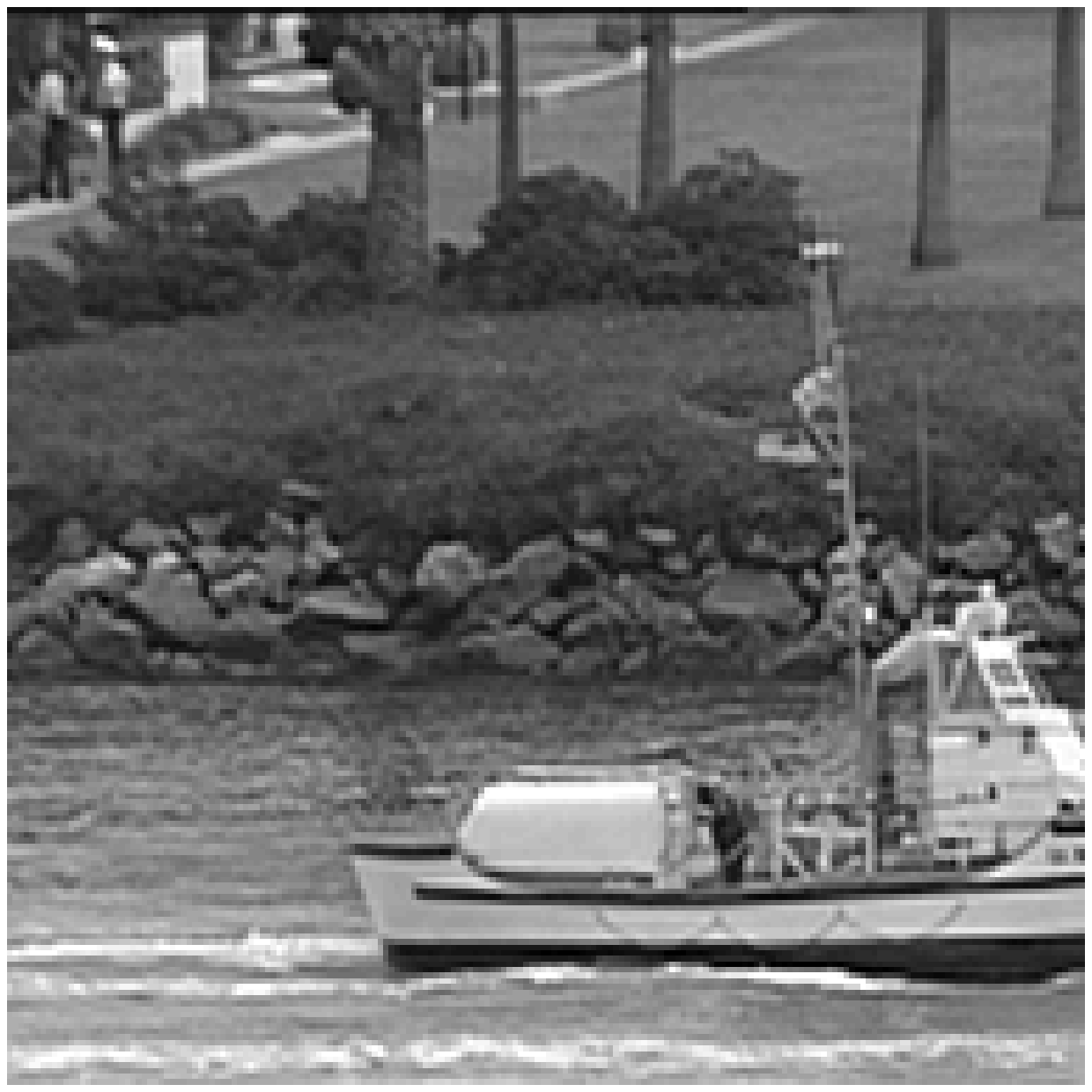}
\includegraphics[width=1.8in]{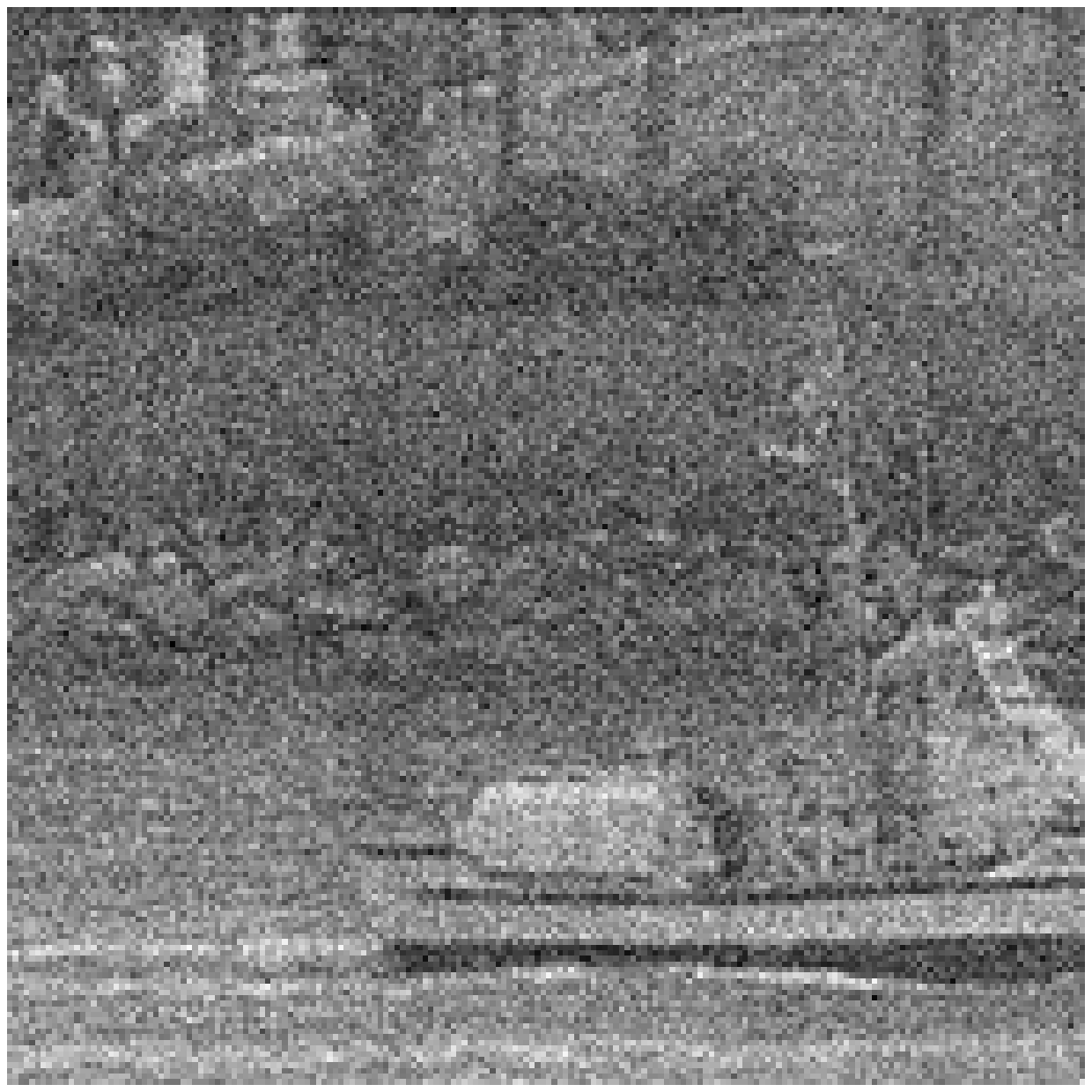}
\includegraphics[width=1.8in]{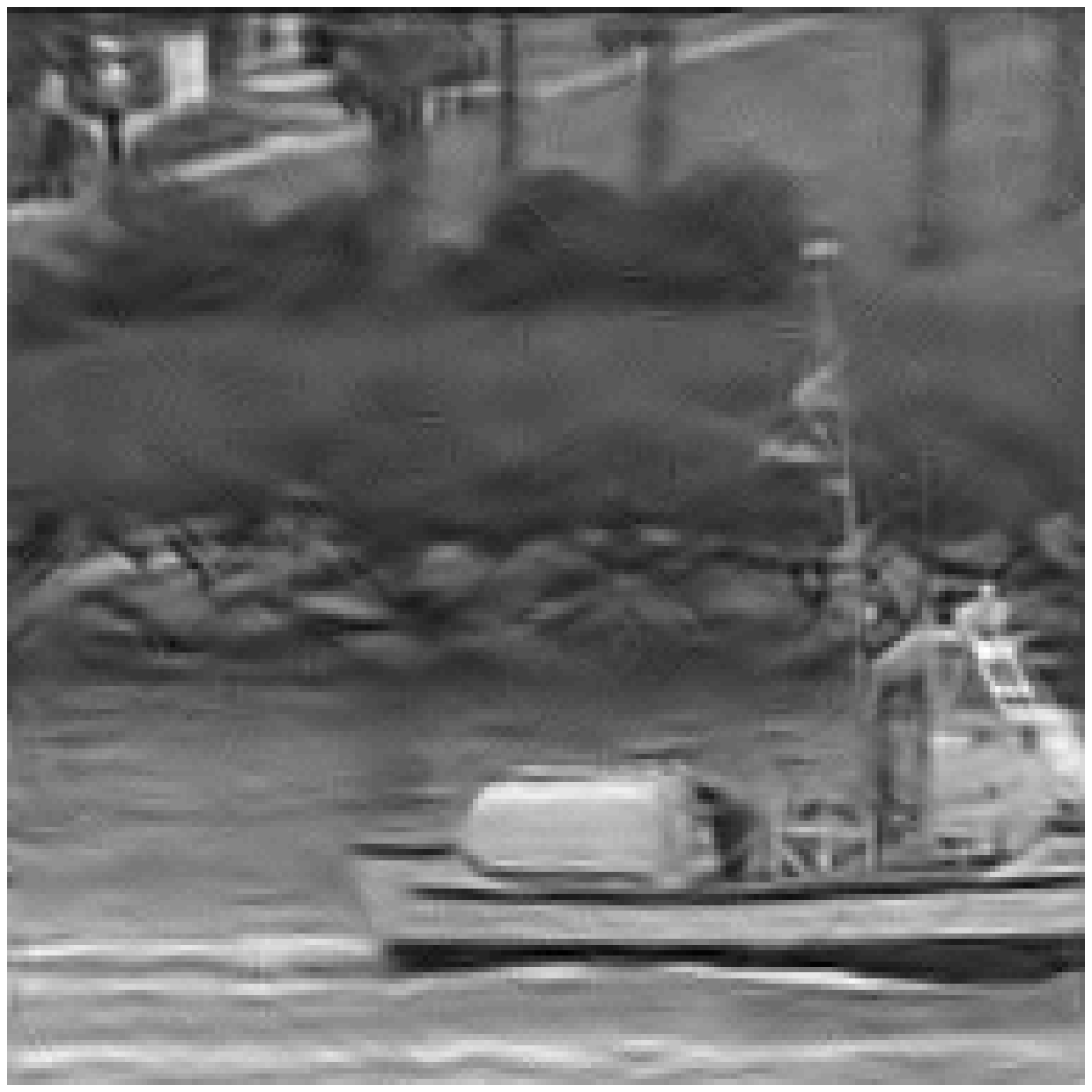}
\put(-350,-15){original}
\put(-260,-15){noisy ($\sigma = 40$ PSNR = $16.06$)}
\put(-115,-15){$SL3D_1$ ($\text{PSNR} = 26.17$)}
\\
\vspace{20pt}
\includegraphics[width=1.8in]{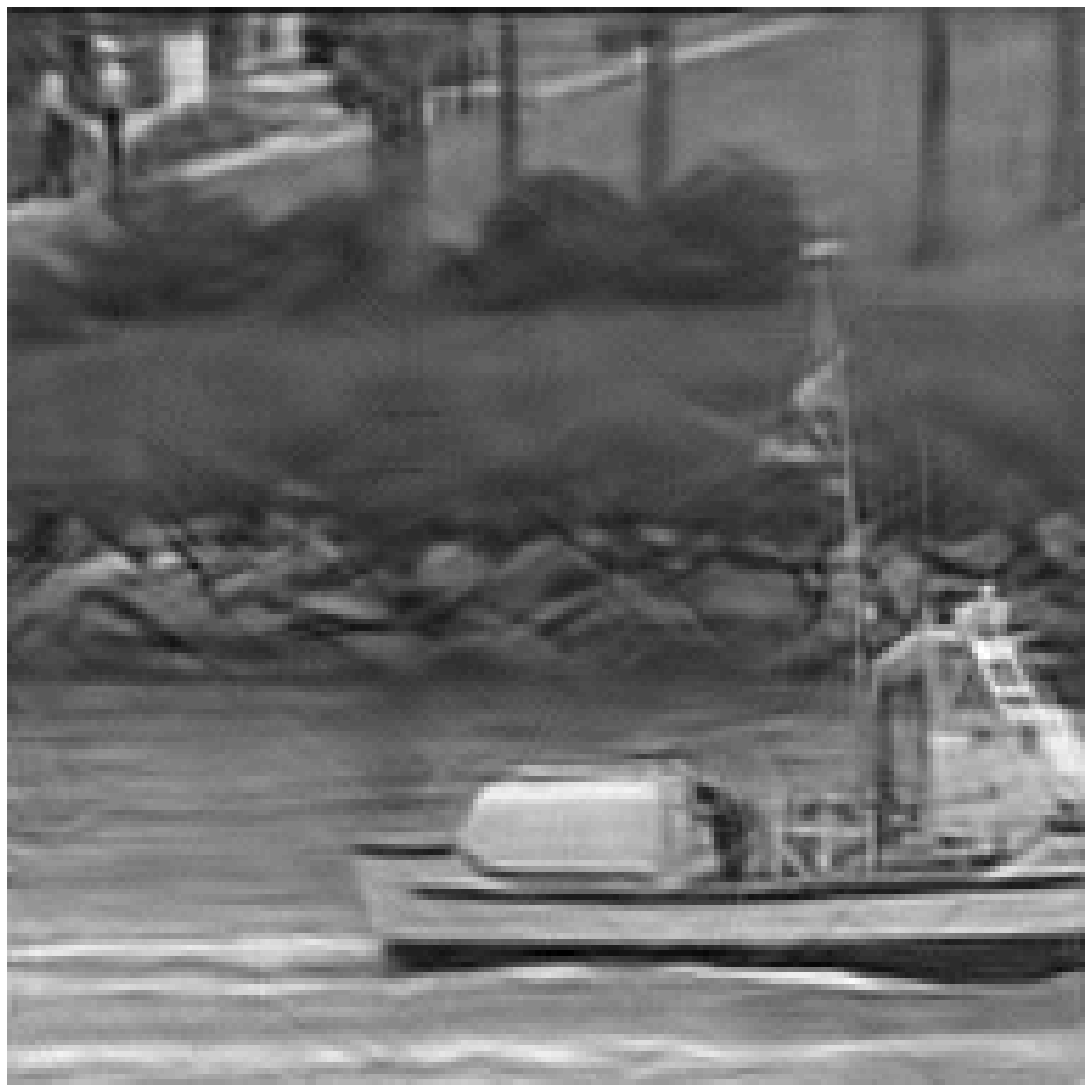}
\includegraphics[width=1.8in]{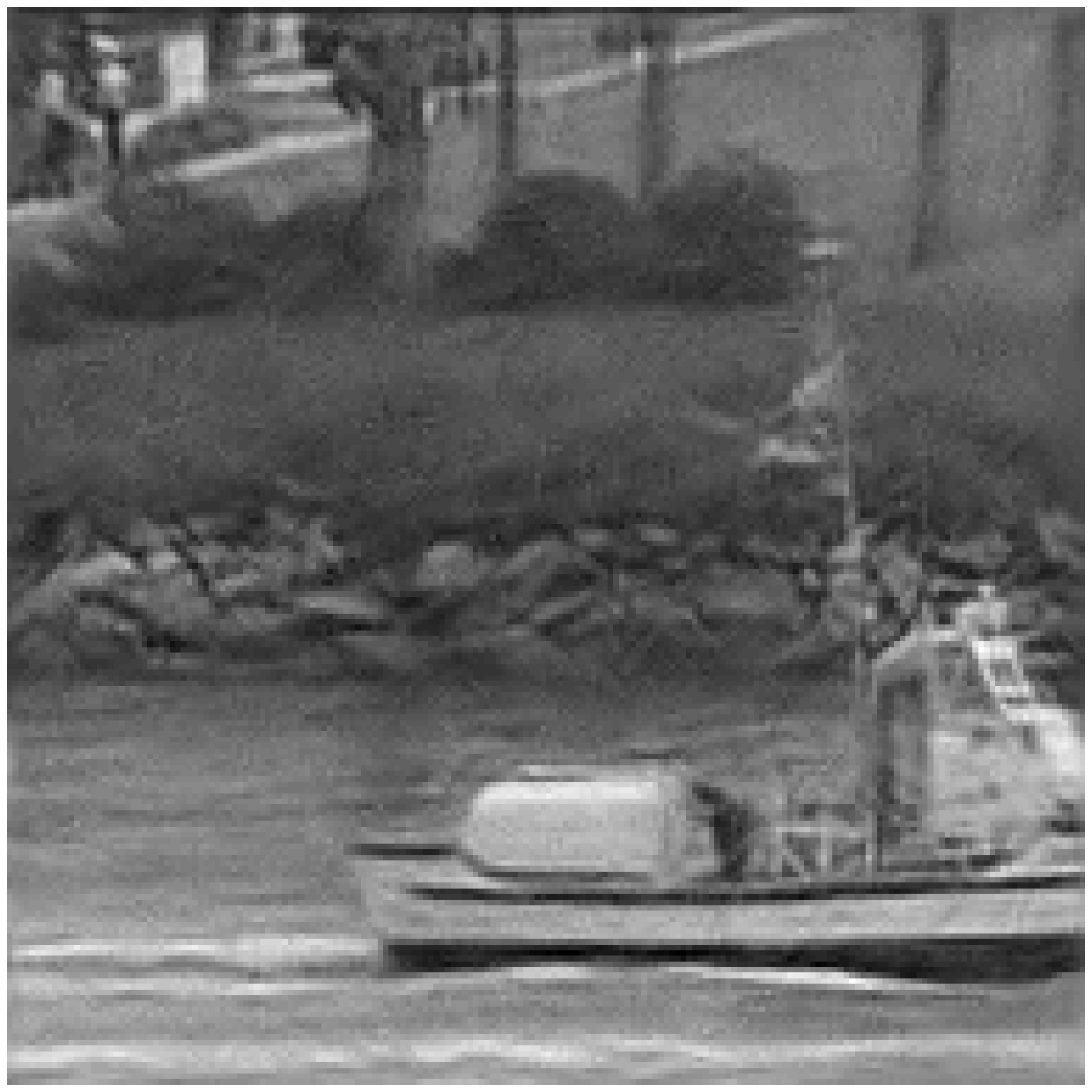}
\includegraphics[width=1.8in]{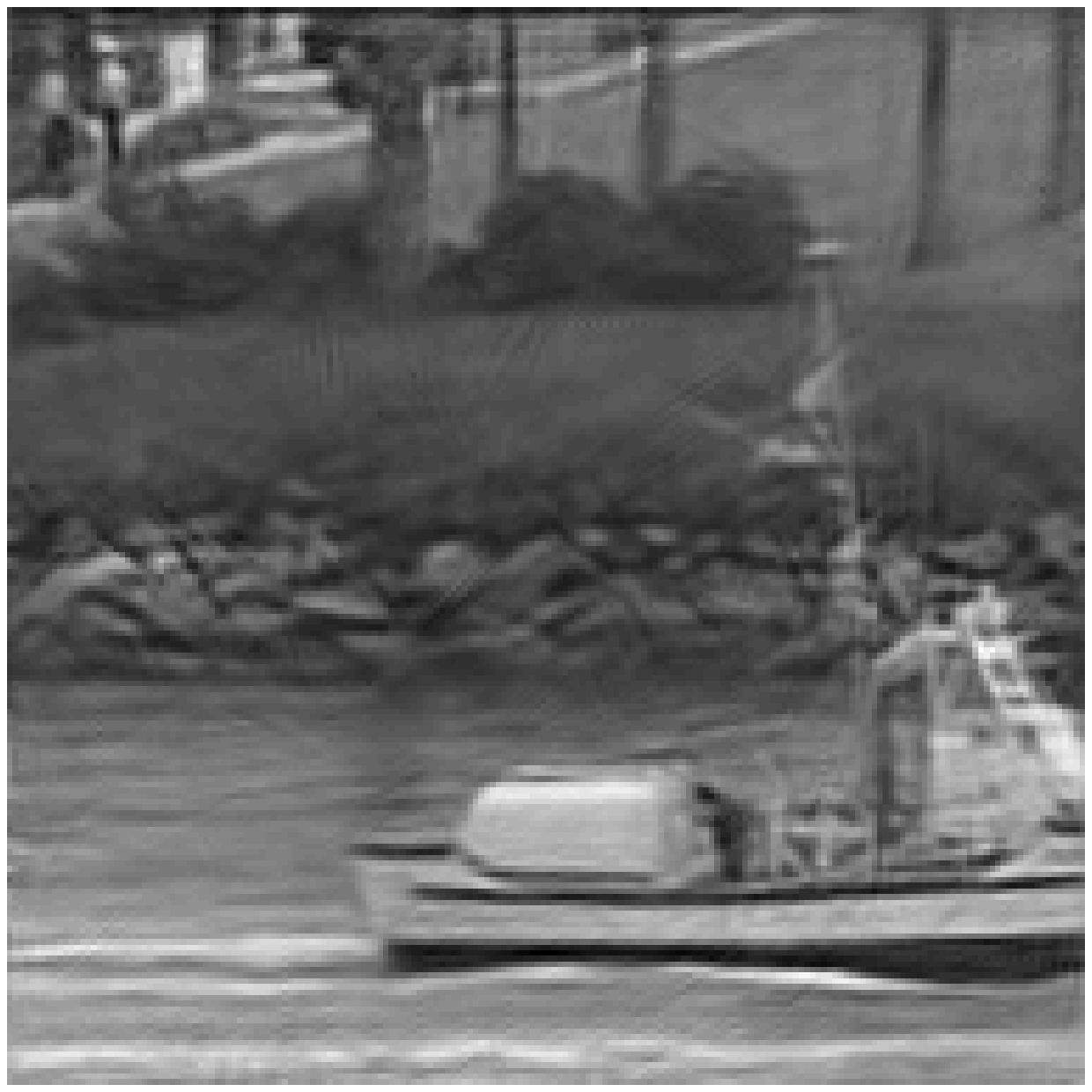}
\put(-385,-15){$SL3D_2$ ($\text{PSNR} = 27.14$)}
\put(-247,-15){$NSST$ ($\text{PSNR} = 25.68$)}
\put(-115,-15){$SURF$ ($\text{PSNR} = 25.91$)}\\
\vspace{20pt}
\includegraphics[width=1.8in]{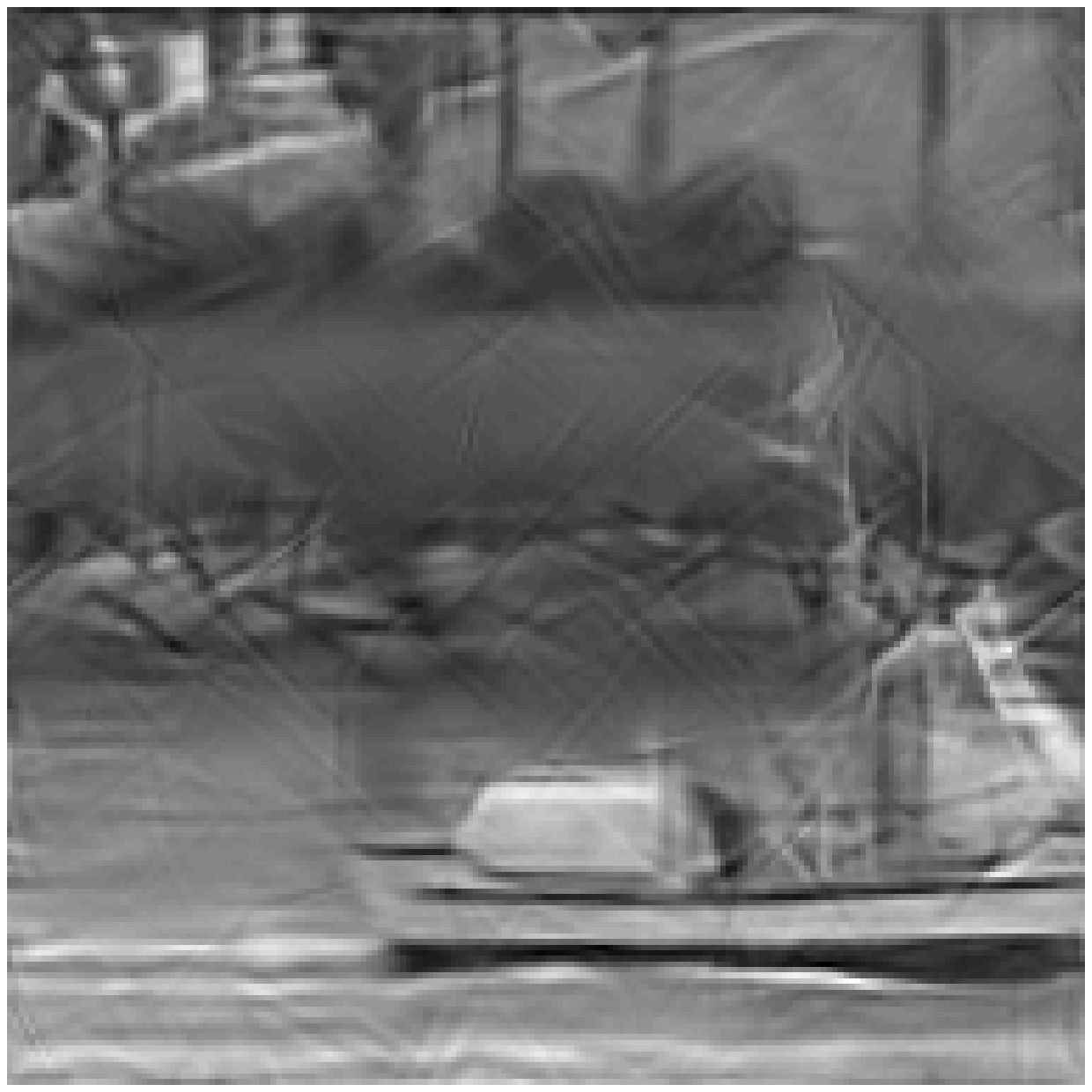}
\put(-115,-15){$SL2D_2$ ($\text{PSNR} = 24.52$)}
\end{center}
\caption{The original, noisy and denoised frame 110 of the coastguard sequence..}
\label{fig:results_video_denoising_images}
\end{figure}

 \begin{figure}[t!]
\begin{center}
\includegraphics[width=1.8in]{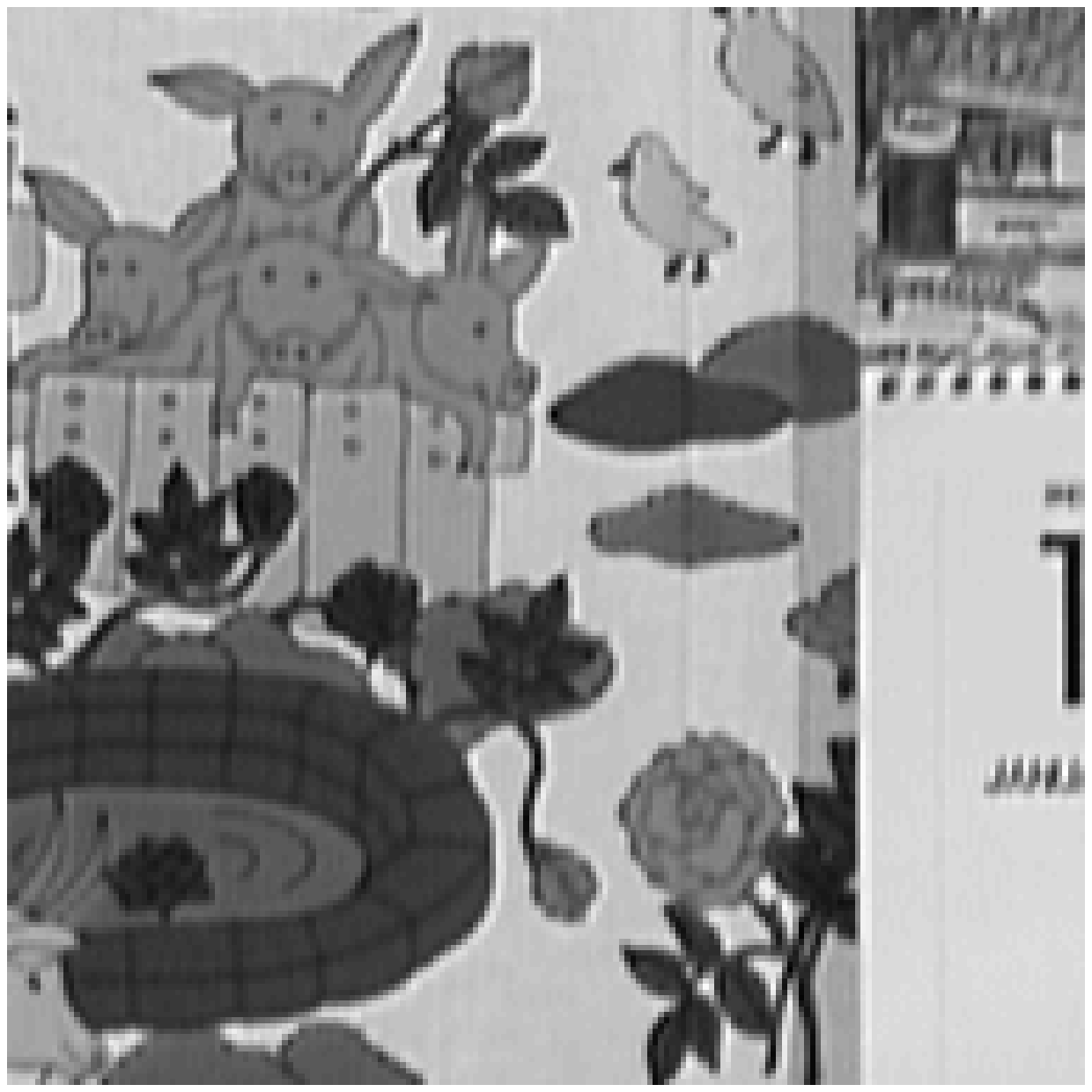}
\includegraphics[width=1.8in]{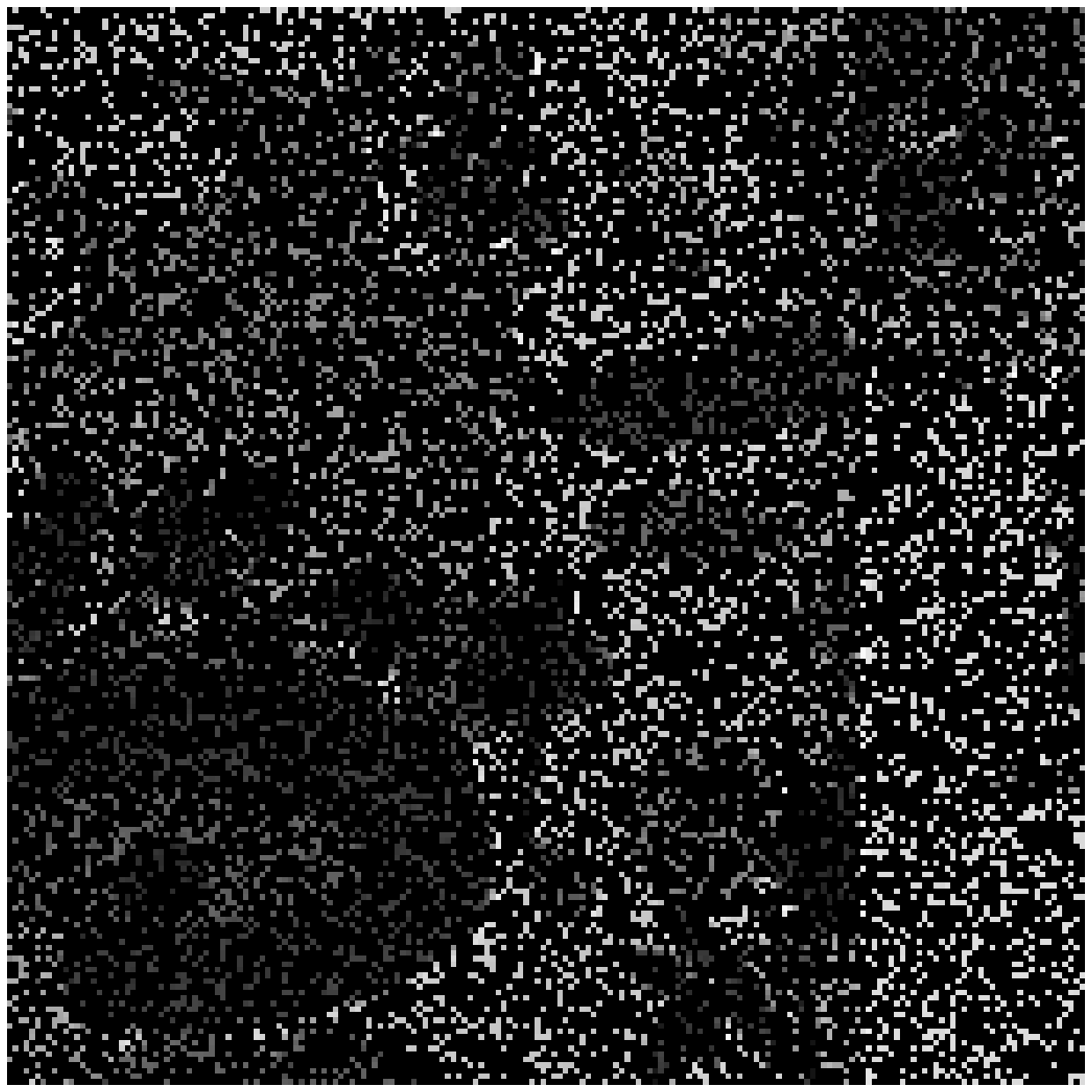}
\includegraphics[width=1.8in]{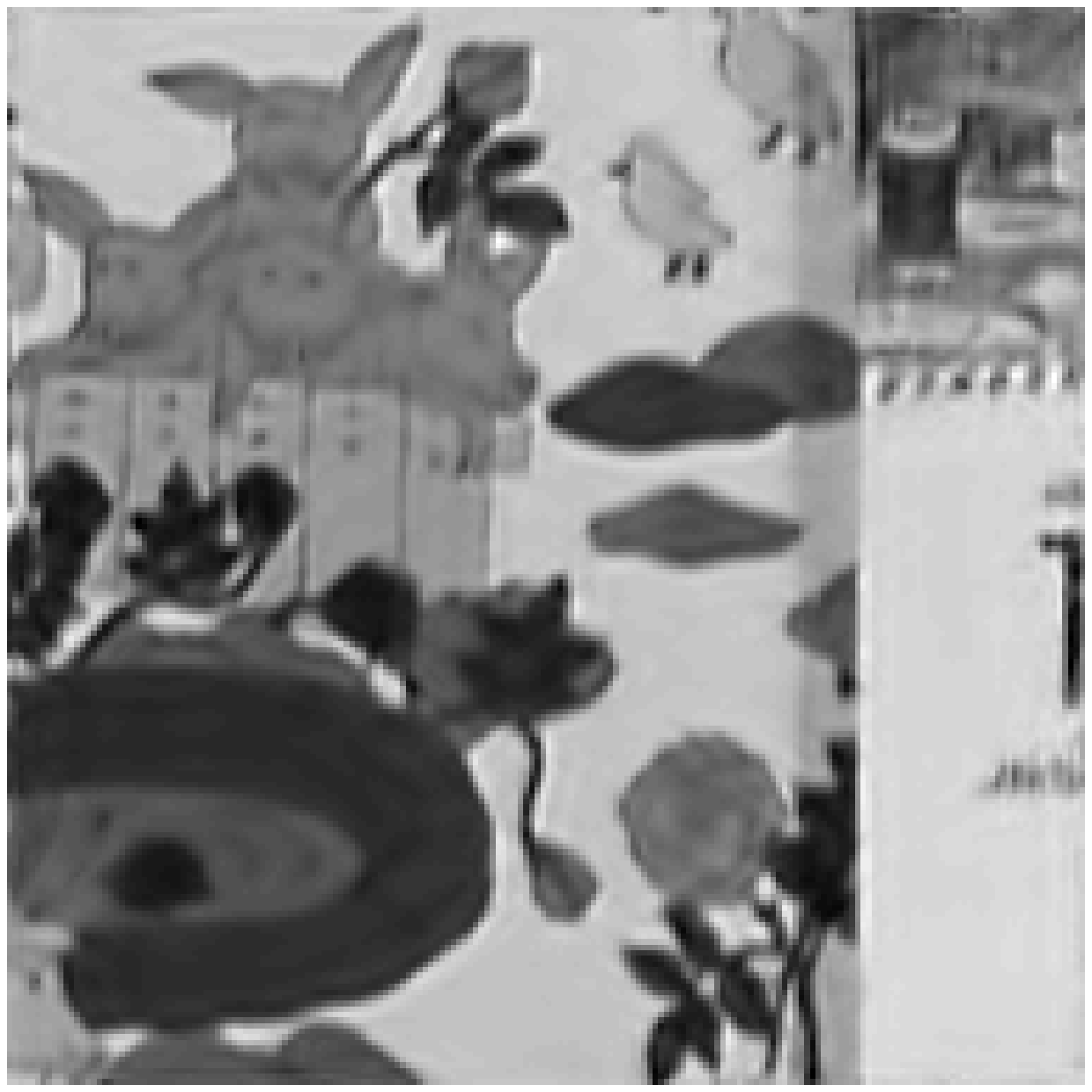}
\put(-350,-15){original}
\put(-260,-15){occluded (\unit{80}{\%} missing)}
\put(-115,-15){$SL3D_1$ ($PSNR = 28.15$)}
\\
\vspace{20pt}
\includegraphics[width=1.8in]{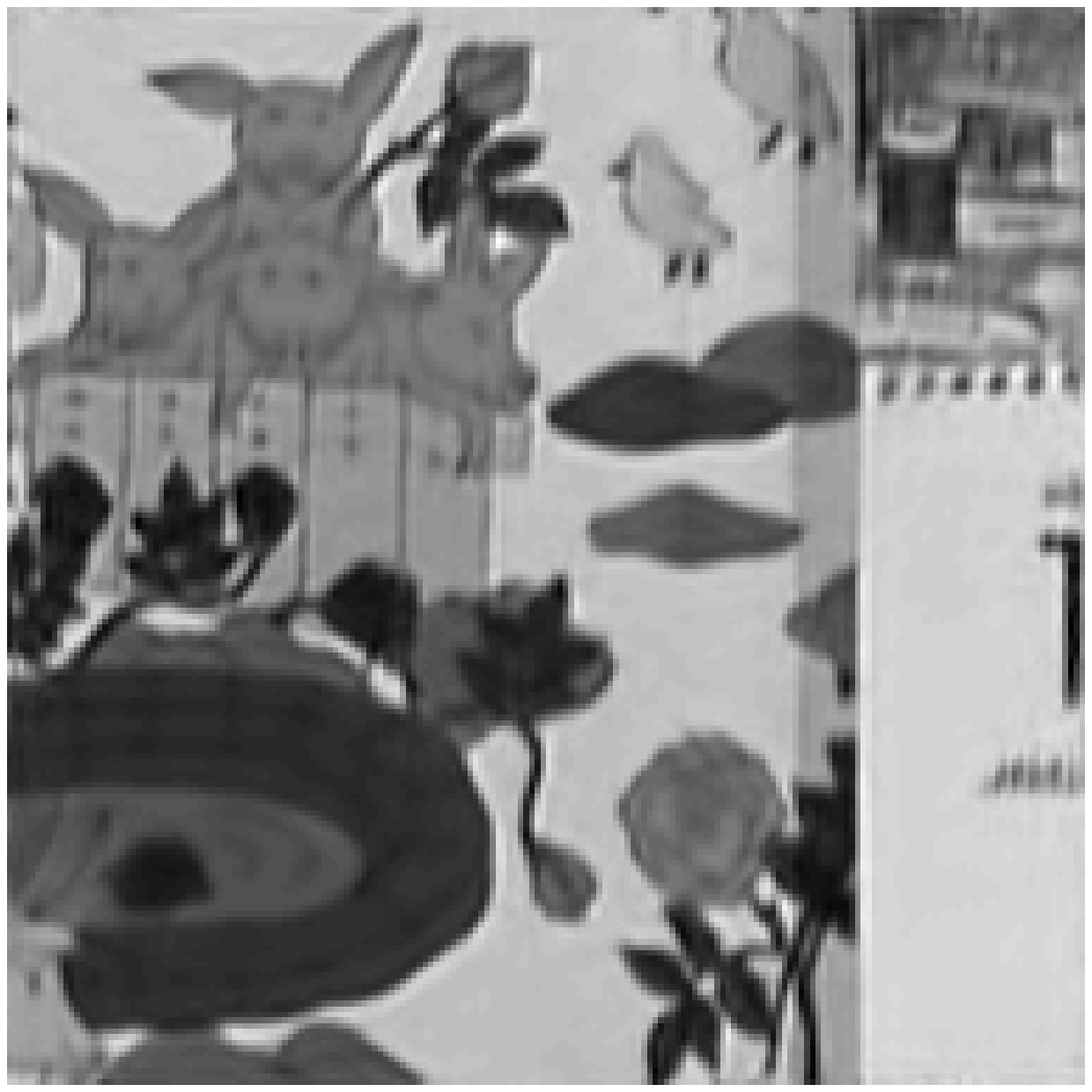}
\includegraphics[width=1.8in]{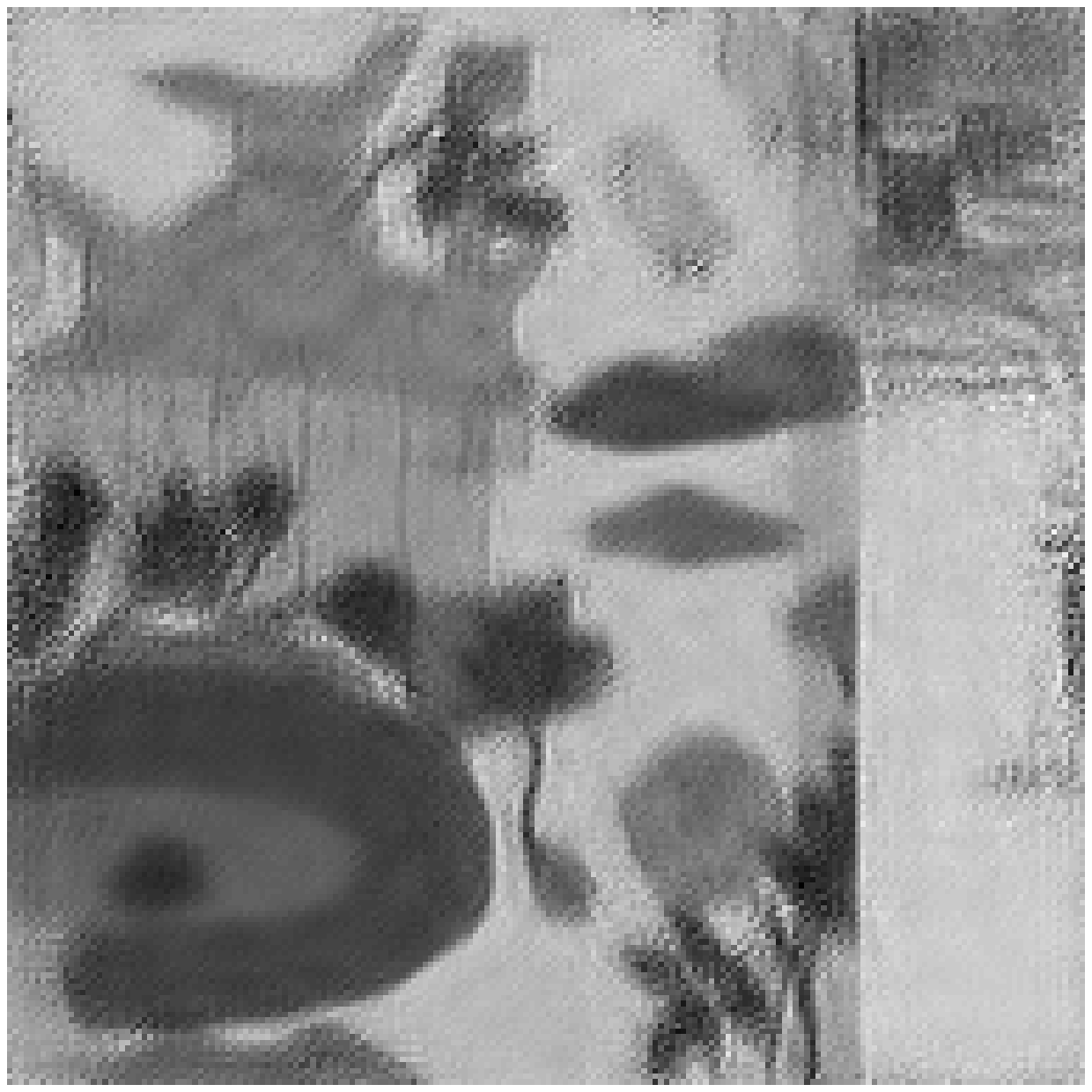}
\includegraphics[width=1.8in]{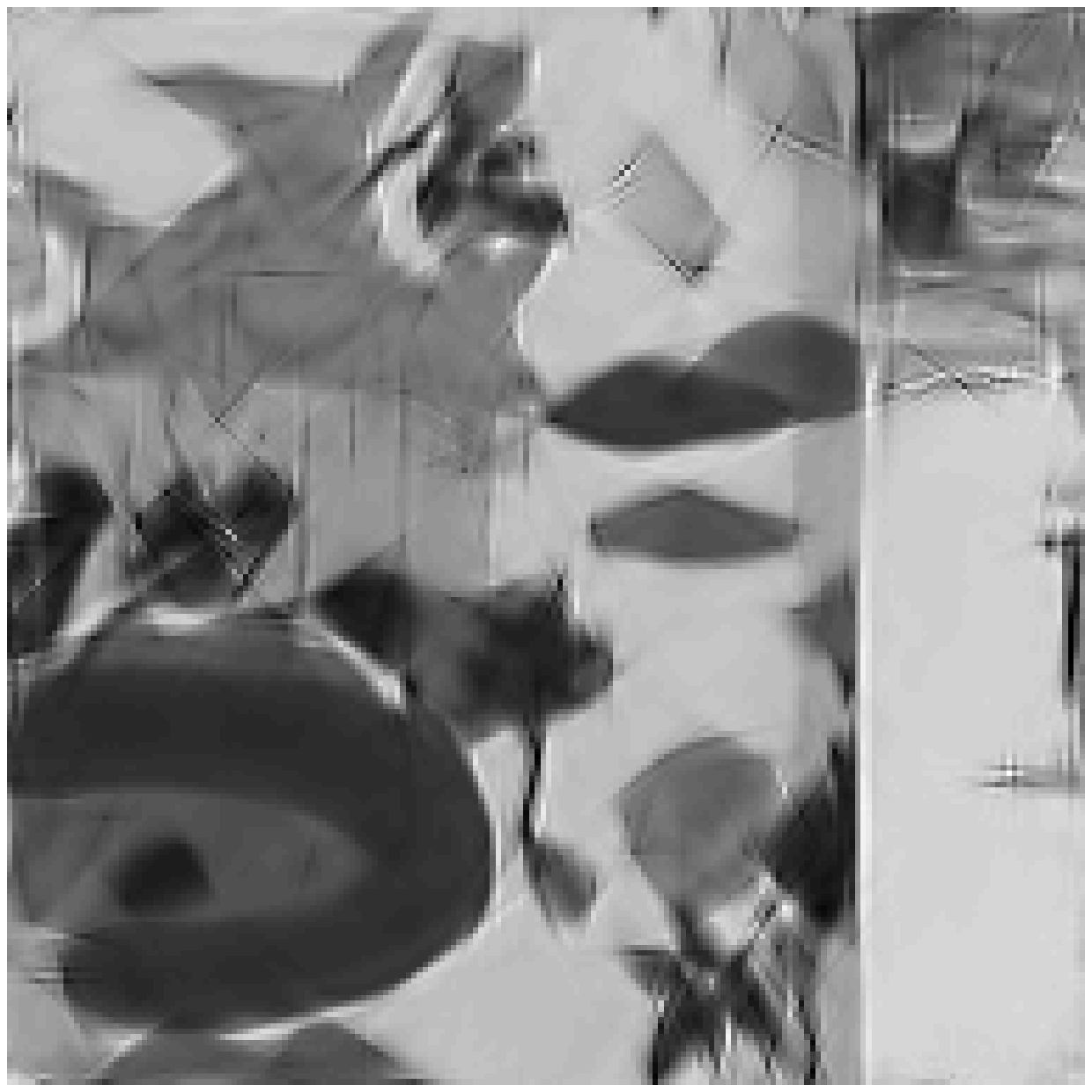}
\put(-385,-15){$SL3D_2$ ($PSNR = 29.98$)}
\put(-247,-15){$SURF$ ($PSNR = 22.40$)}
\put(-115,-15){$SL2D_2$ ($PSNR = 23.63$)}\\
\end{center}
\caption{The original, occluded and inpainted frame 37 of the mobile sequence.}
\label{fig:results_video_inpainting_images}
\end{figure}

\end{appendix}

\end{document}